\documentclass[preprint,review,10pt]{elsarticle}

\usepackage{amsmath,amssymb,amsthm}  % 数学基础
\usepackage{bm,mathrsfs,cases}       % 数学增强
\usepackage{graphicx,epstopdf}       % 图形支持
\usepackage[percent]{overpic}         % 图片标注
\usepackage{subcaption}               % 子图排版
\usepackage{tikz,pgfplots}            % 矢量绘图
\usepackage{float}                    % 定位选项
\usetikzlibrary{decorations.pathreplacing}
\usepackage{comment} 
\usepackage[colorlinks,linkcolor=blue]{hyperref}

\makeatletter
\renewcommand{\maketag@@@}[1]{\hbox{\m@th\normalsize\normalfont#1}}%
\makeatother

\usepackage{multirow,booktabs,makecell,caption}

\usepackage{algorithm}
\usepackage{algorithmic}

\usepackage{geometry}
\usepackage{xcolor}
\usepackage{listings}
\newtheorem{remark}{Remark}[section]

\theoremstyle{definition}
\newtheorem{example}{Example}[section]

\makeatletter
\def\ps@pprintTitle{%
 \let\@oddhead\@empty
 \let\@evenhead\@empty
 \def\@oddfoot{}%
 \let\@evenfoot\@oddfoot}
\makeatother

\makeatletter\@addtoreset{equation}{section} \makeatother

\begin{document}
\begin{frontmatter}

\title{Energy-Equidistributed Moving Sampling Physics-informed Neural Networks for Solving Conservative Partial Differential Equations}

% ===== 作者信息 =====
\author[label1,label2]{Qinjiao Gao} 
 \ead{09110180026@fudan.edu.cn}
\author[label1]{Longzhe Xu}
\ead{23020040160@pop.zjgsu.edu.cn}
\author[label1]{Dongjiang Wang}
\ead{23020040065@pop.zjgsu.edu.cn}
\author[label3]{Ran Zhang\corref{cor1}}
\ead{zhang.ran@mail.shufe.edu.cn}

% ===== 单位地址 =====
\address[label1]{School of Statistics and Mathematics, Zhejiang Gongshang University, Hangzhou, China}
\address[label2]{Collaborative Innovation Center of Statistical Data Engineering, Technology $\&$ Application, Zhejiang Gongshang University, Hangzhou, China}
\address[label3]{School of Mathematics, Shanghai University of Finance and Economics, Shanghai, China}

% ===== 通讯作者与基金声明 =====
 \cortext[cor1]{Corresponding author.\\
 This work is supported by the National Natural Science Foundation of China (Grants No. 12001487 and No. 11901377),  and the Characteristic  $\&$  Preponderant Discipline of Key Construction Universities in Zhejiang Province (Zhejiang Gongshang University- Statistics).}

% ===== 摘要 =====
 \begin{abstract} 
This paper presents a novel Energy-Equidistributed adaptive sampling framework for multi-dimensional conservative PDEs, introducing both location-based and velocity-based formulations of Energy-Equidistributed moving mesh PDEs (EMMPDEs). The framework utilizes the energy density function as the monitor function, ensuring that mesh adaptation dynamically tracks energy evolution during temporal integration. These theoretical developments are integrated with deep neural networks to establish the Energy-Equidistributed Moving Sampling Physics-Informed Neural Networks (EEMS-PINNs), which integrate physics-informed learning with energy-adaptive mesh optimization.  Extensive numerical experiments demonstrate that EEMS-PINNs effectively maintain solution accuracy in long-time simulations while preserving conserved energy. The framework's robustness is further evidenced by its stable performance in non-conservative systems. The code for this paper can be found at \href{https://github.com/sufe-Ran-Zhang/EMMPDE}{https://github.com/sufe-Ran-Zhang/EMMPDE}.

%This paper presents an innovative energy-equidistributed adaptive sampling strategy for simulating conservative partial differential equations (PDEs). a location-based and  a modified velocity-based energy-equidistributed moving mesh PDEs (EMMPDEs) are introduced. These theoretical advances are integrated with deep neural networks to create Energy-Equidistributed Moving Sampling Physics-informed neural networks (EEMS-PINNs), a novel framework that combines physics-informed learning with energy-adaptive mesh optimization. Comprehensive numerical experiments demonstrate EEMS-PINNs' advantages: maintaining solution precision during long-time simulations while preserving the system's conserved energy. The robustness of our approach is further validated through its stable performance in systems with non-conservative total energy dynamics. 

%This paper presents an innovative energy-equidistributed adaptive sampling strategy for simulating Hamiltonian partial differential equations (PDEs). A novel framework, named the Energy Equidistribution Principles (EEPs) and three distinct groups of energy-equidistributed moving mesh PDEs (EMMPDEs) are introduced, highlighting the critical role of energy conservation in achieving accurate simulations. Building on EEPs and Physics-informed Neural Networks(PINNs), come up with the innovative combination Energy-equidistributed moving sampling PINNs(EEMS-PINNs). Finally, four typical numerical examples are conducted to demonstrate the effectiveness and conservation character of the proposed method.

\end{abstract}

% ===== 关键词 =====
\end{frontmatter}

\label{Maintext}
\section{Introduction}
Physics-informed machine learning combines physics-based models with data-driven techniques to solve partial differential equations (PDEs) more effectively recently \cite{han2017deep, sirignano2018dgm, chiu2022can}. By embedding physical constraints into machine learning frameworks, it achieves better accuracy and efficiency than purely data-driven approaches, particularly for high-dimensional, inverse, or uncertain problems \cite{han2017deep, EIDNES2024112738, Boltzmann}. In this paper, we focus on conservative PDEs that possess a first integral or an invariant energy with respect to time. 
%Conservative PDEs play a vital role in modeling various physical phenomena in science and engineering, such as elasticity equations in elastodynamics, the shallow-water equations in climate modeling, and the Kuramoto–Sivashinsky equation in chemical reaction dynamics. 
%Thus looking for numerical schemes with conservation properties are essential in solution approximation of the conservative PDEs. 
%These equations are characterized by their conservative nature and underlying symplectic structure, %which are fundamental to their discretization and numerical treatment. 

The development of numerical schemes preserving conservation properties dates back to pioneering work in the 1920s, when Courant, Friedrichs and Lewy first established the conservation-preserving finite difference approximation for wave equations \cite{Lewy1928}. 
These conservation properties have proven fundamental to demonstrating solution existence, uniqueness, and global stability in numerical analysis \cite{feng2010symplectic}. Subsequent decades witnessed significant theoretical and computational advancements through various conservative numerical approaches, including: finite difference methods, spectral methods, geometric integrators, variational integrators, wavelet collocation methods, Hamiltonian boundary value methods and
meshless methods \cite{Hairer2006, Runge-Kutta, wavelet, sun2024energy}.

In recent years, deep neural networks have emerged as powerful tools for solving PDEs, where the training process fundamentally relies on optimizing a physics-informed loss function that ensures the network's output satisfies the governing equations \cite{raissi2019physics}. 
This loss function is computed at residual collocation points distributed throughout the computational domain, analogous to the grid points used in conventional numerical discretizations. The spatial arrangement and sampling strategy of these residual points significantly influence the solver's approximation accuracy and convergence behavior \cite{lu2021deepxde}. Common sampling approaches include uniform grid distributions, completely random sampling, and quasi-random techniques such as Latin hypercube sampling \cite{McKay200055}, each offering distinct advantages in terms of point coverage and computational efficiency.

%Primarily the deep network solvers are constructed by optimizing PDE loss function, which ensures that the trained network aligns with the PDE being solved. This loss is evaluated at a set of scattered residual collocation points, similar to the grid points used in traditional numerical methods. Therefore, the location and distribution of these residual points are crucial for the approximation performance. The most common sampling methods include equispaced uniform point sampling,  uniformly random sampling, and Latin hypercube sampling etc.

To improve training efficacy, the collocation points can be dynamically updated during optimization iterations rather than remaining fixed. 
The Residual-based Adaptive Refinement (RAR) method \cite{lu2021deepxde} progressively adds high-residual points to the training set, while Residual-based Adaptive Refinement with Distribution (RAD) \cite{wu2023comprehensive} promotes balanced resamples across the whole computational domain. DAS-PINN method \cite{TANG2023111868} utilizes a residual-based generative model to generate new collocation points to refine the training set. 
\cite{nabian2021efficient} studied the performance of an importance sampling approach for efficient training of PINN.   While these adaptive methods effectively minimize PDE residuals by focusing on high residual error regions, they may inadvertently introduce two limitations: increased solution approximation errors in other parts of the domain due to uneven error distribution, and elevated computational costs from the continual point augmentation.

Alternatively, \cite{hou2023enhancing} introduces WAM-PINN which relies on monitor functions as an indicator of the variation of the solution function to guide collocation point adjustments. 
This approach bears conceptual similarity to the moving mesh (i.e. $r$-adaptation) method introduced by Huang et al. \cite{huang1994moving}, where mesh nodes are dynamically repositioned while preserving topological connectivity. The movement of nodes is governed by a moving mesh PDE (MMPDE) derived from variational equidistribution principles \cite{de1973good}, with a user-defined monitor function serving as the crucial link between the underlying PDE solution and mesh adaptation dynamics. In \cite{tang}, Tang et al. also proposed a conservative-interpolation formula to keep the mass-conservation
of the underlying numerical PDE solution at each redistribution step.  

Recent advances have integrated these concepts with deep learning frameworks. \cite{M2N} designed a
Neural-Spline based model and a graph attention network (GAT) based
model for mesh deformation.  
The Data-free Mesh Mover (DMM) introduced in \cite{DMM} utilizes a two-branch neural architecture with learnable interpolation that preserves solution information during mesh deformation. 
Flow2Mesh \cite{Flow2Mesh} combines a perceptual feature network (PFN) with a graph-based mesh movement network (MMN) for fluid dynamics simulation. 
MMPDE-Net \cite{MMPDE_NET} investigates an adaptive sampling neural network based on MMPDE principles. Building on this, the Moving Sampling PINN (MS-PINN) framework employs an initial PINN's solution to inform the monitor function construction in MMPDE-Net, which then generates optimized sampling points to accelerate PINN's loss convergence.

The selection of monitor functions has been a fundamental yet challenging aspect of moving mesh methods, as it directly governs mesh concentration and critically impacts solution accuracy \cite{Monge}. While conventional approaches often rely on ad hoc choices, our recent work \cite{gao2025energy} introduces Energy Equidistribution Principles (EEPs) for one-dimensional Hamiltonian PDEs. This framework utilizes the energy density function as the monitor function, ensuring that mesh adaptation dynamically tracks energy evolution during temporal integration. Building on the Energy Equidistribution Principles (EEPs), we develop three distinct Energy-Equidistributed moving mesh PDEs (EMMPDEs) derived through root-searching iteration, Lagrangian variational formulation, and Taylor expansion approximation, respectively. The numerical implementation employs an alternating solution strategy, where the governing PDE and EMMPDE are discretized via finite difference and solved iteratively—a mesh at the 
new time level $t^{n+1}$ is first generated using the mesh and the physical solution
at the current time level $t^{n}$, and the underlying PDE solution is then obtained at the new time level $t^{n+1}$. 

This paper extends our previous Energy-Equidistribution Principles \cite{gao2025energy} to multi-dimensional conservative PDEs through two key contributions: a location-based \cite{hou} EMMPDE derived from variational principles in computational space, and a modified velocity-based \cite{gao2018adaptive,gao2019moving} EMMPDE formulation. These theoretical advances are  then integrated with deep neural networks to create Energy-Equidistributed Moving Sampling PINNs (EEMS-PINNs), a novel framework that combines physics-informed learning with energy-adaptive mesh optimization. Comprehensive numerical comparisons with PINNs and WAM-PINNs demonstrate EEMS-PINNs' advantages: maintaining solution precision during long-time simulations while  preserving the system's conserved energy. 
The robustness of our approach is further validated through its stable performance in systems with non-conservative total energy dynamics. 

The structure of the paper is as follows: In  the following Section \ref{sec:intro}, we outlines the fundamental concept of PINNs and conservative PDEs. Section \ref{sec:EEMPE} is the central part of the paper. We present the Energy-Equidistributed moving mesh strategies and the algorithmic framework of EEMS-PINNs. Numerical experiments are implemented in Section \ref{sec:num} to validate the feasibility of the methodology. Finally, the last section  summarizes the key conclusions and address the the future researches.

 \section{Preliminary}\label{sec:intro}
 
 \subsection{Physics-informed neural networks (PINNs)}
 
 %Physics-Informed Neural Networks (PINNs) represent a novel paradigm in scientific computing of PDE that seamlessly integrates physical laws with deep learning architectures. This approach leverages neural networks to solve PDEs while incorporating physical constraints directly into the learning process.

By encoding the governing PDEs directly into the neural network's loss function through automatic differentiation, PINNs enforce physical consistency while maintaining the flexibility and approximation power of modern machine learning architectures. This paradigm shift enables the solution of PDE systems without relying on traditional discretization methods, while preserving fundamental physical principles through mathematical constraints embedded in the optimization process.

 Consider the PDEs of the form
 \begin{equation}\label{PDE}
\begin{aligned}
&\boldsymbol{u}_t + \mathcal{N}[\boldsymbol{u}](\boldsymbol{x},t) =\boldsymbol{f}(\boldsymbol{x},t),   \ \boldsymbol{u}=\boldsymbol{u}(\boldsymbol{x},t)\in \mathbb{R}^m,  \  \boldsymbol{x} \in \Omega\subseteq \mathbb{R}^d,    \  t \in [0,T],  \\
&\boldsymbol{u}(\boldsymbol{x},0) = \boldsymbol{u}_0(\boldsymbol{x}), \quad x \in \Omega, \\
&\boldsymbol{u}\left( \boldsymbol{x},t \right) =\boldsymbol{u}_b\left(\boldsymbol{x},t \right) , \  (\boldsymbol{x},t) \in  \partial \varOmega \times  \left[ 0,T \right],
\end{aligned}
\end{equation}
where
$\boldsymbol{u}_0(\boldsymbol{x})$ is the initial condition,  $\boldsymbol{u}_b(\boldsymbol{x},t)$ represents boundary conditions and $f(\boldsymbol{x},t)$ is the corresponding right hand side, and
$\mathcal{N}[\cdot]$ denotes a nonlinear differential operator up to $p$-order acting on the unknown solution $\boldsymbol{u}(\boldsymbol{x},t)$ respect to the spatial variable $\boldsymbol{x}$.
%i.e., 
%\[
%\boldsymbol{u}^{(\alpha)}=\left\{\frac{\partial^{|\alpha|} \boldsymbol{u}}{\partial x_1^{\alpha_1} x_2^{\alpha_2} \cdots x_d^{\alpha_d}}:|\alpha| \leq p\right\}, \quad \alpha \in\left(\mathbb{Z}^{\geq}\right)^d. \]

Typically, PINNs employ a fully connected neural network $\widehat{\boldsymbol{u}}\left( \boldsymbol{x}, t;\Theta \right)$ to approximate the solution of the function $\boldsymbol{u}(\boldsymbol{x},t)$, where $\Theta$ represents the parameter set of the neural network. A deep neural network (DNN) $\widehat{\boldsymbol{u}}\left( \boldsymbol{x}, t;\Theta \right)$ possessing  $L-1$ hidden layers with widths $l_k$ ($1 \leq  k \leq  L- 1$), and an input and output layer is denoted as follows:	
	\begin{align}\label{eq:network}
		\begin{split}
			&\boldsymbol{z}^{(0)}=(\boldsymbol{x},t),               \\
			&\boldsymbol{z}^{(k)}=\mathscr{L}_k(\boldsymbol{z}^{(k-1)}) =\phi( W^{\left( k \right)}\boldsymbol{z}^{(k-1)}+b^{\left( k \right)}),\ k=1,\cdots,L-1,\\
            & \widehat{\boldsymbol{u}}\left( \boldsymbol{x}, t;\Theta \right):=
            \mathscr{L}_L(\boldsymbol{z}^{(L-1)})= W^{\left( L \right)}\boldsymbol{z}^{(L-1)}+b^{\left( L \right)},
			%&f\left(x,t;\theta \right) =\left( \mathscr{L} _L\cdots 
			%\sigma \mathscr{L} _1 \right) \left( x^0 \right),  
		\end{split}
	\end{align}
where $\mathscr{L}_k: \mathbb{R}^{l_{k-1}} \to \mathbb{R}^{l_{k}}$ denotes the mapping of the $k$-th layer with $l_0=d+1$ and $l_L=m$. $\phi(\cdot)$ denotes the activation function, which allows a neural network to map nonlinear relationship. 
A critical requirement for the activation function in PINNs is at least twice continuously differentiable, as this regularity is essential for properly computing the differential operators. The hyperbolic tangent function, $tanh(\cdot)$, is frequently employed due to its smoothness and bounded gradient properties. The complete set of trainable parameters is denoted by $\Theta=\{(W^{\left( k \right)}, b^{\left( k \right)})\}_{k=1}^L$, where $W^{\left( k \right)}$ and $b^{\left( k \right)}$ represent the weight matrices and bias vectors for each of the $L$ hidden layers, respectively. %These parameters are iteratively optimized during training to minimize the physics-informed loss function through gradient-based methods.
%In PINNs, the activation function is required at least twice continuously differentiable and is often chosen as the hyperbolic tangent $tanh(\cdot)$.  $\Theta=\{(W^{\left( k \right)}, b^{\left( k \right)})\}_{k=1}^L$ represents the weight and bias of all the $L$ layers, which will be updated during training. 
%The automatic differentiation is utilized to calculate the derivatives of the DNN $f\left( x;\Theta \right)$ with respect to $S$ and $t$ by the backward chain rule.

In PINNs, the neural network solution approximant $\widehat{\boldsymbol{u}}(\boldsymbol{x}, t; \Theta)$ generates the following fundamental residual quantities that measure the discrepancy from the exact solution of the PDE system (\ref{PDE}): 
\begin{equation}
\begin{aligned}
&\boldsymbol{r}_{p}(\boldsymbol{x}, t; \Theta) := \widehat{\boldsymbol{u}}_t\left( \boldsymbol{x}, t;\Theta \right) + \mathcal{N}[\widehat{\boldsymbol{u}}\left( \boldsymbol{x}, t;\Theta \right)]-\boldsymbol{f}(\boldsymbol{x},t), \\
&\boldsymbol{r}_{i}(\boldsymbol{x}; \Theta) :=  \widehat{\boldsymbol{u}}(\boldsymbol{x},0) - \boldsymbol{u}_0(\boldsymbol{x}),\\
&\boldsymbol{r}_{b}(\boldsymbol{x}, t; \Theta) := \widehat{\boldsymbol{u}}\left( \boldsymbol{x},t \right) -\boldsymbol{u}_b\left(\boldsymbol{x},t \right).  
\end{aligned}
\end{equation} 
%These residuals characterize how well the neural network fits the interested PDE. 
The minimization of these residuals during training ensures that the neural network approximant simultaneously satisfies the governing PDE, initial conditions, and boundary constraints. To discrete the loss function defined by the residuals, we define the training dataset as $\mathcal{X}$, which includes the sampling points $\mathcal{X}_i=\{(\boldsymbol{x}_i^n,t_i^n=0) \}_{n=1}^{N_i}$ in the initial domain, the sampling points $\mathcal{X}_b=
\{(\boldsymbol{x}^n_b,t^n_b)\}_{n=1}^{N_b}$ 
in the boundary domain, and collocation points $\mathcal{X}_p=\{(\boldsymbol{x}^n_p,t^n_p)\}_{n=1}^{N_p}$ in the equation domain. During training, these points are typically distributed uniformly across the spatiotemporal domain.
%The approximation function $f\left(x;\Theta \right)$ i.e., the  parameter $\Theta$ is learned with the input points by minimizing the loss function

Then PINNs are aimed at finding $\widehat{\boldsymbol{u}}\left( \boldsymbol{x}, t;\Theta \right)$ by minimizing the empirical $L_2$ loss function

	\begin{align}\label{eq:loss}
		\begin{split}
			&\mathcal{L} \left( \Theta ;\mathcal{X}\right) =
            \lambda _p\mathcal{L} _p\left( \Theta ;\mathcal{X}_p \right) +
            \lambda _b\mathcal{L} _b\left( \Theta ;\mathcal{X}_b \right) +\lambda _i\mathcal{L} _i\left( \Theta ;\mathcal{X}_i\right),
		\end{split}
	\end{align}
	where
    \begin{small}
	\begin{equation}\label{eq:loss3}
			\mathcal{L} _p\left( \Theta ;\mathcal{X} _p \right) =\frac{1}{N_{p}}\sum_{n=1}^{N_p}{| r_p\left(\boldsymbol{x}^n_p,t^n_p;\Theta \right) 
            |^2},\ 
			\mathcal{L}_b\left( \Theta ;\mathcal{X} _b \right) =\frac{1}{N_{b}}\sum_{i=1}^{N_b}{|    
            r_b\left(\boldsymbol{x}_b^n;\Theta \right)
            |^2},\ 
			\mathcal{L} _i\left( \Theta ;\mathcal{X}_i \right) =\frac{1}{N_{i}}\sum_{i=1}^{N_i}{|
            r_i\left(\boldsymbol{x}_i^n,t_i^n;\Theta \right)
            |^2}. 
	\end{equation}
	\end{small}  
%in which the PDE residual is defined as
Here $\lambda_f$, $\lambda_b$, and $\lambda_i$ denote the positive scalar weights. 
%Thus PINNs leverage the inherent constraints and boundary conditions embedded in partial differential equations to expedite network training. 
Through gradient-based optimization of the loss function, PINNs learn parameters $\Theta \text{ that produce }\widehat{\boldsymbol{u}}(\boldsymbol{x},t)$ converging to the PDE solution $\boldsymbol{u}(\boldsymbol{x},t)$.

\subsection{Conservative PDEs}

 In this paper, we consider the energy-conservative PDEs \cite{feng2010symplectic} of the following general form 
 \begin{equation}\label{H_general}
\boldsymbol{u}_t = S \frac{\delta H}{\delta {\boldsymbol{u}}}[\boldsymbol{u}] + \boldsymbol{f}(\boldsymbol{x},t), \quad \boldsymbol{x} \in \Omega \subseteq \mathbb{R}^d, \ t \in [0,T],
\end{equation}
where $S$ denotes a skew-symmetric operator with respect to the $L^2$ inner product, and  the energy functional $H[\boldsymbol{u}]$ is given by 
\begin{equation}\label{equ:enegy}
H[\boldsymbol{u}](t) = \int_{\Omega} \omega(\boldsymbol{x}; \boldsymbol{u}^{(\alpha)}) \, \mathrm{d}\boldsymbol{x}.
\end{equation}
i.e., the integral of the energy density function $\omega$ over the spatial domain with 
$$
\boldsymbol{u}^{(\alpha)}=\left\{\frac{\partial^{|\alpha|} \boldsymbol{u}}{\partial x_1^{\alpha_1} x_2^{\alpha_2} \cdots x_d^{\alpha_d}}: |\alpha| \leq p\right\}, \quad \alpha \in\left(\mathbb{Z}^{+}\right)^d.
$$ 
$\frac{\delta H}{\delta \boldsymbol{u}}[\boldsymbol{u}]$ represents the variational derivative of $H$ defined by 
$$
\left\langle\frac{\delta {H}}{\delta \boldsymbol{u}}[\boldsymbol{u}], \boldsymbol{v}\right\rangle_{L_2}=\left.\frac{\mathrm{d}}{\mathrm{~d} \epsilon}\right|_{\epsilon=0} {H}[\boldsymbol{u}+\epsilon \boldsymbol{v}], 
\quad \forall \boldsymbol{v} \in H^p(\Omega). 
$$
%, where ${u}^{(l)}$ indicates that the dependence extends through the $\alpha$-th order partial derivatives of ${u}$.
Under zero forcing ($f(\boldsymbol{x},t) \equiv 0$) with compatible boundary conditions (e.g., homogeneous Dirichlet or periodic constraints), the PDE system (\ref{H_general}) reduces to a proper Hamiltonian formulation that exactly preserves the Hamiltonian integral $H$ \cite{cano2006conserved, damelin2009energies}. This follows from the skew-symmetry of $S$:
\begin{equation}\label{equ:dH_dt}
\frac{dH}{dt} = \left\langle \frac{\delta H}{\delta {u}} [\boldsymbol{u}], \frac{\partial \boldsymbol{u}}{\partial t} \right\rangle_{L^2} = \left\langle \frac{\delta H}{\delta \boldsymbol{u}} [\boldsymbol{u}], S \frac{\delta H}{\delta \boldsymbol{u}} [\boldsymbol{u}] \right\rangle_{L^2} = 0.
\end{equation}
In our theoretical analysis, we will focus primarily on the integral-preserving Hamiltonian formulation, where the conservation properties hold exactly. For the numerical experiments, we will consider both the energy-conservative case and the generalized PDE form (\ref{H_general}) with non-zero forcing terms, thereby demonstrating the robustness of our methodology across different physical scenarios. 

\begin{example}\label{ex:wave}
Consider the following nonlinear wave (NLW) equation  as an illustrative example. For any given smooth function $G$, we have  
\begin{equation}\label{equ:wave}
\begin{aligned}
& u_{t t}-\Delta u+G^{\prime}(u)=0, \quad(\boldsymbol{x}, t) \in \Omega \times[0, T], \\
& u(\boldsymbol{x}, 0)=u_0(\boldsymbol{x}), \ 
u_t(\boldsymbol{x}, 0)=u_1(\boldsymbol{x}), \quad x \in \Omega, \\
& u(\boldsymbol{x},t)=u_b(\boldsymbol{x},t), \quad (\boldsymbol{x},t) \in \partial \Omega \times [0,T].
%&{u}\left(\boldsymbol{x},t \right) ={B}\left(\boldsymbol{x},t \right) , \  (\boldsymbol{x},t) \in  \partial \varOmega \times  \left[ 0,T \right].
\end{aligned}
\end{equation}
By introducing $u_t=v$, we have $\boldsymbol{u}=(u,v)$. Then the equation can be transformed as a first-order system of PDEs 
$$
\boldsymbol{u}_t=
\left[\begin{array}{cc}
0 & 1 \\
-1 & 0
\end{array}\right]\left[\begin{array}{c}
-\Delta u+G^{\prime}(u) \\
v
\end{array}\right], 
$$
which has a first integral
$$
H[\boldsymbol{u}](t)=\int_{\Omega}\left(\frac{1}{2} v^2+\frac{1}{2}|\nabla u|^2+G(u)\right) \mathrm{d} \boldsymbol{x}. 
$$
%for any constant $C \in \R$. 
It can be written into the form (\ref{H_general}) by noting that
$$
\frac{\delta H}{\delta u}[u]=-\Delta u+G^{\prime}(u), \quad \frac{\delta H}{\delta v}[v]=v.
$$
\end{example}

%\section{Energy-equidistributed moving mesh PINNs (EE-PINNs)}

\section{Energy-Equidistributed moving mesh strategies}\label{sec:EEMPE}

\begin{comment}
We establish the framework in one-dimensional space by considering a vector field $\boldsymbol{u}(x,t)$. Under the canonical parametrization where the spatial coordinate satisfies $x\in[0,1]$, a diffeomorphic mapping between the computational domain and physical domain is constructed as follows:
 \begin{equation}\label{coordinate}
 x=x(\xi,~t),~~\xi \in [0,~1],
 \end{equation}
 subject to the boundary conditions:
 \begin{equation}\label{boundary}
 x(0,~t)=0, ~~x(1,~t)=1.
 \end{equation}
 Using the chain rule, the original PDE (\ref{general}) transforms into:
\begin{equation}
\frac{\partial \boldsymbol{u}}{\partial t} - \frac{\partial \boldsymbol{u}}{\partial x}\cdot\frac{\partial x}{\partial t} = S \frac{\delta H}{\delta \boldsymbol{u}} [\boldsymbol{u}].
\end{equation}
\end{comment}

In this central part of the paper, building upon our previous work on Energy-Equidistributed Principles (EEPs)  \cite{gao2025energy}, we first extend these concepts to multi-dimensional conservative PDEs by developing two key formulations: a location-based Energy-Equidistributed Moving Sampling PDE (EMMPDE) derived through variational methods in computational space \cite{hou}, and then a modified velocity-based EMMPDE designed for enhanced dynamic mesh adaptation. Leveraging this theoretical foundation, we introduce Energy-Equidistributed Moving Sampling PINNs (EEMS-PINNs), a novel deep learning framework that integrates our energy-based adaptive sampling strategy with PINNs.

%we first generalize our previous EEPs to multi-dimensional conservative PDEs. We present one location-based EMMPDE motivated by the variational approach in the computational space \cite{hou}. In addition, a modification velocity-based EMMPDE has also been provided. Then we combine this novel strategy with the most widely used deep learning solver PINNs, resulting in a new deep learning framework called Energy-equidistributed Moving Sampling PINNs(EEMS-PINNs).

%Mesh movement methods include velocity-based and location-based methods. In this paper, we employ a location-based method that utilizes moving mesh PDEs (MMPDEs) to dynamically adjust the mesh via the gradient flow equation derived from an adaptation functional.  

\subsection{Energy-Equidistributed Principles (EEPs)}
Conceptually, MMPDE generates a moving mesh that continuously maps a suitable computational space $\Omega_C$ into the physical domain $\Omega$. The mapping is defined by a bijective coordinate transformation from the computational coordinate  $\boldsymbol{\xi} \in \Omega_C$  to the physical coordinate $\boldsymbol{x} \in \Omega$ through the relation:
\begin{equation}\label{equ:mapping}
  \boldsymbol{x}=\boldsymbol{x}(\boldsymbol{\xi},t).  
\end{equation}
%While the computational domain is typically discretized using a uniform mesh for simplicity, the MMPDE framework permits complete flexibility in selecting alternative mesh configurations for $\Omega_C$.
The moving mesh methodology fundamentally relies on the equidistribution principle. For multi-dimensional problems, this approach naturally extends to a variational formulation, where the coordinate transformation (\ref{equ:mapping}) is derived as the minimizer of a mesh functional. The resulting moving mesh PDE system then emerges as the Euler-Lagrange equations of this optimization problem. This variational perspective provides a unified framework for generating optimally adapted meshes while maintaining desirable geometric characteristics across multiple dimensions, such as uniformity, alignment, and orthogonality. 

In moving mesh methods based on equidistribution, the monitor function bridges the physical solution and the MMPDEs. This function is designed to take larger values in regions where the physical solution requires higher mesh resolution, such as near singularities or sharp gradients like shock waves in the Burgers equation. %Consequently, the selection of the monitor function has been a fundamental research topic. 

While this approach effectively resolves fine-scale features, it inherently overlooks the conservation of energy, a fundamental property of Hamiltonian systems. 
To address this limitation, our proposed Energy-Equidistributed Principles introduce the use of the energy density function as the monitor function. This innovation ensures that mesh adaptation not only tracks spatial solution features but also preserves the Hamiltonian system’s energy dynamics throughout temporal evolution. For clarity, we first reformulate the one-dimensional formulation  in \cite{gao2025energy} through a variational framework before extending to multidimensional cases.

%In this paper, motivated by the work of  one-dimension case in \cite{gao2025energy}, we propose an Energy-equidistributed moving mesh PDE (EMMPDE) for high-dimensional Hamiltonian PDEs that explicitly incorporates energy dynamics. 

%\subsection{Moving mesh PDE (MMPDE)}\label{sec:MMPDE}

%The location based methods directly control the location of mesh points. In this paper we focus on the typical variational method. 

%, and in certain cases(such as when important behavior occurs close to boundaries) it may be useful toconsider a nonuniform but fixed mesh in $\Omega_C$ which is then mapped to a moving mesh in $\Omega$.
%between the physical domain $\Omega$ and , and the mesh is determined by the minimizer of a functional

%Here we will present the mesh moving strategy motivated by the variational approach in the computational space \cite{hou}, which is opposed to the variational approach defined in the physical space. Since in practice, the physical domain may have a very complex geometry, and as a result, we will get more complicated than the Euler–Lagrange equation (\ref{equ:E-L}), which requires more computational effort in obtaining numerical approximations.

\hspace*{\fill}
  
\noindent {\bf 1D case}

Without loss of generality, we set $\Omega_C =\Omega = [0,1]$ and the boundary conditions for the coordinate mapping are explicitly enforced as: 
\begin{equation} \label{boundary}
x(0,t)=0, \ x(1,t)=1, \quad \forall t \in [0,T].
\end{equation}
Assuming that ${u}(x)$ is the solution of  of Equ.(\ref{H_general}) with initial condition ${u}_0(x)$, then the total energy remains invariant over time, satisfying $H(t)= H_{0}$ for all $t>0$, where the initial energy $H_0$ is given as  
$$H_0 \doteq H(0) =  \int^1_0 \omega({x}; {u}_0^{(\alpha)})  \ \mathrm{d}x.$$ 

Then motivated by the equidistribution principle, we start from the uniform distribution samples $\xi_{j}=\frac{j}{N+1} \ (j=0,1,\cdots, N)$, and intend to find a mesh $x_j=x(\xi_j,t) \in [0,1]$ such that 
\begin{equation}\label{equ:H_j}
   H_j(t) \doteq \int^{x_{j+1}}_{x_{j}}\omega({x}; {u}^{(\alpha)})\ \mathrm{d}x=\frac{1}{N}H(t), \ \ j=0, 1,\cdots, N-1. 
\end{equation}
This Energy-Equidistributed Principle ensures preservation of the Hamiltonian structure for each semi-discrete ODE system by employing the energy density function $\omega$ as the monitor function. 
%, maintaining the constant energy $H_j(t)=\frac{1}{N} H_0$, This conservation aligns intrinsically with the symplectic geometry of Hamiltonian systems. 
Rewrite (\ref{equ:H_j}), we have 
\begin{equation}
      \int^{x(\xi_j, t)}_0 \omega({x}; {u}^{(\alpha)}) \ \mathrm{d}x = 
      \xi_j H_0, \ \ j=1,2,\cdots,N-1.
\end{equation}

Consequently, for the continuous transformation of the mesh, $x=x(\xi, t)$,  the Energy-Equidistributed Principle then goes to:  
\begin{equation}\label{eq:energy_d1}
     \int^{x(\xi,~t)}_0 \omega({x}; {u}^{(\alpha)}) \ \mathrm{d}x =\xi H_0, \ \ \xi \in (0,1).
\end{equation}
By computing the first and second derivatives of Equ. (\ref{eq:energy_d1}) with respect to the computational coordinate $\xi$, we derive the following two commonly used one-dimensional EEPs:
\begin{subequations} 
\begin{align}
& \omega({x}; {u}^{(\alpha)})\frac{\partial x}{\partial\xi}=H_0, \\
 & E(x(\xi,t),t) \doteq \frac{\partial}{\partial \xi} \left(\omega({x}; {u}^{(\alpha)}) \frac{\partial x}{\partial\xi}\right)=0. \label{eep2}
 \end{align}
 \end{subequations}
Through the calculus of variations \cite{huang2010adaptive}, we can obtain that Equ. (\ref{eep2}) is exactly the Euler-Lagrange (E-L) equation of the variational problem 
\begin{equation}\label{equ:functional_1D}
\mathcal{I}[u](t)=\frac{1}{2} \int_0^1\left(\omega(x;{u}^{(\alpha)}) \frac{\partial x} {\partial \xi}\right)^2 \  \mathrm{d} \xi
\end{equation}
with respect to the coordinate mapping $x(\xi, t)$. 

On the other hand, the function approximation perspective in adaptive mesh generation,  motivated by Ren and Wang \cite{REN2000246}, focuses on constructing coordinate transformations where solutions exhibit improved regularity properties in computational space. Thus, \cite{hou} propose to look for the mesh map $x=x(\xi, t)$ that minimizes a gradient-dependent functional: 
\begin{equation}\label{equ:H_1D}
   \min_{x} \int_0^1 \omega({x}(\xi,t); {u}^{(\alpha)})\  \mathrm{d}\xi \ \  \text{with} \ 
    \omega=\sqrt{1+ u_x^2x_{\xi}^2}. 
\end{equation}
This variational formulation leads to the Euler-Lagrange equation:
\begin{equation}\label{eq:E-L_H1}
\left(\frac{u_x^2}{\sqrt{1+u_x^2 x_{\xi}^2}} x_{\xi}\right)_{\xi}=\frac{u_x u_{x x} x_{\xi}^2}{\sqrt{1+u_x^2 x_{\xi}^2}} .
\end{equation}
The direct numerical solution of Equ. (\ref{eq:E-L_H1}) becomes computationally intractable due to the strongly nonlinear source term. A simplified formulation emerges by setting the right-hand side of (\ref{eq:E-L_H1}) to zero and replacing $u_x$ with the regularized quantity $u_x x_\xi$, which reduces to the Energy-Equidistribution Principle Equ. (\ref{eep2}) with the particular monitor function $\omega = \sqrt{1 + u_x^2 x_\xi^2}$. This simplified approach maintains essential mesh adaptation properties while improving numerical feasibility. As the generalization to Equ. (\ref{equ:H_1D}), the variational problem can be formulated as minimizing the energy integral
\begin{equation}\label{equ:H_1D_general}
\min_{x} H(t) = \int_0^1 \omega(x(\xi,t); u^{(\alpha)}) \ \mathrm{d}\xi,
\end{equation}
where $\omega$ represents the energy density function associated with the conserved quantity $H$.%, providing a physical foundation for mesh adaptation while ensuring computational tractability.

%However, the extremely large nonlinear source term makes the numerical solution of Equ. (\ref{eq:E-L_H1}) computationally infeasible. So a candidate for a good mesh generator is obtained by setting to zero the right-hand side of Equ.(\ref{eq:E-L_H1}) and replace $u_x$ by the smoother quantity $u_xx_{\xi}$. Then it simply goes to Equ. (\ref{eep2}) with $ \omega=\sqrt{1+ u_x^2x_{\xi}^2}$ as a particular choice of monitor function. More generally, for $t>0$ and any given solution $u$, if we take the energy density function of the energy $H$, then the variational problem is 
%\begin{equation}\label{equ:H_1D_general}
%   \min_{x}  H(t) =\int_0^1 \omega({x}(\xi,t); {u}^{(\alpha)})\  \mathrm{d}\xi.
%\end{equation}
%a nonlinear elliptic equation with a very stiff source term. 

%Besides, in \cite{hou}, for a particular choice of monitor function, Equ.(\ref{eep2}) can alse be seen as a relaxation Euler-Lagrange equation of minimization problem 

\hspace*{\fill}

\noindent{\bf 2D case}

We now extend the aforementioned methodology to multi-dimensional conservative PDE systems. While we present the theory specifically for the two-dimensional case ($d=2$) for clarity, all results generalize naturally to higher dimensions. Given a scalar energy density function $\omega$, we construct the symmetric positive definite matrix-valued monitor functions
$G_k=\omega(\boldsymbol{x}) I_d$ ($k=1,2,\cdots,d$), where $I_d$ denotes the $d$-dimensional identity matrix. The variational formulations (\ref{equ:functional_1D}) and (\ref{equ:H_1D_general}) extend naturally to higher dimensions, providing the theoretical basis for conservative mesh adaptation in multi-dimensional problems.

%Then we will generalize the aforementioned methods to multi-dimensional conservative PDEs. For convenience, we will consider two-dimensional ($d=2$) as the example, but the theory  applies to higher dimensional problems as well.   Given a general scalar energy density function $\omega$, we can define the symmetric positive definite matrix-valued  monitor function  $G_k=\omega(x) I_d$ with identity matrix $I_d$. Then the variational formulations (\ref{equ:functional_1D}) and (\ref{equ:H_1D})  are particularly useful and fundamental when considering the high-dimensional case.  

 %the variational functional (\ref{funcational}) is the generalization of   the 1D (\ref{equ:functional_1D}), and the multi-dimensional counterpart of the simply the Hamiltonian energy (\ref{equ:enegy}). The minimizations of the functional with respect to $x=x(\xi,t)$ go to Equ.(\ref{equ:E-L}) which is the generalizaiton of E-L equation (\ref{eep2}) in 1D case. 
Specifically, the multidimensional generalization of the variational functional (\ref{equ:functional_1D}) takes the form: 
\begin{equation}\label{funcational}
\mathcal{I}[\boldsymbol{u}](t)=\frac{1}{2} \int_{\Omega_C}\left[\nabla x_1^T G_1 \nabla x_1 +\nabla x_2^T G_2 \nabla x_2 \right] \  \mathrm{d} \xi_1 \mathrm{d} \xi_2, 
%\mathcal{I}[x](t)=\frac{1}{2} \sum_{k=1}^d \int_{\Omega_C}\nabla x_k^T G_k \nabla x_k d \xi, 
\end{equation}
where $\nabla=\left(\partial_{\xi_1}, \partial_{\xi_2}\right)^T$. The multi-dimensional counterpart of (\ref{equ:H_1D_general}) is simply minimizing the Hamiltonian energy (\ref{equ:enegy}):
\begin{equation}\label{equ:H_2D}
   \min_{\boldsymbol{x}} \int_{\Omega_C} \omega(\boldsymbol{x}; \boldsymbol{u}^{(\alpha)}) \ \mathrm{d} \boldsymbol{\xi}. 
\end{equation}
%The Energy-Equidistribution Principle (\ref{eep2}) generalizes accordingly \textemdash  in two dimensions yielding a coupled system of elliptic partial differential equations  governing the coordinate transformation $\boldsymbol{x}(\boldsymbol{\xi},t)$:
%The minimizations of the both of the functionals with respect to $x=x(\xi,t)$ go to %Equ.(\ref{equ:E-L}) which is the generalizaiton of E-L equation (\ref{eep2}) in 1D case. %In 2D case, we have the following EMMPDE
Thus, we can naturally extend the Euler-Lagrange equation (\ref{eep2}) to higher dimensions, and get the Energy-Equidistribution Principle in the two-dimensional case: 
\begin{equation}\label{equ:E-L_2d}
\begin{aligned}
%(x_1)_{t} = & -\frac{1}{\tau}E_1 \doteq -\frac{1}{\tau} \left[ \frac{\partial}{\partial \xi_1}\left(G_1 \frac{\partial x_1}{\partial \xi_1}\right)+\frac{\partial}{\partial \xi_2}\left(G_1 \frac{\partial x_1}{\partial \xi_2}\right) \right], \\
%(x_2)_t=& --\frac{1}{\tau} E_2 \doteq -\frac{1}{\tau} \left[ \frac{\partial}{\partial \xi_1}\left(G_2 \frac{\partial x_2}{\partial \xi_1}\right)+\frac{\partial}{\partial \xi_2}\left(G_2 \frac{\partial x_2}{\partial \xi_2}\right) \right]. 
& E_1 \doteq \frac{\partial}{\partial \xi_1}\left(G_1 \frac{\partial x_1}{\partial \xi_1}\right)+\frac{\partial}{\partial \xi_2}\left(G_1 \frac{\partial x_1}{\partial \xi_2}\right)=0, \\
& E_2 \doteq \frac{\partial}{\partial \xi_1}\left(G_2 \frac{\partial x_2}{\partial \xi_1}\right)+\frac{\partial}{\partial \xi_2}\left(G_2 \frac{\partial x_2}{\partial \xi_2}\right)=0.
\end{aligned} 
\end{equation}

We have introduced a mesh adaptation strategy based on a computational-space variational formulation (as opposed to physical-space approaches), following the framework developed in \cite{hou}. 
By operating in the computational domain rather than the physical space, this strategy avoids the numerical complexities inherent in solving the Euler-Lagrange equations (\ref{equ:E-L_2d}) on intricate domains. %, while maintaining precise control over mesh quality through energy-based equidistribution principles. 

%Here we present the mesh moving strategy motivated by the variational approach in the computational space \cite{hou}, which is opposed to the variational approach defined in the physical space. Since in practice, the physical domain may have a very complex geometry, and as a result, we will get more complicated than the Euler–Lagrange equation (\ref{equ:E-L_2d}), which requires more computational effort in obtaining numerical approximations. 

%To provide a mechanism for dynamically adjusting the mesh to possible rapid changes of time-dependent solutions, Huang \cite{huang1994moving} proposed a very robust class of MMPDEs derived from the gradient flow equations associated with the above mesh variational principle. A standard method to solve the diffusion equation (\ref{equ:E-L_2d}) is to consider the following EMMPDE equation
Huang \cite{huang1994moving} introduced a class of robust MMPDEs derived from gradient flow formulations of variational principles to enable dynamic mesh adaptation for time-dependent solutions. The standard implementation for solving (\ref{equ:E-L_2d}) employs the location-based EMMPDE
\begin{equation}\label{EMMPDE}
\boldsymbol{x}_t=-\frac{1}{\tau}\boldsymbol{E},     
\end{equation}
where $\tau > 0$ is a relaxation parameter controlling mesh movement timescales and
$\boldsymbol{E}=\left(E_1, E_2\right)^T$ represents the variational derivatives of the energy functional. The mesh evolves according to this system from an initial configuration toward equilibrium, producing an optimized mesh distribution according to the specified energy-based monitor function while maintaining numerical stability through the gradient flow structure. 

%${E}=\left(E_1, E_2\right)^T$ and $\tau$ being a prescribed parameter which could control the speed of the mesh movement. Then, beginning with an initial guess, we march in “time” to steady state. 

For the mesh evolution over time step  $\tau$, when $\boldsymbol{x}(\boldsymbol{\xi}, t+\tau)$ approaches steady state, the equilibrium condition  $\boldsymbol{E}(\boldsymbol{x}(\boldsymbol{\xi},t+\tau),t+\tau)=0$ holds.  
%Since we have the initial function data of the underlying Hamiltonian PDE, 
Suppose $\boldsymbol{E}$ is smooth enough, a modified EMMPDE can be derived through first-order Taylor expansion of the steady state condition at time $t$:
$$0=\boldsymbol{E}(\boldsymbol{x}(\boldsymbol{\xi},t+\tau),t+\tau)=\boldsymbol{E}(\boldsymbol{x}(\boldsymbol{\xi},t),t)+\tau \left[A \boldsymbol{x}_t(\boldsymbol{\xi},t)+ \boldsymbol{E}_t(\boldsymbol{x}(\boldsymbol{\xi},t),t)\right],  $$
where the Jacobian matrix $A$ is given by:
$$
A=\left(\begin{array}{ll}
\left(E_1\right)_{x_1} & \left(E_1\right)_{x_2} \\
\left(E_2\right)_{x_1} & \left(E_2\right)_{x_2}
\end{array}\right).
$$
Consequently, we have the following velocity-based EMMPDE
\begin{equation}\label{equ:velocity}
    \boldsymbol{x}_t = -A^{-1}\left(  \frac{1}{\tau} \boldsymbol{E}+ \boldsymbol{E}_t  \right).
\end{equation}
The variational strategy EMMPDE (\ref{EMMPDE}) exemplifies location-based mesh movement by directly controlling mesh point positions. In contrast, EMMPDE (\ref{equ:velocity}) incorporates a velocity-based approach through the term 
$\boldsymbol{E}_t$,  targeting mesh velocity explicitly and deriving point locations by integrating the velocity field. For further modifications and extensions of MMPDE-based methods, we refer to \cite{huang1994moving}. 

%Compared to the variational strategy EMMPDE (\ref{EMMPDE}) which is a special examples of the so-called location-based mesh movement strategies, since it controls directly the location of mesh points. EMMPDE (\ref{equ:velocity}) also involves a velocity-based term $E_t$, because it targets directly the mesh velocityand obtains the location of mesh points by integrating the velocity field. For more modification methods involved in  MMPDEs, one can refer \cite{huang1994moving}.    
\begin{remark}
The selection of an appropriate monitor function $\omega$ in MMPDEs remains a challenging and problem-dependent aspect of mesh adaptation, with common choices including solution-derived quantities such as the generalized arc-length measure $\omega=\sqrt{1+c^2\left|\nabla_{\boldsymbol{x}} \boldsymbol{u}(\boldsymbol{x})\right|^2}$ or alternatively $\omega=\sqrt{1+c^2\left|\nabla_{\boldsymbol{\xi}} \boldsymbol{u}(\boldsymbol{x}(\boldsymbol{\xi}))\right|^2}$, 
vorticity-based functions for fluid flows, or curvature-based monitors for geometric features \cite{Monge}. In contrast, Energy-Equidistributed MMPDEs  employ the energy density function as the monitor to guarantee exact energy equidistribution among mesh points. A key requirement for all monitor functions is maintaining strict positivity, which can be ensured through changing $\omega$ into $\omega+C$ using a sufficiently large positive constant $C$. 

%In MMPDE, the choice of an appropriate monitor function $\omega$ is difficult, problem dependent, and the subject of much research. The monitor function  can be determined by apriori considerations of the geometry or of the physics of the PDE solution $u$. An example is the generalized solution arc-length given by
%$M=\sqrt{1+c^2\left|\nabla_x u(x)\right|^2}$ or alternatively $M=\sqrt{1+c^2\left|\nabla_{\xi} u(x(\xi))\right|^2}$  as mentioned earlier. 
%It is often used to construct meshes which can follow moving fronts with locally high gradients \cite{huang2010adaptive}. It is also common to use monitor functions based on the (potential) vorticity, or curvature of the solution. 

%By contrast, in EMMPDE, the energy density function is used as the monitor function, which guarantees energy equidistribution among nodal point. However, it is important to note that the positiveness of the monitor function $\omega$ is required. So one can change $\omega$ into $\omega+C$, where $C$ is a sufficiently large positive constant to ensure $\omega + C > 0$. 
\end{remark}

\subsection{Energy-Equidistributed Moving Sampling PINNs (EEMS-PINNs)}

%Traditionally, the discrete underlying PDE and the mesh equation for the moving mesh method can be solved either simultaneously or alternately by finite difference,  finite element methods and kernel methods etc. In this paper, we will use the neural networks to learn the EMMPDE as well as the underlying PDE with PINN. 

Conventional implementations of moving mesh methods typically solve the coupled system of the underlying PDE and mesh equations either simultaneously or alternately using standard discretization approaches such as finite difference, finite element, or kernel-based methods. Building on the framework of \cite{MMPDE_NET}, we propose a novel alternating learning paradigm where PINNs iteratively solve both the EMMPDE and the underlying physical PDE, termed Energy-Equidistributed Moving Sampling PINNs (EEMS-PINNs). The framework first trains an initial PINNs to approximate the solution, which then informs the construction of the energy-based monitor function for the EMMPDE. This adaptive mesh generation guides subsequent PINNs training through optimized collocation points, accelerating loss convergence while preserving the Hamiltonian structure of the system. Numerical experiments in Section \ref{sec:num} demonstrate that the Energy-Equidistributed principle provides such effective guidance that optimal meshes often emerge after just one iteration, eliminating the need for repeated refinements.

%In this work, motivated by \cite{MMPDE_NET}, we propose a novel paradigm where both the EMMPDE and the underlying physical PDE are learned  alternately through PINN.  The framework employs an initial PINNssolution to inform the energy functional construction in learning the EMMPDE, which then generates optimized sampling points to accelerate PINNsloss convergence. From the numerical experiments, we can see that the energy-equidistributed principle is a strong guiding such that we even do not need multiply iterations of this process to get the optimal meshes.  

%Thus, firstly, as mentioned in Section \ref{sec:intro}, given the initial training collocation points $\mathcal{X}^{ini}=\mathcal{X}^{ini}_p \cup  \mathcal{X}^{ini}_i \cup \mathcal{X}^{ini}_b $, with $\mathcal{X}^{ini}_p=\{(\xi^n_p,t^n_p)\}_{n=1}^{N_p}$, $\mathcal{X}^{ini}_i=\{(\xi^n_i,t^n_i=0)\}_{n=1}^{N_i}$ and  $\mathcal{X}^{ini}_b=\{(\xi^n_b,t^n_b)\}_{n=1}^{N_b}$. One may generate a set of randomly distribution. 
%Without loss of generality,  for better demonstration, we will take the nonlinear wave equation (\ref{equ:wave}) as an example. We can do the pre-training of the underlying Hamiltonian PDE with PINN. Assuming the approximation solution of the Hamitonian PDE being the multilayer network $\left(\hat{{u}}, \hat{{v}}\right)\doteq \widehat{\boldsymbol{u}}(\boldsymbol{x},t;\Theta_H)$ as defined in (\ref{eq:network}),  the training PDE residuals of the nonlinear wave equation are reformulated into two equations

Following the initialization procedure described in Section \ref{sec:intro}, we begin with a set of uniformly distributed collocation points across the computational domain: $\mathcal{X}^{ini}=\mathcal{X}^{ini}_p \cup  \mathcal{X}^{ini}_i \cup \mathcal{X}^{ini}_b $, with $\mathcal{X}^{ini}_p=\{(\boldsymbol{\xi}^n_p,t^n_p)\}_{n=1}^{N_p}$, 
$\mathcal{X}^{ini}_i=\{(\boldsymbol{\xi}^n_i,t^n_i=0)\}_{n=1}^{N_i}$ and 
$\mathcal{X}^{ini}_b=\{(\boldsymbol{\xi}^n_b,t^n_b)\}_{n=1}^{N_b}$. Using the nonlinear wave equation (\ref{equ:wave}) as our model system, we pre-train an approximation of the solution using PINN with a fully connected multilayer network architecture $\left(\hat{{u}}, \hat{{v}}\right)\doteq \widehat{\boldsymbol{u}}(\boldsymbol{\xi},t;\Theta_H)$  defined in (\ref{eq:network}). The PDE residuals are formulated as two coupled components: 
\begin{equation}
\begin{aligned}
&{r}_{p,1}(\boldsymbol{\xi}, t; \Theta_H) = \widehat{u}_t\left( \boldsymbol{\xi}, t;\Theta_H \right)  - \widehat{v} \left( \boldsymbol{\xi}, t;\Theta_H \right),\\
&{r}_{p,2}(\boldsymbol{\xi}, t; \Theta_H) = \widehat{v}_t\left( \boldsymbol{\xi}, t;\Theta_H \right) - \Delta \widehat{u}\left( \boldsymbol{\xi}, t;\Theta_H \right) +G'(\widehat{u}\left( \boldsymbol{\xi}, t;\Theta_H \right)).
\end{aligned}
\end{equation} 
The optimal parameter $\Theta^*_H$ of the output $\left(\hat{{u}}, \hat{{v}}\right)$ can be learned by  minimizing the following loss function 
 \begin{small}
\begin{align}\label{loss_wave}
		\begin{split}
			\mathcal{L}_{H} \left( \Theta_H ;\mathcal{X}^{ini}\right) = &
            \lambda_{p,1}\mathcal{L}_{p,1}\left( \Theta_H ;\mathcal{X}^{ini}_p \right) +
            \lambda_{p,2}\mathcal{L}_{p,2}\left( \Theta_H ;\mathcal{X}^{ini}_p \right) +
            \lambda_{b}\mathcal{L}_{b}\left( \Theta_H ;\mathcal{X}^{ini}_b \right) + \\
           & \lambda_{i,1}\mathcal{L}_{i,1}\left( \Theta_H ;\mathcal{X}^{ini}_i\right)+
            \lambda_{i,2} \mathcal{L}_{i,2}\left( \Theta_H ;\mathcal{X}^{ini}_i\right),
		\end{split}
	\end{align}
    \end{small}
	where
    \begin{small}
    \begin{align}\label{loss5_wave}
		\begin{split}
			&\mathcal{L} _{p,1}\left( \Theta_H ;\mathcal{X}^{ini}_p \right) =\frac{1}{N_{p}}\sum_{n=1}^{N_p}{| r_{p,1}\left(\boldsymbol{\xi}^n_p,t^n_p;\Theta_H \right) 
            |^2},\  \mathcal{L} _{p,1}\left( \Theta_H ;\mathcal{X}^{ini}_p \right) =\frac{1}{N_{p}}\sum_{n=1}^{N_p}{| r_{p,1}\left(\boldsymbol{\xi}^n_p,t^n_p;\Theta_H \right) 
            |^2}, \\ 
			&\mathcal{L}_{i,1}\left(\Theta_H ;\mathcal{X}^{ini}_i \right) =
            \frac{1}{N_{i}}\sum_{i=1}^{N_i}{|    
            %r_i\left(x_i^n;\Theta_H \right)
           \hat{u}(\boldsymbol{\xi}_i^n,0)-u_0(\boldsymbol{\xi}_i^n)
            |^2},\ 
			\mathcal{L} _{i,2}\left( \Theta_H ;\mathcal{X}^{ini}_i \right) =\frac{1}{N_{i}}\sum_{i=1}^{N_i}{|
            %r_i\left(x_i^n,t_i^n;\Theta_H \right)
            \hat{v}(\boldsymbol{\xi}_i^n,0)-u_1(\boldsymbol{\xi}_i^n)
            |^2}, \\
            & \mathcal{L}_b\left( \Theta_H ;\mathcal{X}^{ini} _b \right) =
            \frac{1}{N_{b}}\sum_{i=1}^{N_b}{|    
            %r_b\left(x_b^n;\Theta_H \right)
            \hat{u}(\boldsymbol{\xi}_b^n,t_b^n)-u_b(\boldsymbol{\xi}_b^n, t_b^n)
            |^2}.
	\end{split}
	\end{align}
\end{small}
% It then will feed EMMPDE the prior information of the physical solution to help to define the energy density function.

The pre-trained PINN solution $(\hat{u},\hat{v})$ provides essential prior information for the Energy-Equidistributed moving sampling method, specifically enabling the construction of the  energy density monitor function:
$$\omega(\boldsymbol{\xi};\hat{u},\hat{v})= \frac{1}{2} \hat{v}^2+\frac{1}{2}|\nabla \hat{u}|^2+G(\hat{u}).$$ 
%with pre-trained network, we proceed to learn the adaptive coordinate transformation 
With this physics-informed $\omega$ and and initial collocation points $\mathcal{X}^{ini}$, we then learn the adaptive coordinate transformation 
 $\boldsymbol{x}=\boldsymbol{x}(\boldsymbol{\xi},t)$ through the EMMPDE deep neural network $\hat{\boldsymbol{x}}(\boldsymbol{\xi},t; \Theta_E)$  with $\Theta_E$ being the EMMPDE network parameter.   

To enforce the geometric boundary condition $\boldsymbol{x}(\boldsymbol{\xi}, t) = \boldsymbol{\xi}$ for $\boldsymbol{\xi} \in \partial \Omega$, during mesh adaptation learning, we design a structured neural network architecture inspired by \cite{SUKUMAR2022114333}. The boundary-constrained coordinate transformation is formulated as:
 %Specially, we should point out that according to the boundary condition like (\ref{boundary}), during the mesh adaptive learning of the network $\hat{\boldsymbol{x}}(\boldsymbol{\xi},t; \Theta_E)$, we always set  $\boldsymbol{x}(\boldsymbol{\xi}, t) = \boldsymbol{\xi}$ for $\boldsymbol{\xi} \in \partial \Omega$. Thus in order to impose this boundary conditions in the deep neural network, we introduce the solution structure for Dirichlet boundary conditions \cite{SUKUMAR2022114333}
\begin{equation}\label{equ:net_boundary}
\hat{\boldsymbol{x}}^{bc}(\boldsymbol{\xi},t; \Theta_E) = \boldsymbol{\xi} +\psi(\boldsymbol{\xi}) \hat{\boldsymbol{x}}(\boldsymbol{\xi},t; \Theta_E),    
\end{equation}
where $\psi$ is the normalizer function satisfying $\psi =0$ on $\partial \Omega$. %The structure of $\phi$ can be defined by the approximation distance function which give the shortest distance between any point in $\mathbb{R}^d$ to $\partial \Omega$. For more details one can refer  \cite{SUKUMAR2022114333}. For example,  we can see that for one dimensional case with  $\Omega_C=\Omega=[0,1]$, the boundary condition (\ref{boundary}) are $\xi=0$ and $\xi=1$. The normalizer function  can be either the following choices:  $$\phi_A(\xi)=\phi_1(\xi) \phi_2(\xi), \quad \phi_B(\xi)=\phi_1(\xi)+\phi_2(\xi)-\sqrt{\phi_1^2(\xi)+\phi_2^2(\xi)}, \quad \phi_C(\xi)=\frac{\phi_1(\xi) \phi_2(\xi)}{\sqrt{\phi_1^2(\xi)+\phi_2^2(\xi)}},$$ with $\phi_1=\xi, \ \phi_2=1-\xi$. 
The normalizer function $\psi(\boldsymbol{\xi})$ can be constructed using the signed distance function to the boundary $\partial \Omega$ \cite{SUKUMAR2022114333}. In the 1D case with $\Omega = [0,1]$, where boundary conditions (\ref{boundary}) require $\xi = 0$ and $\xi = 1$, we consider three common normalizer constructions:
\begin{align}
\psi_A(\xi) &= \xi(1-\xi) \quad \text{(Product form)}  \notag\\ 
\psi_B(\xi) &= \xi + (1-\xi) - \sqrt{\xi^2 + (1-\xi)^2} \quad \text{(R-equivalence)}  \notag \\
\psi_C(\xi) &= \frac{\xi(1-\xi)}{\sqrt{\xi^2 + (1-\xi)^2}} \quad \text{(Normalized product)}  \notag
\end{align}
where $\psi^{(1)}(\xi) = \xi$ and $\psi^{(2)}(\xi) = 1-\xi$ represent the signed distance functions to each boundary. These constructions satisfy the essential properties:
$\psi(0) = \psi(1) = 0$ (boundary annihilation) and
$\psi(\xi) > 0$ for $\xi \in (0,1)$ (interior positivity). For multi-dimensional domains, the construction generalizes via R-functions:
\begin{equation}
\psi(\boldsymbol{\xi}) = \prod_{k=1}^d \psi_k(\xi_k)
\end{equation}
where $\psi_k$ are 1D normalizers for each dimension.

%Then $\hat{\boldsymbol{x}}_t^{bc}(\boldsymbol{\xi}_p^n,t_p^n; \Theta_E^*)$ will be learned by minimizing the loss function of the EMMPDE (\ref{EMMPDE}) defined by 
The optimal boundary-constrained mesh $\hat{\boldsymbol{x}}_t^{bc}(\boldsymbol{\xi},t; \Theta^*_E)$ is learned by minimizing the EMMPDE loss function 
    \begin{align}\label{loss_EMMPDE}
		\begin{split}
			&\mathcal{L}_{E} \left(\Theta ;\mathcal{X}^{ini}\right) =
            \frac{1}{N_{p}}\sum_{i=1}^{N_p}{\left|    
           \hat{\boldsymbol{x}}_t^{bc}(\boldsymbol{\xi}_p^n,t_p^n)-\frac{1}{\tau}E(\hat{\boldsymbol{x}}^{bc}(\boldsymbol{\xi}_p^n,t_p^n),t_p^n)
            \right|^2}.\ 
            %\mathcal{L}_{p}\left(\Theta;\mathcal{X}^{ini}_p \right)
            %\omega_{p,2}\mathcal{L}_{p,2}\left( \Theta ;\mathcal{X}^{ini}_p \right) +
            %\omega_{b}\mathcal{L}_{b}\left( \Theta ;\mathcal{X}^{ini}_b \right) +
           % \omega_{i,1}\mathcal{L}_{i,1}\left( \Theta ;\mathcal{X}^{ini}_i\right)+
           % \omega_{i,2} \mathcal{L}_{i,2}\left( \Theta ;\mathcal{X}^{ini}_i\right),
		\end{split}
	\end{align}
This formulation naturally extends to alternative EMMPDE variants like (\ref{equ:velocity}) through appropriate modification. 

In the last step, the EEMS-PINN framework completes its adaptive cycle by feeding the optimized sampling points $\mathcal{X}^{new} = \hat{\boldsymbol{x}}^{bc}(\mathcal{X}^{ini}; \Theta_E^*)$ back into the physics-informed neural network. These energy-adapted collocation points, now concentrated in regions of high dynamical activity, enable the PINN to produce an improved solution approximation $\widehat{\boldsymbol{u}}(\boldsymbol{x},t; \Theta_H^{*})$ through minimization of the loss function (\ref{loss_wave}). %The entire process typically requires only 1-2 iterations to reach optimal sampling distribution for most benchmark problems.  
Numerical experiments in Section \ref{sec:num} demonstrate that the Energy-Equidistributed principle provides such effective guidance that optimal meshes often emerge after just one or two iterations for most benchmark problems.

%EMMPDE will feed new sampling points back to PINNs to output better approximate solutions. The input is the new sampling points $\mathcal{X}^{new}=\hat{x}^{bc}(\mathcal{X}^{ini}; \Theta_E^*)$. Then we obtain the new approximation $\widehat{\boldsymbol{u}}(\boldsymbol{x},t; \Theta_H^{*})$ of the underlying Hamiltonian PDE by minimizing the loss function defined in (\ref{loss_wave}) again.

To rigorously validate the energy conservation law during network training, we implement a numerical quadrature scheme that evaluates the energy functional $H$ using the PINN solution approximation $\widehat{\boldsymbol{u}}(\boldsymbol{x},t; \Theta_H^{*})$. For any given time snapshot $t^*$, we define the quadrature point set:
$
\mathcal{X}^{t^*} = \{ (\boldsymbol{x},t) \in \Omega \times [0,T] \, | \, t=t^* \}.
$
The discrete energy approximation is then computed as:
\begin{equation}
    H_d[\widehat{\boldsymbol{u}}](t^*)=\sum_{(\boldsymbol{x}_n,t^*)\in \mathcal{X}^{t^*}} \alpha_n  w(\boldsymbol{x}_n;\widehat{\boldsymbol{u}}^{(\alpha)}), 
\end{equation}
where 
$\boldsymbol{x}_n$ are quadrature points
  and   $\alpha_n$ are corresponding quadrature weights. 
The relative energy error metric:
\begin{equation}
\Delta H_d(t^*) = \frac{|H_d(t^*) - H_{d,0}|}{H_{d,0}}
\end{equation}
quantifies conservation law preservation, where $H_{d,0}$ is the reference energy computed using the exact initial condition $\boldsymbol{u}_0(\boldsymbol{x})$.

\begin{figure}[htp]
        \centering
    \includegraphics[width=1\linewidth]{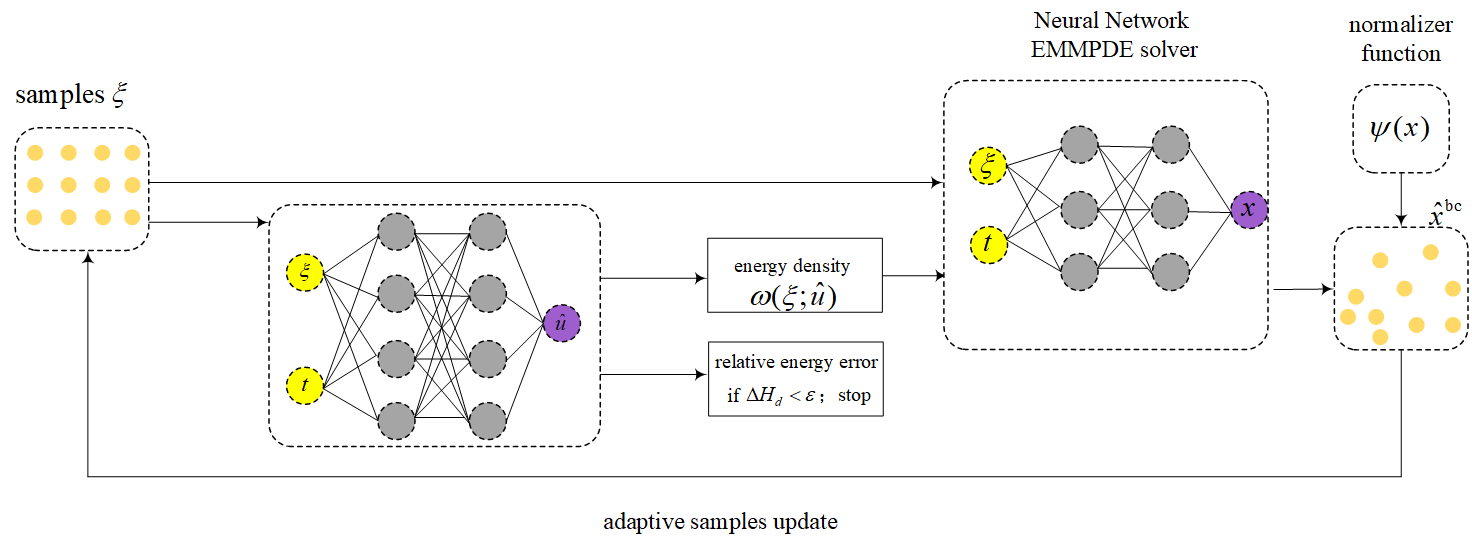}
    \caption{Schematic of EEMS-PINNs. %$N(\cdot)$ and $\partial_{t}$ represents the differential operators; $r$ represents the residual of the PDEs.
    }
    \label{fig:flowchart}
\end{figure}

The EEMS-PINN framework establishes a self-consistent iterative cycle that couples neural network approximation with variational mesh adaptation. As detailed in Algorithm~1 and visualized in Figure ~\ref{fig:flowchart}, each iteration performs three fundamental operations. This adaptive process automatically achieves optimal point concentration in high-energy regions while maintaining the Hamiltonian structure.

% % ===== 算法部分 =====

\begin{algorithm}[htp]
		\caption{Energy-Equidistributed Moving Sampling PINNs (EEMS-PINNs)}
		\label{alg:AOS}
		\renewcommand{\algorithmicrequire}{\textbf{Input:}}
		\renewcommand{\algorithmicensure}{\textbf{Output:}}
		\begin{algorithmic}[1]
			\REQUIRE   Number of collocation points $N_{p}$, $N_i$ and $N_b$;  
                        loss weights for conservative PDE $\lambda_{p}$, $\lambda_{b}$ and $\lambda_i$; EMMPDE time step $\tau$; maximum iterations of pre-training, mesh moving and re-training  $M_{1}, M_2,$ and $M_3$; initial parameters $\Theta_E^0$ and $\Theta_H^0$ for EMMPDE and underlying conservative PDE respectively; 
                       energy tolerance $\epsilon$.  
            %$N_{rn}$ , $N_{rm}$ , the nonlinear Black-Scholes equation , RAR-count.  %%input
			\ENSURE $\Theta^*_H$    %%output
            \STATE Sample the initial collocation points $\mathcal{X}^{ini}=\mathcal{X}^{ini}_p \cup \mathcal{X}^{ini}_i \cup \mathcal{X}^{ini}_b$ with the prior uniform distribution. 
            \STATE  {\bf Pre-training:} \\
            \STATE  With initial parameter $\Theta_H^0$, train a fully connected neural network $\widehat{\boldsymbol{u}}(\boldsymbol{x},t;\Theta_H)$  for a limited
                    number of iterations $M_1$ by minimizing underlying PDE loss function $\mathcal{L}_H(\Theta_H; \mathcal{X}^{ini})$ as defined in Equ.(\ref{loss_wave}).  \\
            \STATE Output the PDE solution estimator $\widehat{\boldsymbol{u}}(\boldsymbol{\xi},t;\Theta^*_H)$.
            \STATE Compute the energy density function $\omega(\boldsymbol{\xi}; \widehat{\boldsymbol{u}}^{(\alpha)}).$ %, and the initial Hamiltonian $H_0$. 
            \STATE{\bf If }the maximum relative energy error:  $\max \Delta H_d < \epsilon$,  {\bf then} stop; {\bf else} continue.  
            \STATE {\bf Mesh Moving: }
            \STATE Train the neural network $\hat{\boldsymbol{x}}^{bc}(\boldsymbol{\xi},t; \Theta_E) = \boldsymbol{\xi} +\psi(\boldsymbol{\xi}) \hat{x}(\boldsymbol{\xi},t; \Theta_E)$  for a limited 
                    number of iterations $M_2$ by minimizing EMMPDE loss function $\mathcal{L}_E(\Theta_E; \mathcal{X}^{ini})$ as defined in Equ.(\ref{loss_EMMPDE}).  \\
            \STATE Output the new mesh points $\mathcal{X}^{new}=\hat{\boldsymbol{x}}^{bc}(\mathcal{X}^{ini},t; \Theta_E^*)$. 
            \STATE {\bf Re-training:} 
            \STATE Update the maximum training iterations $M_1 \leftarrow M_3$.
            \STATE Update the samples $\mathcal{X}^{ini} \leftarrow \mathcal{X}^{new}.$
            \STATE Transfer conservative PDE solution network parameter:  
            $\Theta^0_H \leftarrow \Theta_H^*$.
            \STATE  Re-train the neural network $\widehat{\boldsymbol{u}}(\boldsymbol{x},t;\Theta_H)$: Back to lines 3-5.  
		\end{algorithmic}
	\end{algorithm}

 \section{Numerical Experiments}\label{sec:num}

%the effectiveness of the proposed EEMS-PINNs is substantiated by the following numerical experiments in one- and two-dimensional cases. We will compare EEMS-PINNs proposed in this paper with PINNs \cite{raissi2019physics} and the adaptive collocation point method for PINNs based on monitor function (WAM-PINN) \cite{hou2023enhancing} in the following examples.  In WAM-PINNs, the sample points follow a prior probability density function $p_0(\boldsymbol{x})$ (e.g., the uniform distribution). The movable collocation points in the computational domain are then resampled based on the probability density function:

In this section, the numerical experiments will validate the effectiveness of our Energy-Equidistributed Moving Sampling PINNs (EEMS-PINNs) through systematic comparisons with standard PINNs \cite{raissi2019physics} and Adaptive Point Movement PINNs (WAM-PINNs) \cite{hou2023enhancing} on both 1D and 2D test cases. In the WAM-PINNs framework, collocation points are initially distributed according to a prior probability density $p_0(\boldsymbol{x})$, typically chosen as a uniform distribution over the computational domain. The adaptive sampling process then dynamically adjusts the point locations through resampling from the modified probability density:
\begin{equation}\label{eq:WAM}
p(\boldsymbol{x}) = \frac{w^k(\boldsymbol{x})}{\mathbb{E}[w^k(\boldsymbol{x})]}, \quad \text{with} \  w(\boldsymbol{x}) = \sqrt{1+|\nabla {u}|^2},
\end{equation}
where the gradient $\nabla {u}=
(\partial_{x_1}, \partial_{x_2},\dots,\partial_{x_d}, \partial_t)^T$ includes both spatial and temporal variables of the PDE system and the hyperparameter exponent $k \geq 1$ controls the degree of point concentration. 
The normalization constant $\mathbb{E}[w^k(\boldsymbol{x})]$ is estimated numerically via Monte Carlo integration of the monitor function over the domain. 

    %with the monitor function defined by $ w(x)=\sqrt{1+|\nabla u|^2}$ and $k$ being a hyperparameter.  $\mathbb{E} [{w^k\left(\boldsymbol{x} \right)}]$ can be approximated by Monte Carlo integration.   

%The collocation points are initialized as uniformly sampled points within the computational domain for both PINN and WAM-PINN. The initial sampling points of EEMS-PINN are uniform grid in the space-time domain. 

For PINN and WAM-PINN, collocation points are initially generated through uniform random sampling across the entire computational domain. In contrast, EEMS-PINN employs a structured initialization with points arranged on a uniform space-time grid before undergoing energy-based adaptation. 
The neural network architecture throughout this study employs the deep fully-connected structure defined in Equ.\eqref{eq:network}, with Tanh activation functions in all hidden layers to ensure $C^2$ continuity required for accurate derivative computations in the PDE residuals. For the loss function of the underlying conservative PDE, we maintain uniform weighting coefficients ($\lambda_p = \lambda_b = \lambda_i = 1$) across all constraint terms - governing equations, boundary conditions, and initial conditions - ensuring balanced enforcement of all physical constraints during optimization. 
%The backbone of all the neural networks used in the training is the deep fully connected network as shown in Equ.(\ref{eq:network}). The Tanh activation function is employed in the hidden layers for the trainings. The weights of the loss function are set to be $\lambda_p=\lambda_b=\lambda_i=1$ in the loss of the underlying conservative PDE. 

%The optimization algorithms are implemented by Adam and then LBFGS.  The LBFGS optimizer can become trapped in local minima. To mitigate this issue, as it is mentioned in \cite{lu2021deepxde} we apply the Adam optimizer for a predetermined number of iterations at the beginning of each training phase. The use of Adam helps guide the parameters towards a more favorable initial state. After this initial phase, we use LBFGS optimizer to expedite convergence. All the hyperparameters are listed in Table \ref{Table:para}.

Our optimization approach combines the benefits of Adam and L-BFGS algorithms in a two-phase training strategy. Initially, we employ the Adam optimizer for a predetermined number of iterations to navigate the parameter space and avoid poor local minima, following the recommendations in \cite{lu2021deepxde}. Subsequently, we switch to the L-BFGS optimizer to exploit its superlinear convergence properties near local optima, significantly accelerating the final stages of training. All hyperparameters including learning rates, iteration counts, and the sizes of the networks are carefully documented in Table~\ref{Table:para}.
\begin{table}[!htp]
		\centering
        \footnotesize
		%\captionsetup{labelsep=period,labelfont=bf}
\begin{tabular}{l|p{3cm}llll}
\hline Example &  PDE NN  &  EMMPDE NN   & Learning rate & Numbers of Iter \\
& (Depth$\times$Width) & (Depth $\times$ Width) & Adam/L-BFGS & $(M_1, M_2, M_3)$ \\
\hline 1 & $5 \times 40 $ & $4 \times 20 $&   0.001/$\sim$ &3000/3500/3500\\
 2 & 5 $\times $20 & 4 $\times$ 20 &  0.001/0.5 &5000/6000/10000\\
3 & 5 $\times 50 $& 4 $\times$ 20 &  0.001/0.005 & 5000/2500/10000\\
4 &  8 $\times$ 80 & 4 $\times$ 30 &   0.001/0.5 & 4000/3000/6000 \\
5 & 6 $\times$ 50 & 4 $\times$ 40 & 0.001/0.5 & 4000/2500/4000 \\
6 & 5 $\times$ 40 & 4 $\times$ 20 & 0.001/0.005& 4000/500/4000 &\\
\hline
\end{tabular}
\caption{Hyperparameters used in the following Examples.}\label{Table:para}
\end{table}

%The fully connected network used in PINNscontains five hidden layers with 20 neurons in each hidden layer, while the residual function $F(x,\{w^{(i)}\})$ in Equ.(\ref{Resnet}) has two hidden layers with 20 neurons in each hidden layer.  The number of iteration rounds for RAM-PINN, WAM-PINNsand AM-PIRN are all set to be $10$. All the hyper-parameters were maintained consistently in all experiments, unless stated otherwise. 

The assessment of prediction accuracy is carried out through the relative $L_2$ error between the reference solution $u(\boldsymbol{x},t)$ and the predicted solution $\widehat{u}(\boldsymbol{x},t)$,  
    
    \begin{align}\label{eq:l2_err}
		\begin{split}
			L_2\text{-}error=\frac{\sqrt{\sum_{i=1}^N{\left| \widehat{u}\left( \boldsymbol{x}_i,t_i \right) -u\left( \boldsymbol{x}_i,t_i \right) \right|}^2}}{\sqrt{\sum_{i=1}^N{\left| u\left( \boldsymbol{x}_i,t_i \right) \right|}^2}}
		\end{split}
    \end{align}    
for given test uniform grid points $\{(\boldsymbol{x}_i,t_i)_{i=1}^N \}\in \Omega \times [0, T]$.

\subsection{One-dimensional equations}

\noindent {\bf Example 1 (Klein-Gordon equation)}  

We firstly consider the nonlinear cubic Klein-Gordon equation \cite{pekmen2012differential}
\begin{equation}\label{equ:KG}
u_{t t}+\beta u_{x x}+\alpha u - \gamma u^3=0,\ \    (x,t) \in [a,b] \times [0,+\infty). 
\end{equation}
We take $\beta=-\alpha^2$, $-10\leq x \leq 10$,  and  the initial conditions
$$
\begin{aligned}
& u(x, 0)=\sqrt{\frac{\alpha}{\gamma}} \tanh (\kappa x), 
\ \  u_t(x, 0)=-c \sqrt{\frac{\alpha}{\gamma}} \kappa \operatorname{sech}^2(\kappa x),
\end{aligned}
$$
where $\kappa=\sqrt{\frac{\alpha}{2\left(c^2-\alpha^2\right)}}$, and $\alpha, \gamma, c^2-\alpha^2>0$. Then we have the exact solution
$$
u(x, t)=\sqrt{\frac{\alpha}{\gamma}} \tanh (\kappa(x-c t)), 
$$
from which the Neumann boundary conditions may be extracted as
$$
\begin{aligned}
& u_x(-10, t)=\kappa \sqrt{\frac{\alpha}{\gamma}} \operatorname{sech}^2(\kappa(-10-c t)), \ 
u_x(10, t)=\kappa \sqrt{\frac{\alpha}{\gamma}} \operatorname{sech}^2(\kappa(10-ct)).
\end{aligned}
$$
Then the energy of the cubic  Klein-Gordon equation (\ref{equ:KG}) is given by
\begin{equation}\label{equ:Klein_energy}
H(t)
=\frac{1}{2} \int_{} \left[u_t^2-\alpha^2 u_x^2+\alpha u^2-\frac{\gamma u^4}{2}\right] {\rm d} {x}.
\end{equation}

In this experiment, we set $\alpha = 0.1$, $\gamma = 1$ and $c=0.3$.  The computational domain spans the spatial interval $[-10, 10]$ and temporal interval $[0, 12]$, with the training dataset constructed through sampling of three distinct point sets: $N_i=100$ initial condition points, $N_b=100$ boundary condition points, and $N_p=1000$ collocation points within the full spatiotemporal domain. 

%The spatial domain is $[-10, 10]$ and the temporal domain is from $t=0$ to $t=12$. The training dataset consists of  initial points ($N_i = $),  boundary points ($ N_b = $), and collocation points ($N_p = $) uniformly sampled across the computational domain. 

%Figure  \ref{fig:Klein_sol} illustrates the exact solution against numerical approximations from PINN, WAM-PINN and EEMS-PINN, respectively. Under identical training conditions, the first row of Figure  \ref{fig:Klein_points} demonstrates that EEMS-PINN achieves a significantly lower absolute error compared to PINN and WAM-PINN.  Quantitatively, the maximum absolute errors for WAM-PINN and EEMS-PINN are of the order of $0.01$ and $0.001$, respectively—all markedly lower than PINN’s maximum absolute error of the order of $0.1$. The comparison between the analytic solutions and the approximation solutions of all the methods at different time is shown in Figure  \ref{fig:Klein_sol_time}. Our method gives the best approximation during the time revolution, even in $t=12$. 

\begin{figure}[!htbp]  
	\centering
	
	% 创建第一行第一列的子图  
	\begin{subfigure}{.24\textwidth}  
		\centering  
		\includegraphics[width=1.00\linewidth]{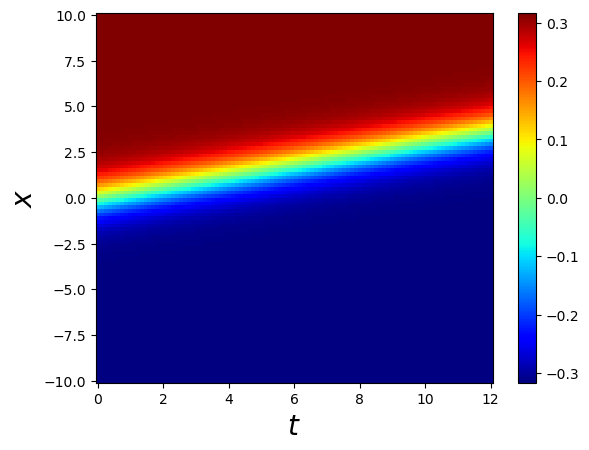}  
		\caption*{Exact}   
	\end{subfigure}%  
	% 创建第一行第二列的子图  
	\begin{subfigure}{.24\textwidth}  
		\centering  
		\includegraphics[width=1.00\linewidth]{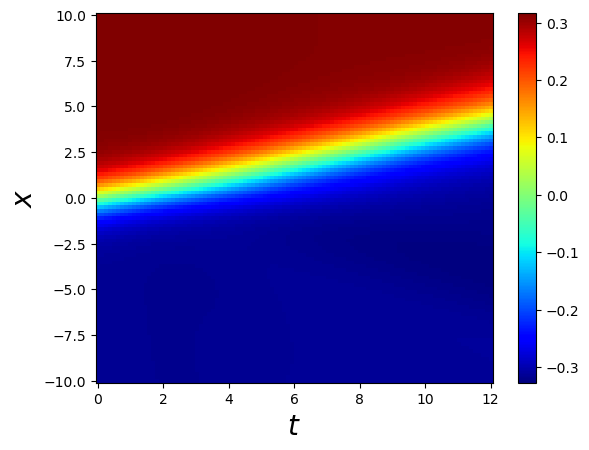}  
		\caption*{PINN}  
	\end{subfigure}  
	% 创建第一行第三列的子图  
	\begin{subfigure}{.24\textwidth}  
		\centering  
		\includegraphics[width=1.00\linewidth]{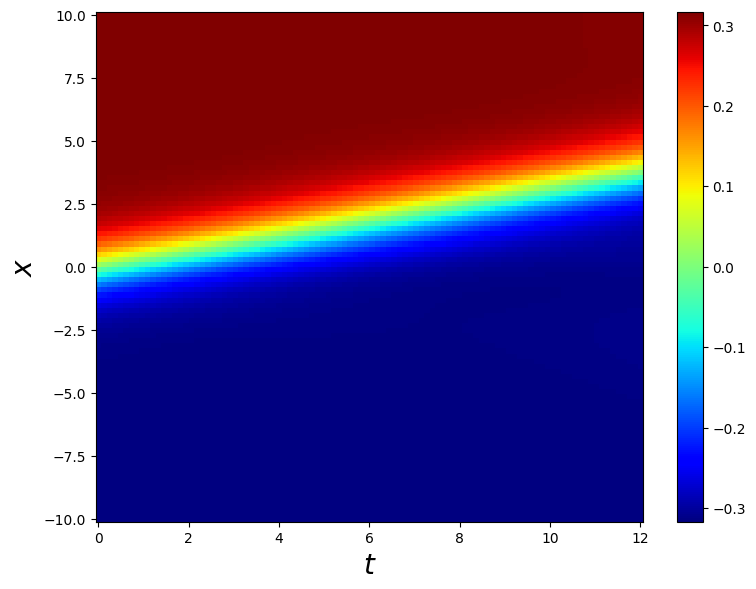}  
		\caption*{WAM-PINN}   
	\end{subfigure}  
		\begin{subfigure}{.24\textwidth}  
		\centering  
		\includegraphics[width=1.00\linewidth]{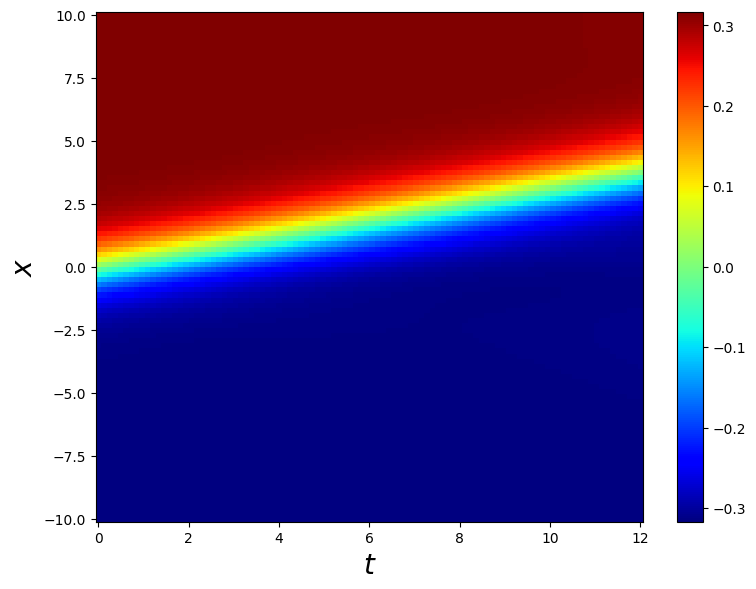}  
		\caption*{EEMS-PINN}  
	\end{subfigure} 
	% 为整个图像阵列添加标题和标签  
	\caption{The exact solution and numerical solutions of PINN, WAM-PINN and EEMS-PINN of  Klein-Gordon equation (\ref{equ:KG}), receptively.} 
    % 整个图像阵列的标题  
	\label{fig:Klein_sol} % 整个图像阵列的标签  
	
\end{figure}

\begin{figure}[!htbp]  
	\centering
	% 创建第一行第一列的子图  
	\begin{subfigure}{.24\textwidth}  
		\centering  
		\includegraphics[width=1.00\linewidth]{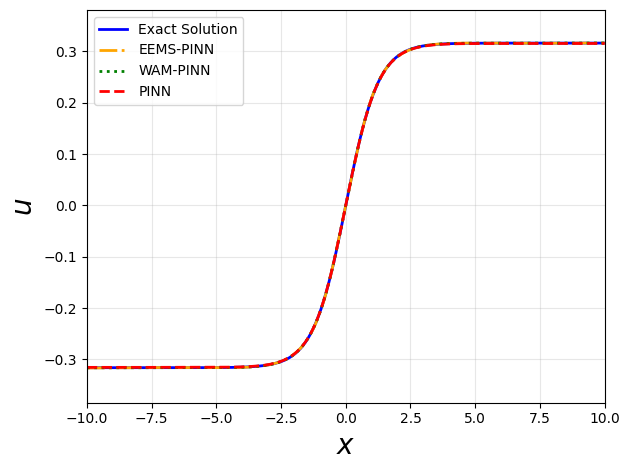}  
		\caption*{$t=0$}   
	\end{subfigure}%  
	% 创建第一行第二列的子图  
	\begin{subfigure}{.24\textwidth}  
		\centering  
		\includegraphics[width=1.00\linewidth]{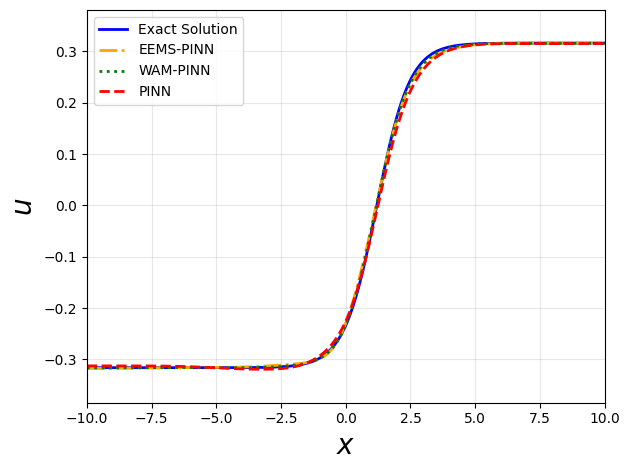}  
		\caption*{$t=4$}  
	\end{subfigure}  
	% 创建第一行第三列的子图  
	\begin{subfigure}{.24\textwidth}  
		\centering  
		\includegraphics[width=1.00\linewidth]{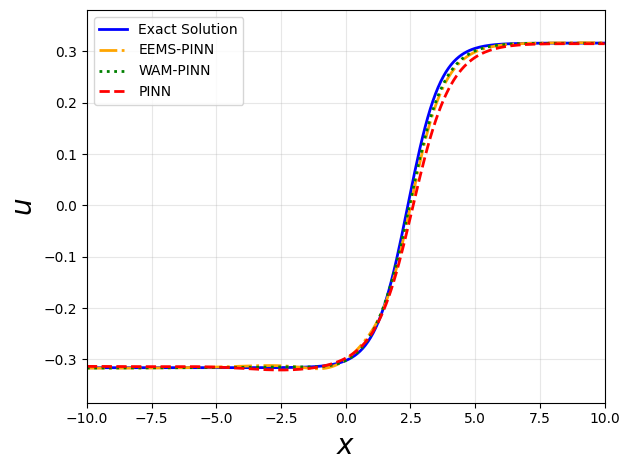}    
		\caption*{$t=8$}   
	\end{subfigure}  
        \begin{subfigure}{.24\textwidth}  
		\centering  
		\includegraphics[width=1.00\linewidth]{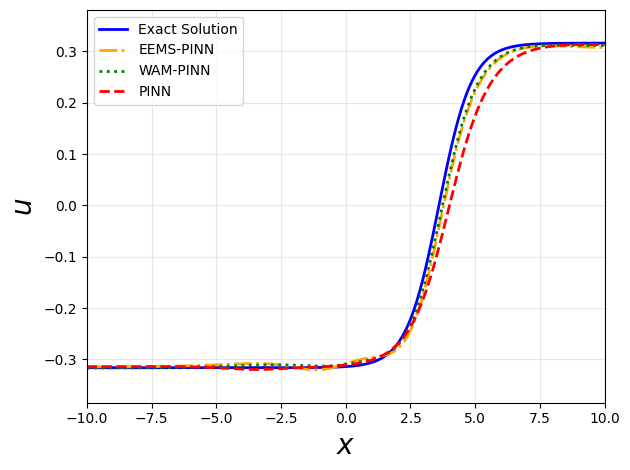}   
		\caption*{$t=12$}  
	\end{subfigure}  
	% 为整个图像阵列添加标题和标签  
	\caption{The exact solution and numerical solutions of PINN, RAM-PINN, WAM-PINN and EEMS-PINN of Klein-Gordon equation (\ref{equ:KG}) at time $t=0$, $t=4$, $t=8$ and $t=12$, respectively.} 
    % 整个图像阵列的标题  
	\label{fig:Klein_sol_time}  
	
\end{figure}

\begin{figure} [!htp]  
	\centering  
	
	% 第一行
	\begin{subfigure}{.32\textwidth}  
		\centering  
		\includegraphics[width=1.00\linewidth]{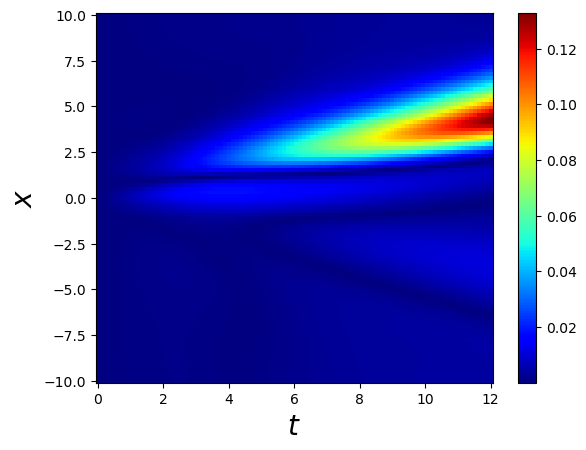}  
		%\caption*{PINN}  
		%\label{fig:14-1}  
	\end{subfigure}%  
	\begin{subfigure}{.32\textwidth}  
		\centering  
		\includegraphics[width=1.00\linewidth]{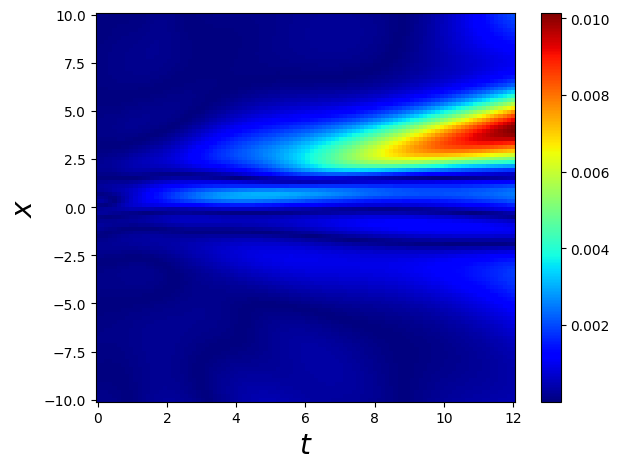}  
		%\caption*{RAR-PINN}  
		%\label{fig:14-2}  
	\end{subfigure}%  
        \begin{subfigure}{.32\textwidth}  
		\centering  
		\includegraphics[width=1.00\linewidth]{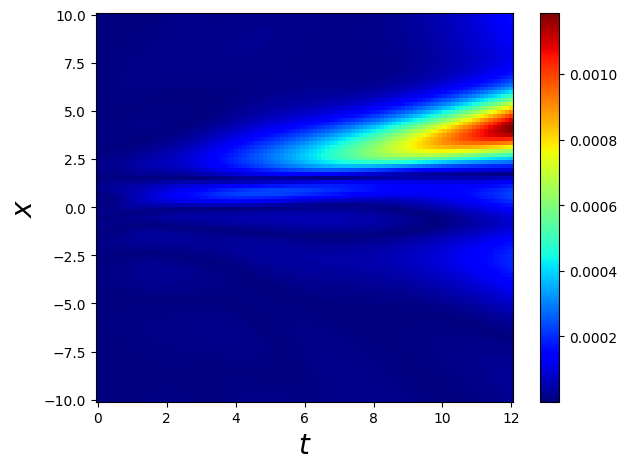}  
		%\caption*{RAR-PIRN}  
		%\label{fig:14-4}  
	\end{subfigure}%

	% 第二行
	\begin{subfigure}{.32\textwidth}  
		\centering  
		\includegraphics[width=1.00\linewidth]{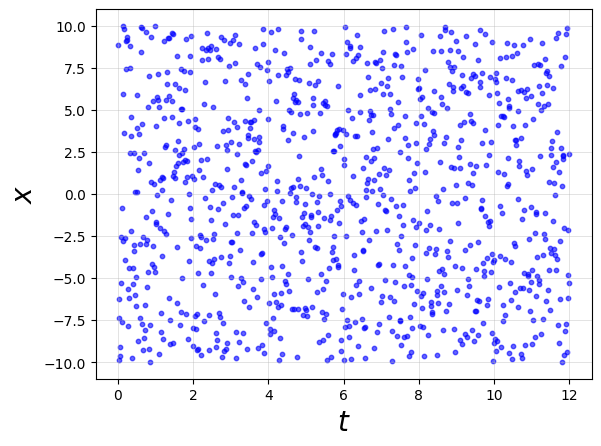}  
		\caption*{PINN}  
	\end{subfigure}%  
        \begin{subfigure}{.32\textwidth}  
		\centering  
		\includegraphics[width=1.00\linewidth]{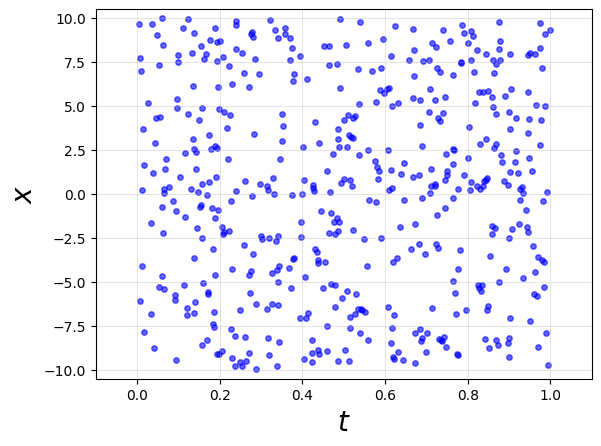}  
		\caption*{WAM-PINN}  
	\end{subfigure}%
	\begin{subfigure}{.32\textwidth}  
		\centering  
		\includegraphics[width=1.00\linewidth]{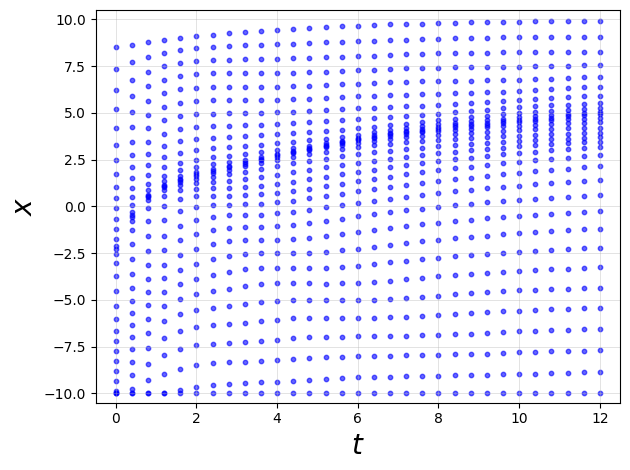}  
		\caption*{EEMS-PINN}  
	\end{subfigure}%
		% 为整个图像阵列添加标题和标签  
	\caption{(first row) The absolute errors of the estimated solutions; (second row) The distribution of sampling points for Klein-Gordon equation (\ref{equ:KG}) after one round mesh moving.}
	\label{fig:Klein_points} 
\end{figure}

Figure~\ref{fig:Klein_sol} presents a direct comparison between the exact solution and numerical approximations obtained from PINN, WAM-PINN, and EEMS-PINN. Under identical training conditions, EEMS-PINN demonstrates superior accuracy, as evidenced by the first row of Figure~\ref{fig:Klein_points} showing its significantly reduced absolute error compared to both PINN and WAM-PINN implementations. Quantitatively, the maximum absolute errors reveal stark differences: WAM-PINN achieves errors of $\mathcal{O}(10^{-2})$, while EEMS-PINN further reduces this to $\mathcal{O}(10^{-3})$---a full order of magnitude improvement over PINN's $\mathcal{O}(10^{-1})$ errors. The temporal evolution of solution accuracy, captured in Figure~\ref{fig:Klein_sol_time}, confirms that EEMS-PINN maintains this precision advantage throughout the entire simulation duration, including at the final time $t=12$, where it continues to provide the most faithful approximation to the analytic solution.

%The second row of Figure  \ref{fig:Klein_points} illustrates the collocation point distributions for all methods. Notably, the movable collocation points $\mathcal{X}_{p}$ in WAM-PINN exhibit a nearly uniform distribution. In contrast, the distribution of $\mathcal{X}_{p}$ for EEMS-PINN is densely clustered in a dense band between $x=0$ to $x=5$, mirroring the spatial coordinates of maximum solution curvature and high energy density. 

The second row of Figure~\ref{fig:Klein_points} provides a detailed comparison of the evolved collocation point distributions, revealing fundamental differences in the adaptive sampling behavior of each method. 
WAM-PINN's movable points show minimal spatial variation with only slight concentration near the final time. In contrast, EEMS-PINN demonstrates significantly more effective adaptation, with the points forming a concentrated band between $x=0$ and $x=5$. This localized clustering corresponds precisely to regions exhibiting both maximum solution curvature (as seen in Figure  \ref{fig:Klein_sol}) and peak energy density, confirming that our Energy-Equidistributed adaptation successfully targets the most dynamically important regions of the solution space.

\begin{figure} [!htp]  
\centering  
\begin{subfigure}{.45\textwidth}  
		\centering  
		\includegraphics[width=1.00\linewidth]{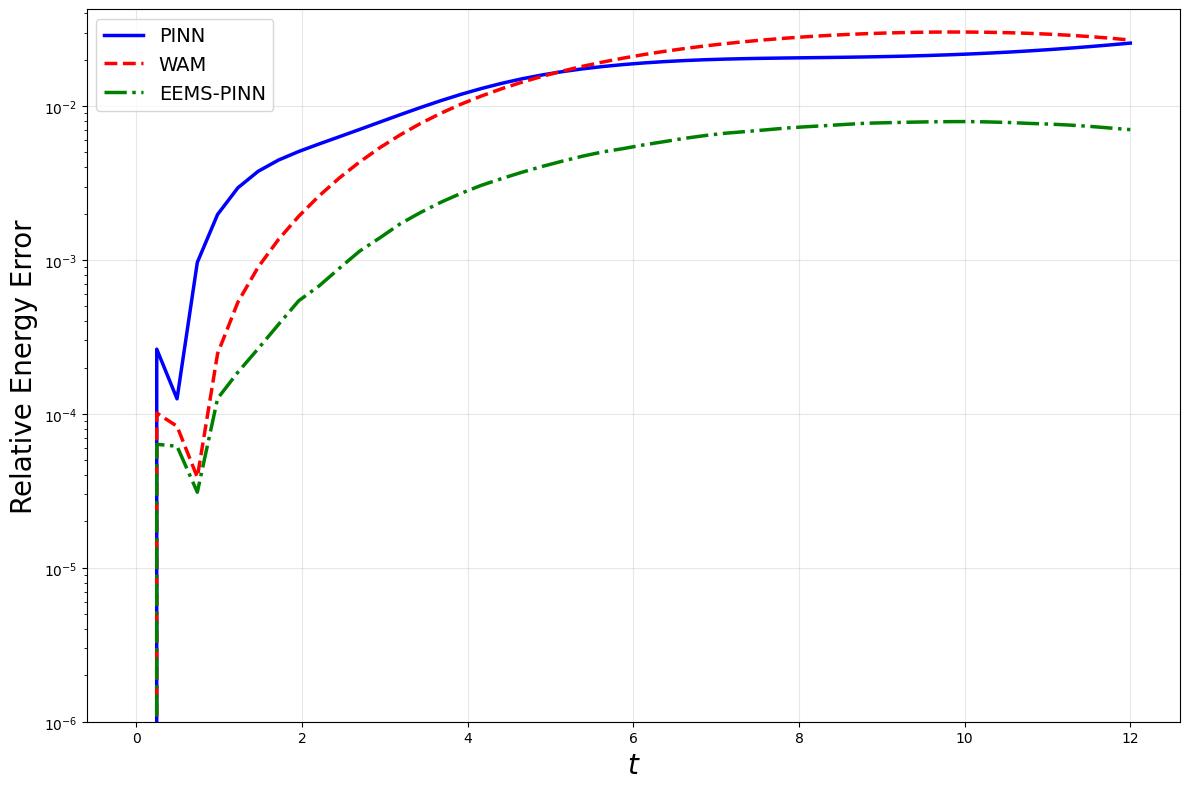}  
		\caption*{(a) Relative energy errors}  
\end{subfigure} 
 \begin{subfigure}{.45\textwidth}  
		\centering  
		\includegraphics[width=1.00\linewidth]{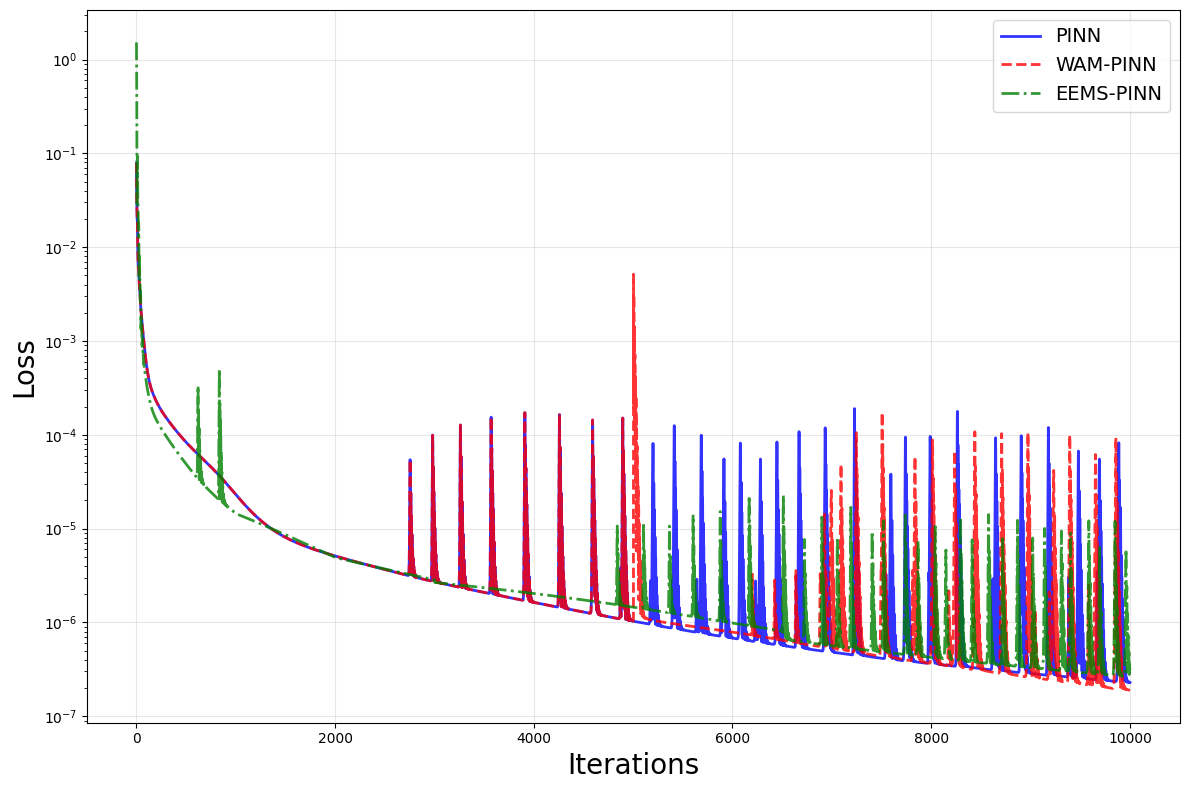}  
		\caption*{(b) Convergence curves of the  loss function}  
\end{subfigure}
	\caption{(a) The relative energy errors (in $\log$ scale); (b) the PDE loss convergence for Klein-Gordon equation (\ref{equ:KG}).} % 整个图像阵列的标题  
	\label{fig:Klein_convergence} % 整个图像阵列的标签
\end{figure}

The relative energy errors and PDE loss convergence behavior are compared across methods in Figure~\ref{fig:Klein_convergence}.  Panel (a) demonstrates that EEMS-PINN maintains the energy conservation with relative errors consistently of order $10^{-3}$, outperforming PINN and WAM-PINN. While WAM-PINN shows intermediate accuracy during the initial phase ($t\in[0,4]$), its energy conservation deteriorates significantly in later stages ($t\in[6,12]$), ultimately exhibiting the largest errors among all methods. Panel (b) reveals the training dynamics, where EEMS-PINN achieves both faster convergence and enhanced stability compared to alternatives---its PDE loss decays exponentially during initial iterations and maintains lower variance throughout optimization. 

The numerical performance comparison in Table~\ref{Tab:Klein} demonstrates EEMS-PINN's superior accuracy in solving the Klein-Gordon equation (\ref{equ:KG}) across different sampling densities. 
EEMS-PINN achieves better accuracy ($4.24\!\times\!10^{-3}$ at $N\!=\!1000$) compared to PINN and WAM-PINN, maintaining errors below $4\!\times\!10^{-3}$ for larger $N$, demonstrating efficient error saturation.

%Figure  \ref{fig:Klein_convergence}(a) shows the relative energy errors of all the numerical methods. Compared to PINN and WAM-PINN, EEMS-PINN keeps the lowest order of 0.01 over the whole time domain. WAM is lower than PINN around the initial time domain $[0, 4]$, but is the highest from $t=6$ to $t=12$. The PDE loss function of all the methods in the final training process after one round sampling moving in shown in Figure  \ref{fig:Klein_convergence}(b).  We can see that the loss of EEMS-PINN  decaying faster compared to the other two methods and EEMS-PINN is more stable during the whole training process. 

%We further investigate the impact of collocation point allocation by conducting experiments with varying numbers of sampling points in Table (\ref{Tab:Klein}). EEMS-PINN can provide a good approximation with the number of samples $N=1000$. 

% 导入图像阵列  

\begin{table}[!htbp]
    \centering
\begin{tabular}{@{}lccc@{}}
\toprule
$N$ & 1000 & 10000 & 20000 \\ \midrule
PINN & $5.22\times 10^{-1}$& $1.76\times 10^{-2}$ & $9.83\times 10^{-3}$ \\
WAM-PINN & $5.89\times 10^{-2}$ & $4.89\times 10^{-3}$ & $4.55\times 10^{-3}$ \\
EEMS-PINN & $4.24\times 10^{-3}$ & $3.88\times 10^{-3}$ & $3.81\times 10^{-3} $\\
\bottomrule
\end{tabular}
 \caption{The relative $L_2$ errors of Klein-Gordon equation (\ref{equ:KG}).}
\label{Tab:Klein}
\end{table}

\hspace*{\fill}

\noindent {\bf Example 2 (Klein-Gordon equation with force term)}  

\begin{figure} [!htp]  
\centering  
\begin{subfigure}{.45\textwidth}  
		\centering  
		\includegraphics[width=1.00\linewidth]{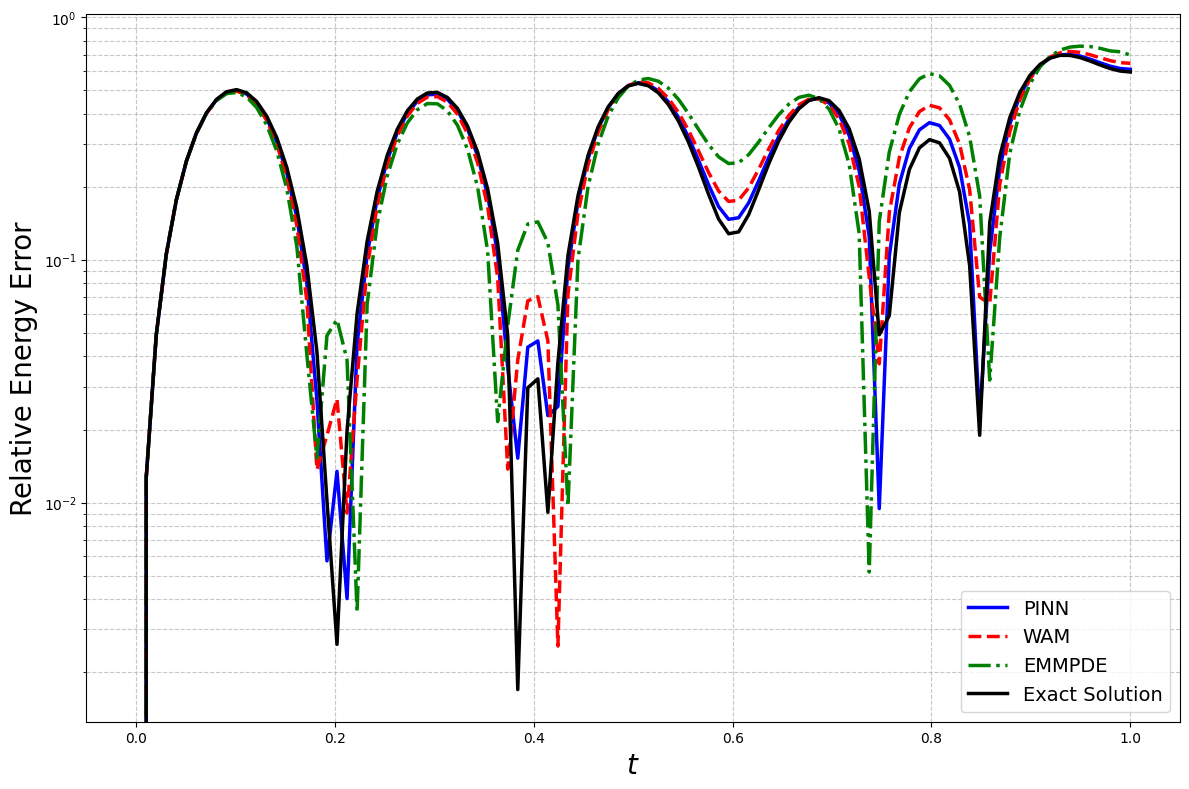}  
		\caption*{(a) Relative energy errors}  
\end{subfigure} 
 \begin{subfigure}{.45\textwidth}  
		\centering  
		\includegraphics[width=1.00\linewidth]{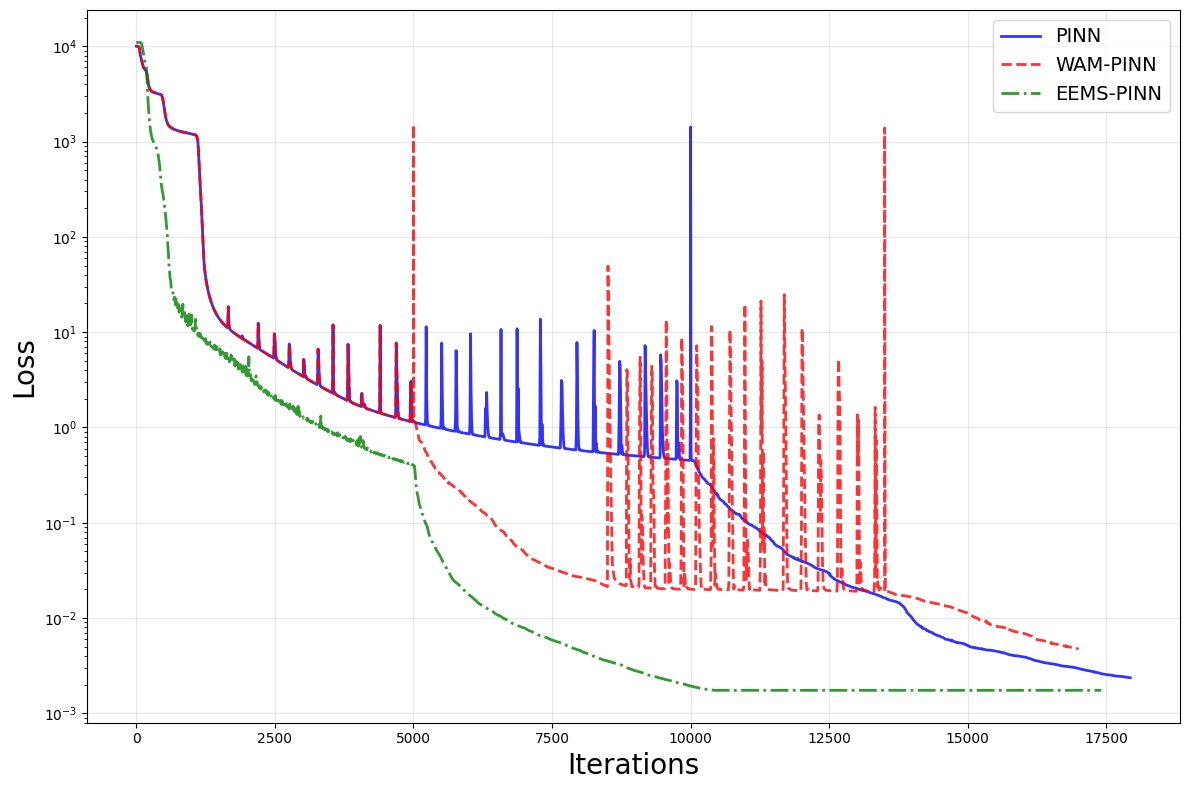}  
		\caption*{(b) the PDE loss convergence}  
\end{subfigure}
	\caption{(a) The relative energy errors (in $\log$ scale); (b) the PDE loss convergence for Klein-Gordon equation (\ref{equ:KG_source}).} % 整个图像阵列的标题  
	\label{fig:Klein_source_convergence} % 整个图像阵列的标签
\end{figure}

%The klein-Gordon equation is a rich mathematical structure that encapsulates various physical phenomena.
We test EEMS-PINN with Klein-Gordon equation with non-zero forcing term of the following form 
\begin{align}\label{equ:KG_source}
\begin{split}
&u_{t t}- u_{xx}+u^3=f(x, t), \ (x,t) \in [0,1] \times [0,1], \\
&u(x, 0)=h_1(x), \ u_t(x, 0)=h_2(x), \\
&u(0, t)=h_3(t),\  u(1, t)=h_4(t). \\
\end{split}
\end{align}
The test solution function is given by
\begin{equation*}\label{solutionHelmholt}
    u(x,t) = x\cos(5\pi t) + (xt)^3.
\end{equation*}
So the external force function $f(x,t)$, and initial/Dirichlet boundary conditions $h_{1}(x)$, $h_{2}(x)$, $h_{3}(t)$ and $h_{4}(t)$ are derived accordingly. 
While the system's total energy is non-conservative due to the external forcing (as confirmed by Figure~\ref{fig:Klein_source_convergence}(a)), we maintain the zero-forcing energy (\ref{equ:Klein_energy}) to validate our method's robustness.  

\begin{figure}[htbp]  
	\centering
	
	% 创建第一行第一列的子图  
	\begin{subfigure}{.24\textwidth}  
		\centering  
		\includegraphics[width=1.00\linewidth]{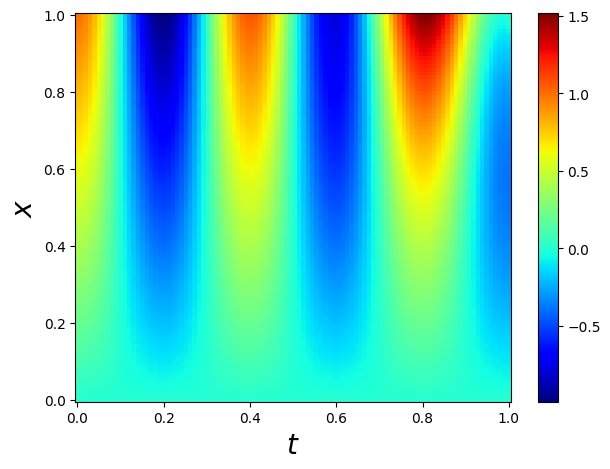}  
		\caption*{Exact}   
	\end{subfigure}%  
	% 创建第一行第二列的子图  
	\begin{subfigure}{.24\textwidth}  
		\centering  
		\includegraphics[width=1.00\linewidth]{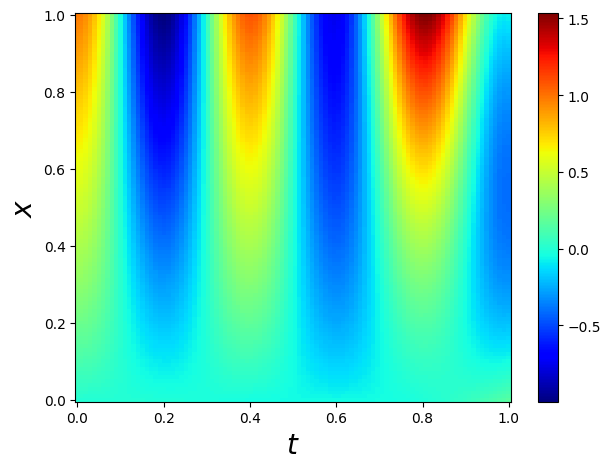}  
		\caption*{PINN}  
	\end{subfigure}  
	% 创建第一行第三列的子图  
	\begin{subfigure}{.24\textwidth}  
		\centering  
		\includegraphics[width=1.00\linewidth]{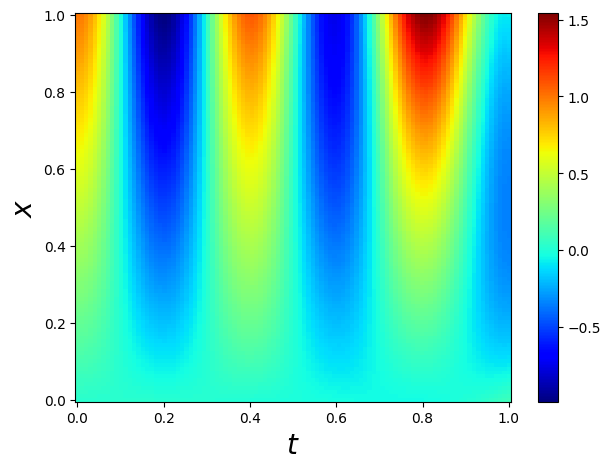} 
		\caption*{WAM-PINN}   
	\end{subfigure}  
		\begin{subfigure}{.24\textwidth}  
		\centering  
		\includegraphics[width=1.00\linewidth]{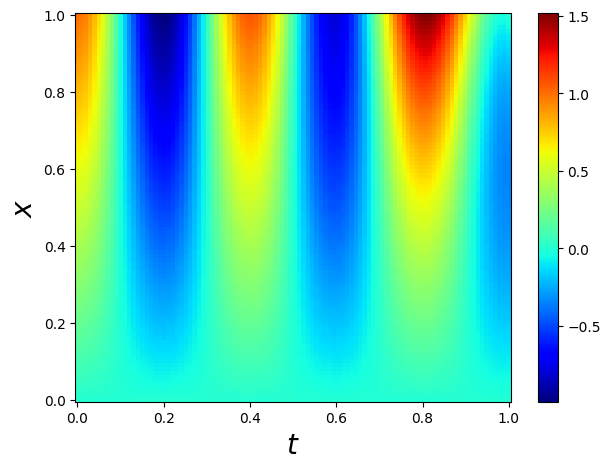}  
		\caption*{EEMS-PINN}  
	\end{subfigure} 
	% 为整个图像阵列添加标题和标签  
	\caption{The exact solution and numerical solutions of PINN, WAM-PINN and EEMS-PINN of  Klein-Gordon equation (\ref{equ:KG_source}), respectively.} 
    % 整个图像阵列的标题  
	\label{fig:Klein_source_sol} % 整个图像阵列的标签  
	
\end{figure}

\begin{figure}[htbp]  
	\centering
	
	% 创建第一行第一列的子图  
	\begin{subfigure}{.24\textwidth}  
		\centering  
		\includegraphics[width=1.00\linewidth]{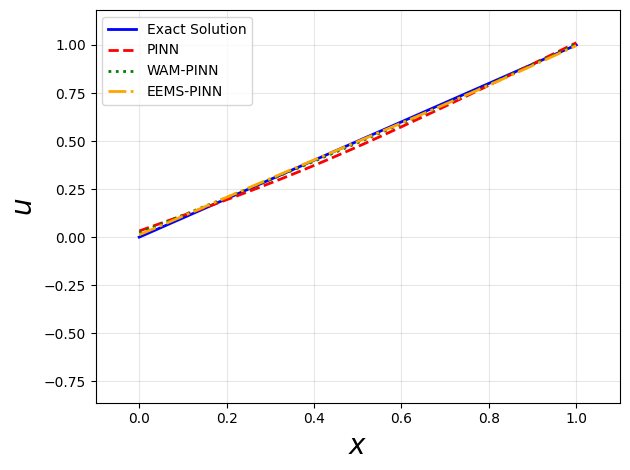}  
		\caption*{$t=0$}   
	\end{subfigure}%  
	% 创建第一行第二列的子图  
	\begin{subfigure}{.24\textwidth}  
		\centering  
		\includegraphics[width=1.00\linewidth]{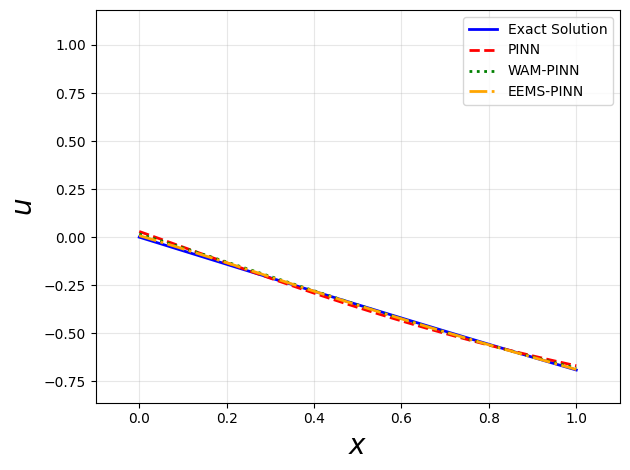}  
		\caption*{$t=0.25$}  
	\end{subfigure}  
	% 创建第一行第三列的子图  
	\begin{subfigure}{.24\textwidth}  
		\centering  
		\includegraphics[width=1.00\linewidth]{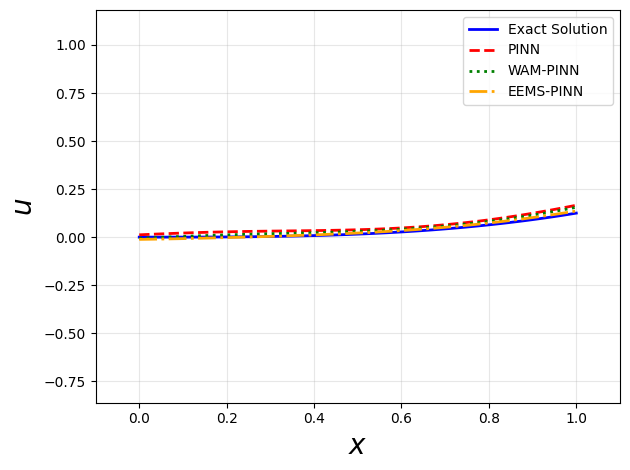}    
		\caption*{$t=0.5$}   
	\end{subfigure}  
        \begin{subfigure}{.24\textwidth}  
		\centering  
		\includegraphics[width=1.00\linewidth]{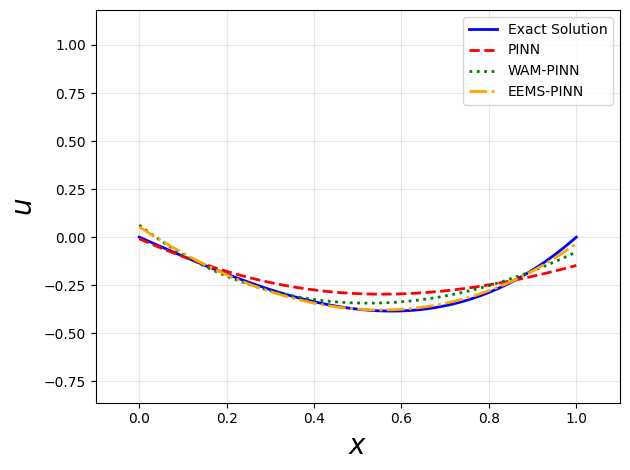}   
		\caption*{$t=1$}  
	\end{subfigure}  
	% 为整个图像阵列添加标题和标签  
	\caption{The exact solution and numerical solutions of PINN, WAM-PINN and EEMS-PINN of equation (\ref{equ:KG_source}) at time $t=0$, $t=0.25$, $t=0.5$ and $t=1$.} 
    % 整个图像阵列的标题  
	\label{fig:Klein_source_sol_time}  
\end{figure}

\begin{figure} [!htp]  
	\centering  
	
	% 第一行
	\begin{subfigure}{.32\textwidth}  
		\centering  
		\includegraphics[width=1.00\linewidth]{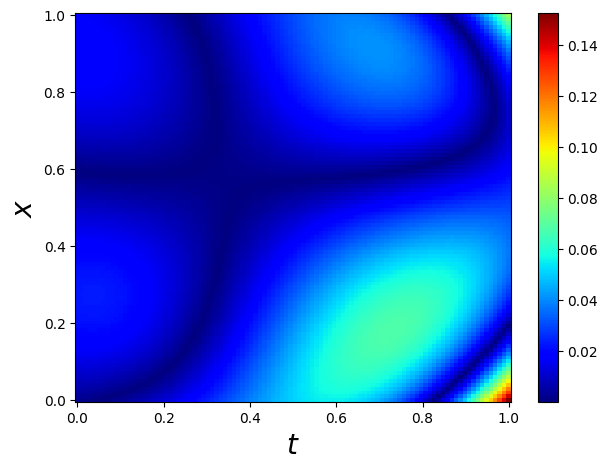}  
		%\caption*{PINN}  
	\end{subfigure}%  
	\begin{subfigure}{.32\textwidth}  
		\centering  
		\includegraphics[width=1.00\linewidth]{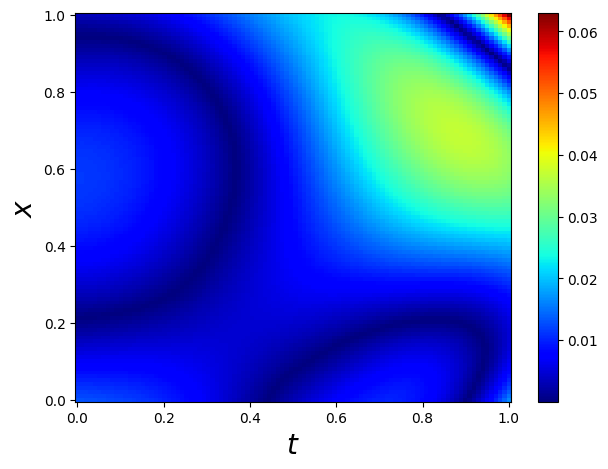}  
		%\caption*{RAR-PINN}    
	\end{subfigure}%  
		\begin{subfigure}{.32\textwidth}  
		\centering  
		\includegraphics[width=1.00\linewidth]{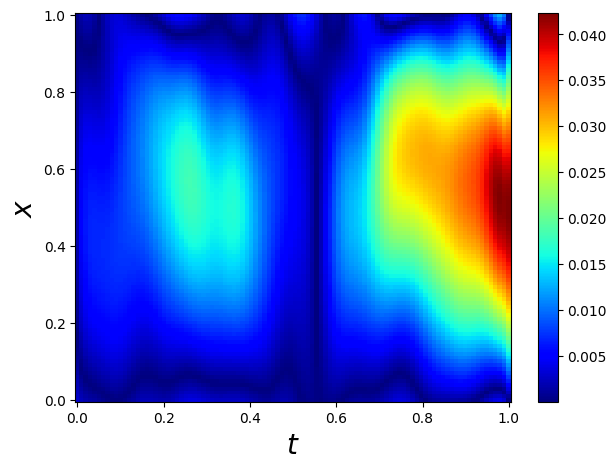}  
		%\caption*{RAR-PIRN}   
	\end{subfigure}%
    
	% 第二行
	\begin{subfigure}{.32\textwidth}  
		\centering  
		\includegraphics[width=1.00\linewidth]{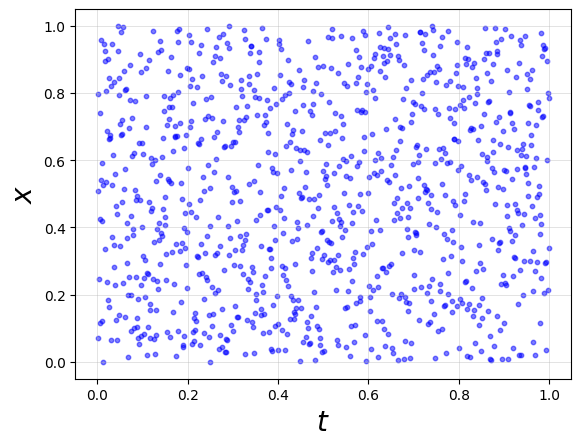}  
		\caption*{PINN}  
	\end{subfigure}%  
        \begin{subfigure}{.32\textwidth}  
		\centering  
		\includegraphics[width=1.00\linewidth]{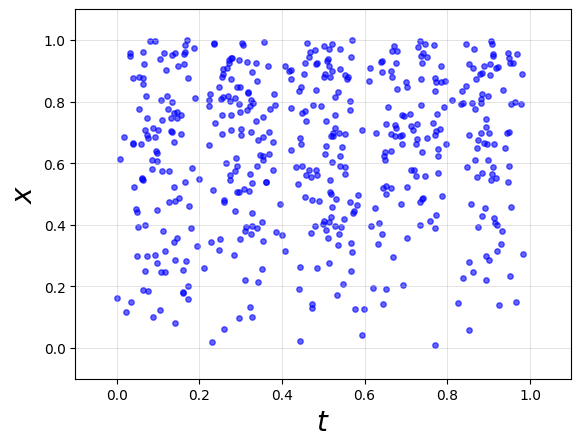}  
		\caption*{WAM-PINN}  
	\end{subfigure}%
	\begin{subfigure}{.32\textwidth}  
		\centering  
		\includegraphics[width=1.00\linewidth]{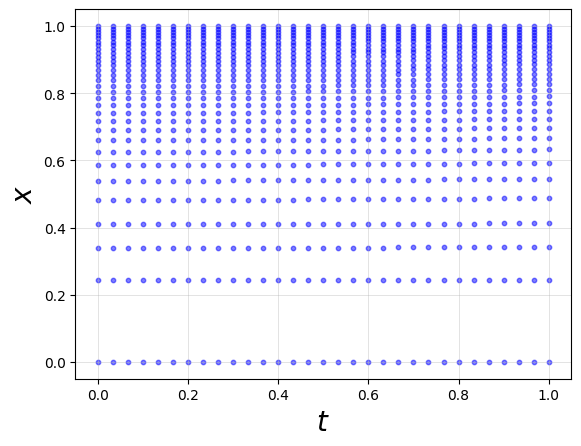}  
		\caption*{EEMS-PINN}  
	\end{subfigure}%
		% 为整个图像阵列添加标题和标签  
	\caption{ (first row) The absolute errors of the estimated solutions; (second row) The distribution of sampling points for Klein-Gordon equation (\ref{equ:KG_source}) after one round mesh moving.}
	\label{fig:Klein_source_points} % 整个图像阵列的标签
\end{figure}

Figures~\ref{fig:Klein_source_sol}--\ref{fig:Klein_source_convergence} present a comprehensive comparison of solution accuracy and point distributions for the forced Klein-Gordon equation using $N_p=1000$ collocation points with $N_i=100$ and $N_b=100$ initial and boundary points. The results demonstrate EEMS-PINN's strong agreement with the exact solution, particularly in regions of  high solution gradient magnitude. The point distribution analysis in Figure~\ref{fig:Klein_source_points} reveals that both WAM-PINN and EEMS-PINN concentrate the points within the boundary layer region $x\in[0.8,1]$. This adaptive clustering corresponds precisely to the zone of peak energy density, demonstrating the method's ability to automatically allocate computational resources to critical regions without manual intervention.

Figure~\ref{fig:Klein_source_convergence}(b) demonstrates the rapid convergence of EEMS-PINN's PDE loss compared to other methods. The numerical results in Table~\ref{Tab:Klein_source} quantitatively confirm EEMS-PINN's higher precision  in solving the forced Klein-Gordon equation (\ref{equ:KG_source}), achieving a relative $L_2$ error of $\mathcal{O}(10^{-4})$ with $N=20000$ sample points---an order of magnitude better than WAM-PINN and conventional PINN. Remarkably, EEMS-PINN maintains this precision advantage across all tested sampling densities, which confirms the method's robustness regardless of the system's non-conservative energy properties.

\begin{comment}
{\color{red}As mentioned in Section \ref{sec:intro}, now the total energy is not conservative, as it is shown in Figure~\ref{fig:Klein_source_convergence}(a). Here we still use the energy density function as the monitor function to demonstrate the
robustness of our methodology across different physical scenarios.

Through 
Figure~\ref{fig:Klein_source_sol} to Figure~\ref{fig:Klein_source_convergence}, we use the sampling number $N_p=1000$, $N_i=100$ and $N_b=100$ to 
present a direct comparison between the exact solution and numerical approximations obtained from PINN, WAM-PINN and EEMS-PINN. Similar to the Klein-Gordon equation with zero-forcing (\ref{equ:KG}) as shown in the second row of Figure  \ref{fig:Klein_source_points}, the adaptive sampling points of WAM-PINN are basically
uniform randomly distributed, while those of EEMS-PINN are densely clustered near the domain boundary  $x=1$. This is corresponding to the high energy density. 

From the Figure  \ref{fig:Klein_convergence}(b), we can see that the PDE loss of EEMS-PINN dedays fastest. Furthermore, the numerical performance comparison in Table~\ref{Tab:Klein_source} demonstrates EEMS-PINN's superior accuracy in solving the Klein-Gordon equation with force term (\ref{equ:KG_source}) across different sampling densities even though we neglect the fact that the energy is not conservative. 
EEMS-PINN achieves superior accuracy ($3.64\!\times\!10^{-3}$ at $N\!=\!1000$) compared to PINN and WAM-PINN, maintaining errors below $4\!\times\!10^{-4}$ for larger $N$.
}
\end{comment}

% 导入图像阵列  

\begin{table}[!htbp]
    \centering
\begin{tabular}{@{}lccc@{}}
\toprule
$N$ & 1000 & 10000 & 20000 \\ \midrule
PINN & $2.88\times 10^{-2}$& $4.86\times 10^{-3}$ & $4.88\times 10^{-3}$ \\
WAM-PINN & $4.21\times 10^{-3}$ & $2.58\times 10^{-3}$ & $1.30\times 10^{-3}$ \\
EEMS-PINN & $3.64\times 10^{-3}$ & $9.83\times 10^{-4}$ & $8.67\times 10^{-4} $\\
\bottomrule
\end{tabular}
 \caption{The relative $L_2$ errors of Klein-Gordon equation (\ref{equ:KG_source}).}
\label{Tab:Klein_source}
\end{table}

\hspace*{\fill}

\noindent {\bf Example 3 (Sine-Gordon equation)}  

We consider the Sine-Gordon equation with a kink-antikink solution \cite{pekmen2012differential}
\begin{align}\label{eq:NonLinearSineGordonEq}
\begin{split}
    & u_{tt} - u_{xx} + \sin(u) = 0, ~~(x,t)  \in \mathbb{R} \times [0,T],\\
        &u(x,0)=u_0(x),~~u_t(x,0)=u_1(x),
\end{split}   
\end{align}
with the periodic boundary conditions. The exact solution is given by 

\begin{equation*}
    u(x,t)=4\tan^{-1}\Big(\frac{\sinh\big(\frac{\gamma t}{\sqrt{1-\gamma^2}}\big)}{\gamma\cosh\big(\frac{x}{\sqrt{1-\gamma^2}}\big)}\Big).
\end{equation*}
It represents two solitons moving in opposite directions with speed $\gamma \in (0,1)$. 
The initial conditions are computed using the exact solution at $t=0$. We can obtain the energy functional of the Sine-Gordon equation (\ref{eq:NonLinearSineGordonEq}) as follows:
\begin{equation*}
H(t)
=\frac{1}{2} \int_{} \left[ u_t^2 + u_x^2 - 2(1-\cos u)\right] {\rm d} {x}.
\end{equation*}
In this experiment, we
set $\gamma= 0.5$ and solve the Sine-Gordon equation (\ref{eq:NonLinearSineGordonEq}) by the proposed method till $T = 1$. 

\begin{figure}[htbp]  
	\centering
	
	% 创建第一行第一列的子图  
	\begin{subfigure}{.24\textwidth}  
		\centering  
		\includegraphics[width=1.00\linewidth]{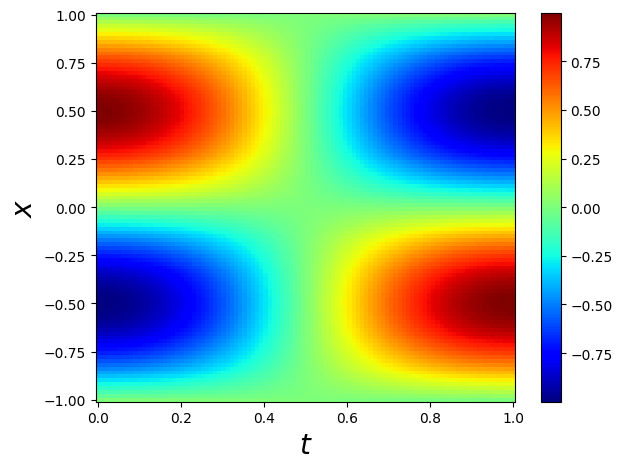}  
		\caption*{Exact}   
	\end{subfigure}%  
	% 创建第一行第二列的子图  
	\begin{subfigure}{.24\textwidth}  
		\centering  
		\includegraphics[width=1.00\linewidth]{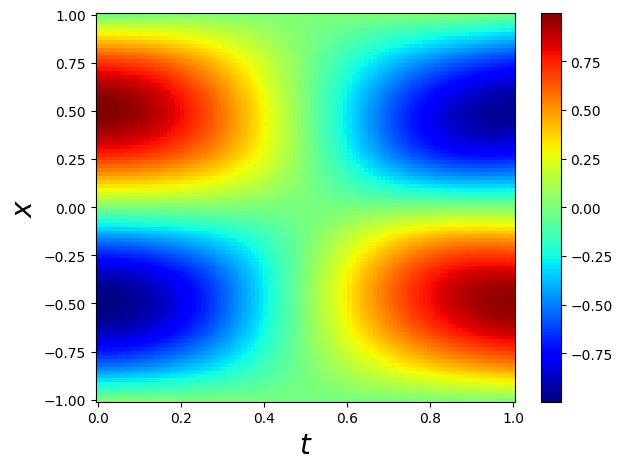}  
		\caption*{PINN}  
	\end{subfigure}  
	% 创建第一行第三列的子图  
	\begin{subfigure}{.24\textwidth}  
		\centering  
		\includegraphics[width=1.00\linewidth]{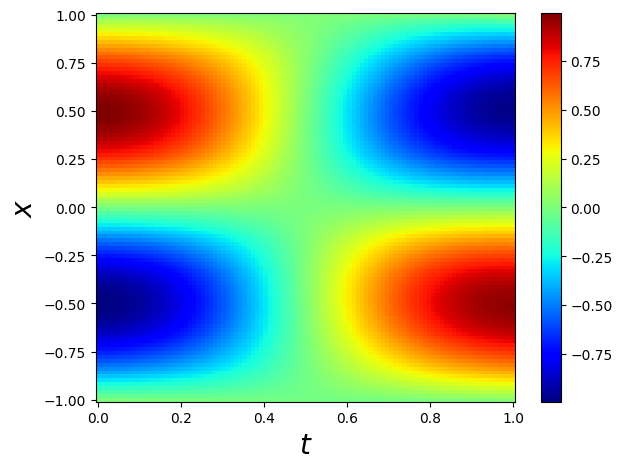}  
		\caption*{WAM-PINN}   
	\end{subfigure}  
		\begin{subfigure}{.24\textwidth}  
		\centering  
		\includegraphics[width=1.00\linewidth]{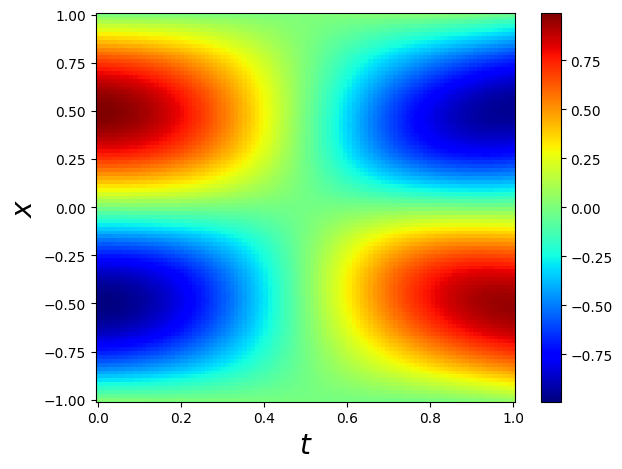}  
		\caption*{EEMS-PINN}  
	\end{subfigure} 
	% 为整个图像阵列添加标题和标签  
	\caption{The exact solution and numerical solutions of PINN, WAM-PINN and EEMS-PINN of Sine-Gordon equation (\ref{eq:NonLinearSineGordonEq}).} 
    % 整个图像阵列的标题  
	\label{fig:Sine_sol} % 整个图像阵列的标签  
	
\end{figure}

\begin{figure}[htbp]  
	\centering
	
	% 创建第一行第一列的子图  
	\begin{subfigure}{.24\textwidth}  
		\centering  
		\includegraphics[width=1.00\linewidth]{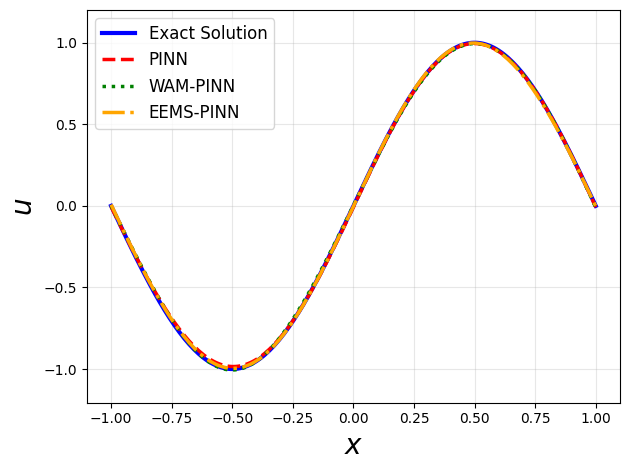}  
		\caption*{$t=0$}   
	\end{subfigure}%  
	% 创建第一行第二列的子图  
	\begin{subfigure}{.24\textwidth}  
		\centering  
		\includegraphics[width=1.00\linewidth]{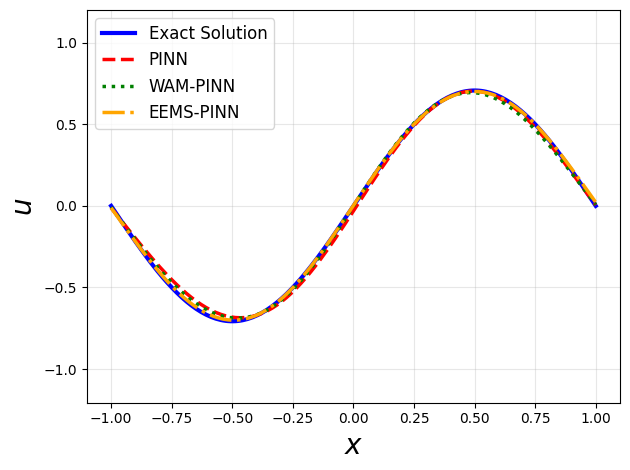}  
		\caption*{$t=0.25$}  
	\end{subfigure}  
	% 创建第一行第三列的子图  
	\begin{subfigure}{.24\textwidth}  
		\centering  
		\includegraphics[width=1.00\linewidth]{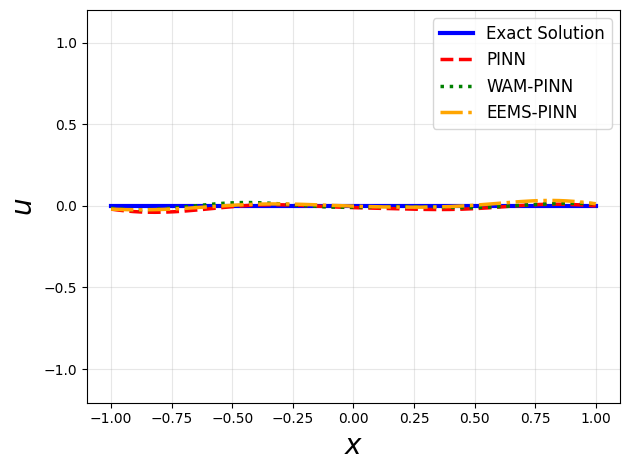}    
		\caption*{$t=0.5$}   
	\end{subfigure}  
        \begin{subfigure}{.24\textwidth}  
		\centering  
		\includegraphics[width=1.00\linewidth]{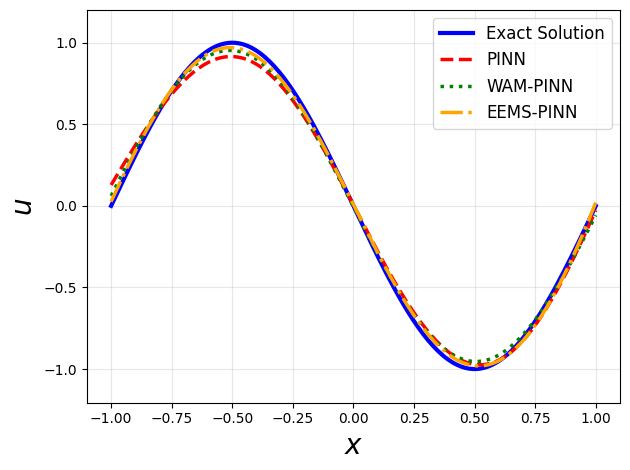}   
		\caption*{$t=1$}  
	\end{subfigure}  
	% 为整个图像阵列添加标题和标签  
	\caption{The exact solution and numerical solutions of PINN, WAM-PINN and EEMS-PINN of Sine-Gordon equation (\ref{eq:NonLinearSineGordonEq}) at time $t=0$, $t=0.25$, $t=0.5$ and $t=1$.} 
    % 整个图像阵列的标题  
	\label{fig:sine_sol_time}  
\end{figure}

\begin{figure} [!htp]  
	\centering  
	
	% 第一行
	\begin{subfigure}{.32\textwidth}  
		\centering  
		\includegraphics[width=1.00\linewidth]{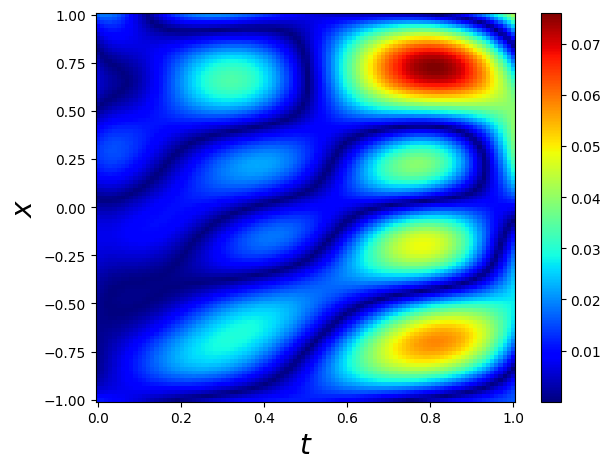}  
		%\caption*{PINN}  
	\end{subfigure}%  
	\begin{subfigure}{.32\textwidth}  
		\centering  
		\includegraphics[width=1.00\linewidth]{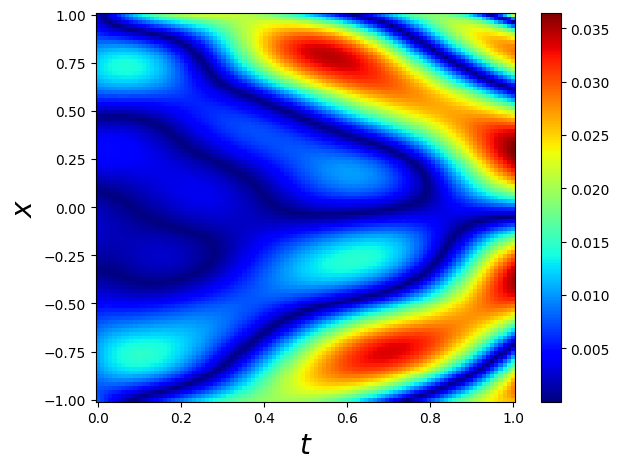}  
		%\caption*{WAM-PINN}    
	\end{subfigure}%  
		\begin{subfigure}{.32\textwidth}  
		\centering  
		\includegraphics[width=1.00\linewidth]{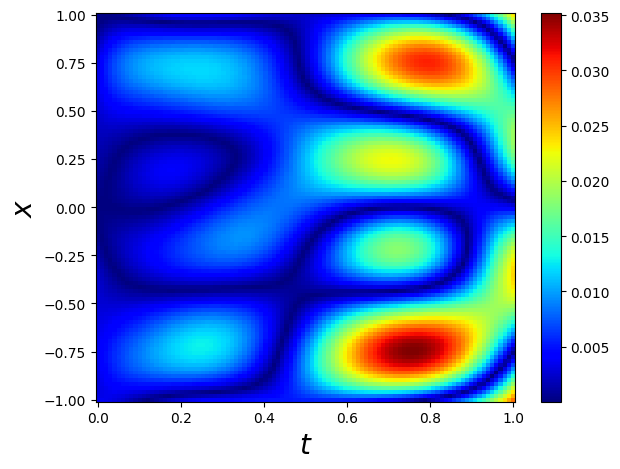}  
		%\caption*{EEMS-PINN}   
	\end{subfigure}%
    
	% 第二行
	\begin{subfigure}{.32\textwidth}  
		\centering  
		\includegraphics[width=1.00\linewidth]{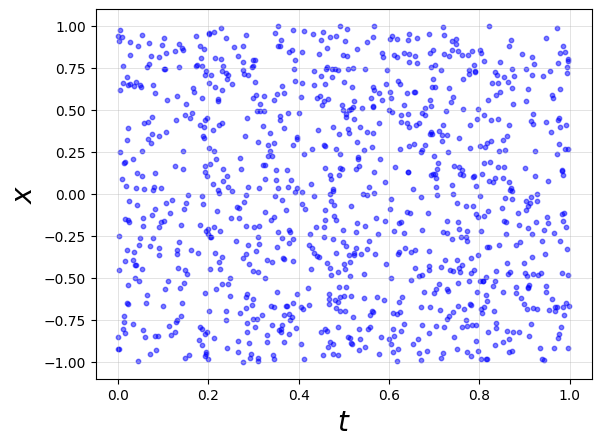}  
		\caption*{PINN}  
	\end{subfigure}%  
        \begin{subfigure}{.32\textwidth}  
		\centering  
		\includegraphics[width=1.00\linewidth]{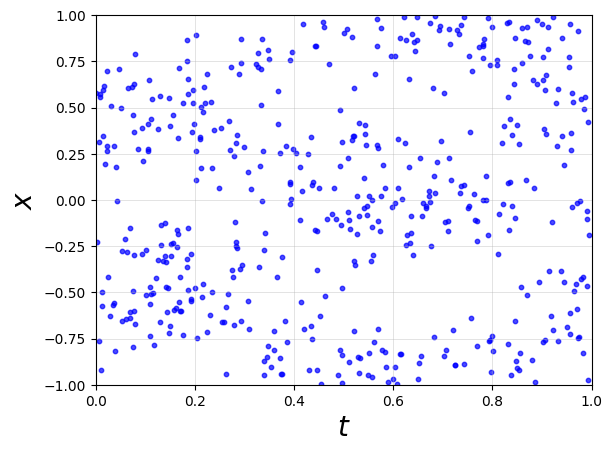}  
		\caption*{WAM-PINN}  
	\end{subfigure}%
	\begin{subfigure}{.32\textwidth}  
		\centering  
		\includegraphics[width=1.00\linewidth]{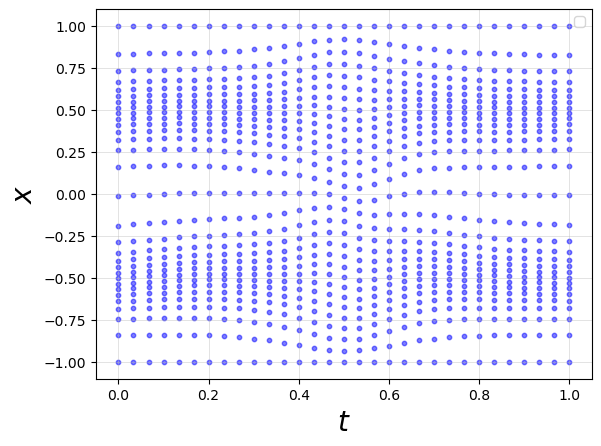}  
		\caption*{EEMS-PINN}  
	\end{subfigure}%
		% 为整个图像阵列添加标题和标签  
	\caption{(first row) The absolute errors of the estimated solutions; (second row) The distribution of sampling points for Sine-Gordon equation (\ref{eq:NonLinearSineGordonEq}) after one round mesh moving.}
	\label{fig:Sine_points} % 整个图像阵列的标签
\end{figure}

\begin{figure} [!htp]  
\centering  
\begin{subfigure}{.45\textwidth}  
		\centering  
		\includegraphics[width=1.00\linewidth]{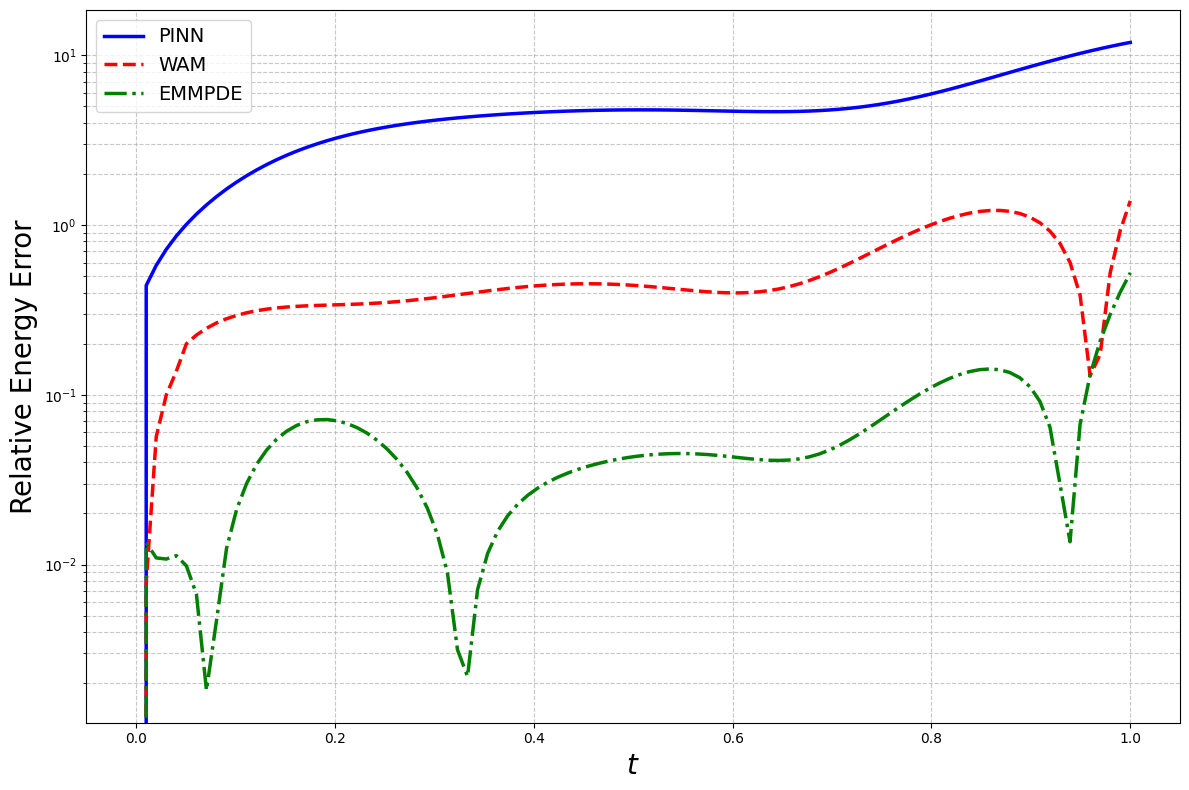}  
		\caption*{(a) Relative energy errors}  
\end{subfigure} 
 \begin{subfigure}{.45\textwidth}  
		\centering  
		\includegraphics[width=1.00\linewidth]{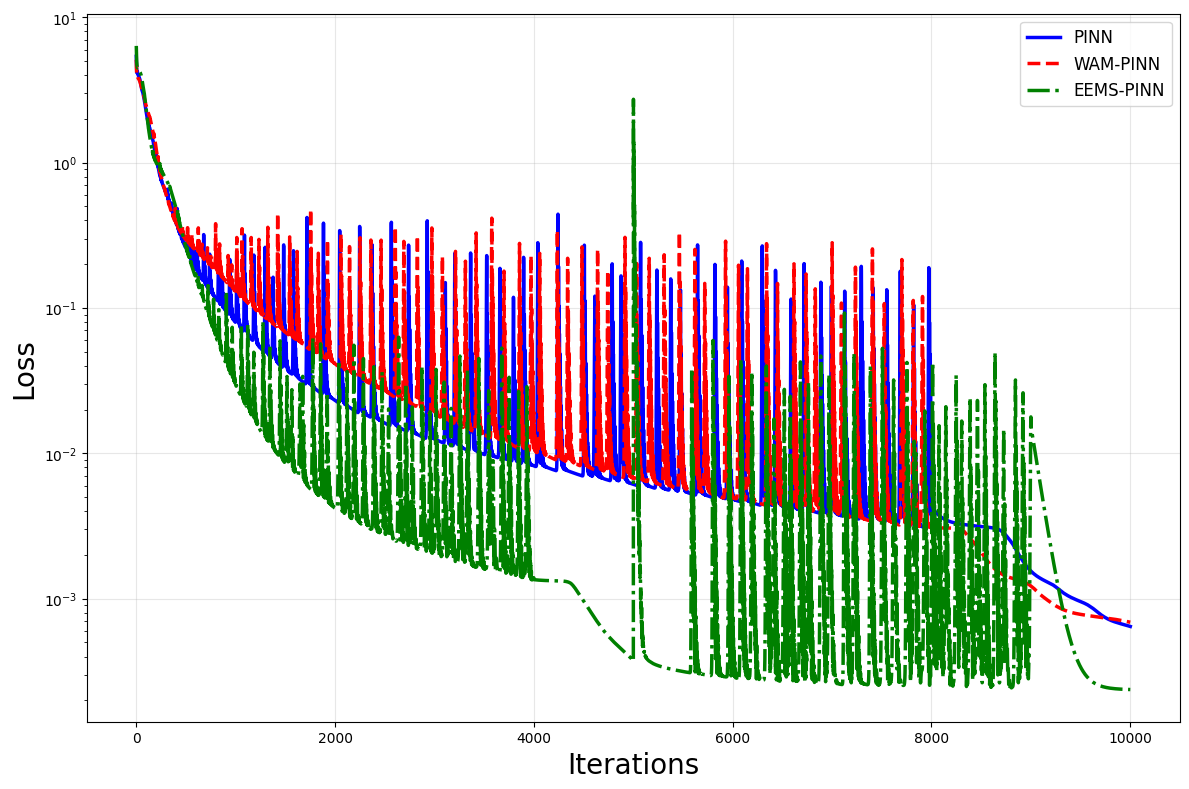}  
		\caption*{(b)  the PDE loss convergence}  
\end{subfigure}
	\caption{(a) The relative energy errors (in $\log$ scale); (b) the PDE loss convergence for Sine-Gordon equation (\ref{eq:NonLinearSineGordonEq}).}
    % 整个图像阵列的标题  
	\label{fig:sine_convergence} % 整个图像阵列的标签
\end{figure}

In Figures~\ref{fig:Sine_sol}--\ref{fig:sine_convergence}, we show a detailed comparison of solution accuracy and adaptive point distributions for the one-dimensional Sine-Gordon equation (\ref{eq:NonLinearSineGordonEq}), evaluated under identical training conditions ($N_p=1000$, $N_i=N_b=100$). 
The results show that EEMS-PINN achieves closer agreement with the exact solution compared to WAM-PINN, particularly in long-time simulations.  Figure~\ref{fig:Sine_points} demonstrates EEMS-PINN's enhanced ability to concentrate collocation points near soliton interactions. While Table~\ref{Tab:Sine} indicates comparable $L_2$ error magnitudes between the two adaptive methods (both $\mathcal{O}(10^{-3})$), EEMS-PINN achieves significantly better relative energy conservation and faster PDE loss convergence.

%The results demonstrate EEMS-PINN's superior agreement with the exact solution, particularly in the long time evolution. The point distribution analysis in Figure~\ref{fig:Sine_points} reveals that EEMS-PINN shows more clearly points concentration near the two solitons than  WAM-PINN. From Table \ref{Tab:Sine}, although EEMS-PINN show  slightly advantages of  $L_2$ errors (in the same order) compared to WAM-PINN, but the former still keeps the better relative energy errors and have faster PDE loss converges. 

\begin{table}[!htbp]
    \centering
\begin{tabular}{@{}lccc@{}}
\toprule
$N$ & 1000 & 10000 & 20000 \\ \midrule
PINN & $2.88\times 10^{-2}$& $5.86\times 10^{-3}$ & $4.53\times 10^{-3}$ \\
WAM-PINN & $8.77\times 10^{-3}$ & $4.32\times 10^{-3}$ & $2.87\times 10^{-3}$ \\
EEMS-PINN & $4.86\times 10^{-3}$ & $2.21\times 10^{-3}$ & $2.13\times 10^{-3} $\\
\bottomrule
\end{tabular}
 \caption{The relative $L_2$ errors of Sine-Gordon equation (\ref{eq:NonLinearSineGordonEq}).}
\label{Tab:Sine}
\end{table}

\hspace*{\fill}

\noindent {\bf Example 4 (KdV equation)}

Consider the following KdV equation:
\begin{equation}\label{equ:kdv}
 u_t  +6 u u_x+u_{xxx}=0,   ~~(x,t)  \in [0,2\pi] \times [0,T]
\end{equation}
with periodic boundary conditions and initial condition $u(x,0)=\sin x$.  
The energy functional is as follows:
$$
{H}(t) =\int\left[\frac{1}{2}\left(u_x\right)^2-u^3\right] {\rm d} x.  
$$

\begin{figure}[htbp]  
	\centering
	
	% 创建第一行第一列的子图  
	\begin{subfigure}{.24\textwidth}  
		\centering  
		\includegraphics[width=1.00\linewidth]{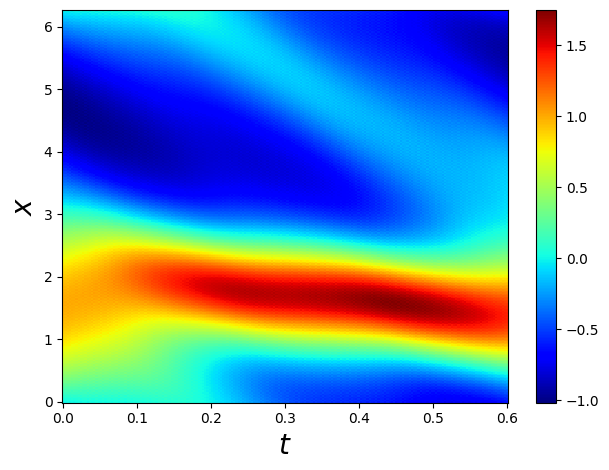}  
		\caption*{Exact}   
	\end{subfigure}%  
	% 创建第一行第二列的子图  
	\begin{subfigure}{.24\textwidth}  
		\centering  
		\includegraphics[width=1.00\linewidth]{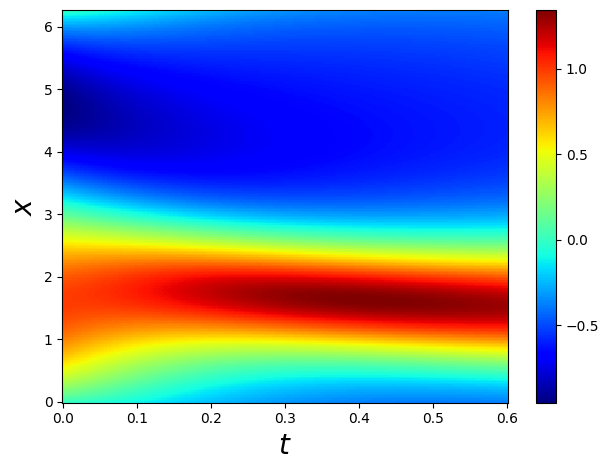}  
		\caption*{PINN}  
	\end{subfigure}  
	% 创建第一行第三列的子图  
	\begin{subfigure}{.24\textwidth}  
		\centering  
		\includegraphics[width=1.00\linewidth]{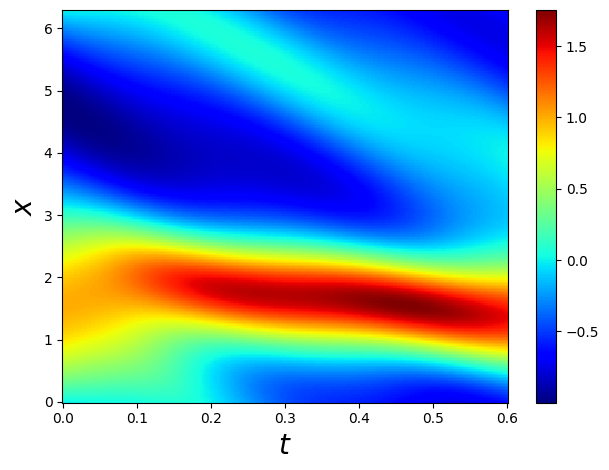}  
		\caption*{WAM-PINN}   
	\end{subfigure}  
		\begin{subfigure}{.24\textwidth}  
		\centering  
		\includegraphics[width=1.00\linewidth]{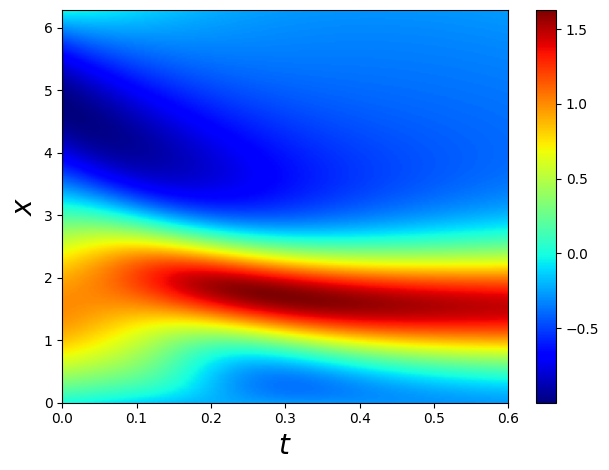}  
		\caption*{EEMS-PINN}  
	\end{subfigure} 
	% 为整个图像阵列添加标题和标签  
	\caption{The exact solution and numerical solutions of PINN, WAM-PINN and EEMS-PINN of KdV equation (\ref{equ:kdv}).} 
    \label{fig:kdv_sol}
    \end{figure}

\begin{figure}[htbp]  
	\centering
	
	% 创建第一行第一列的子图  
	\begin{subfigure}{.24\textwidth}  
		\centering  
		\includegraphics[width=1.00\linewidth]{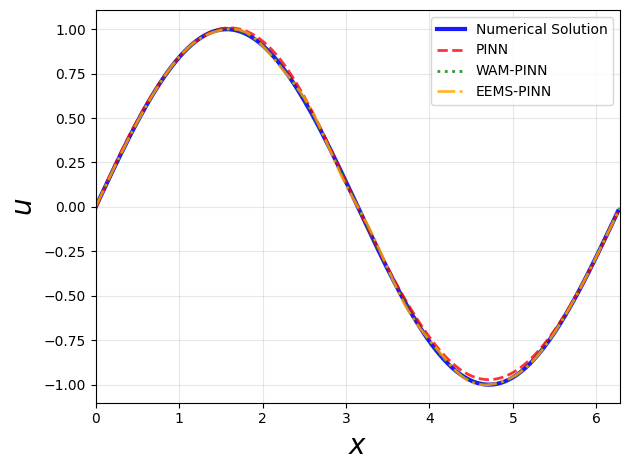}  
		\caption*{$t=0$}   
	\end{subfigure}%  
	% 创建第一行第二列的子图  
	\begin{subfigure}{.24\textwidth}  
		\centering  
		\includegraphics[width=1.00\linewidth]{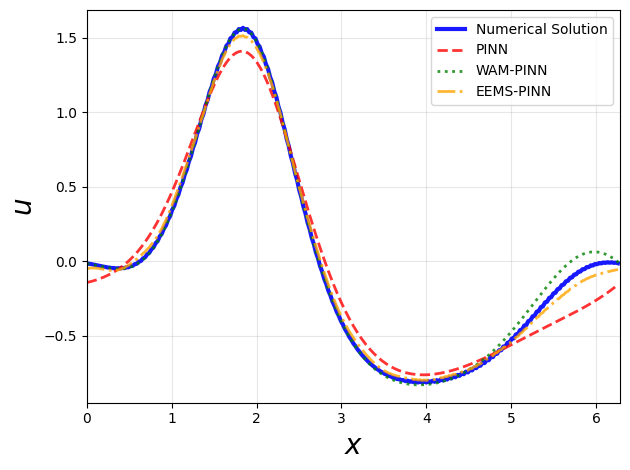}  
		\caption*{$t=0.2$}  
	\end{subfigure}  
	% 创建第一行第三列的子图  
	\begin{subfigure}{.24\textwidth}  
		\centering  
		\includegraphics[width=1.00\linewidth]{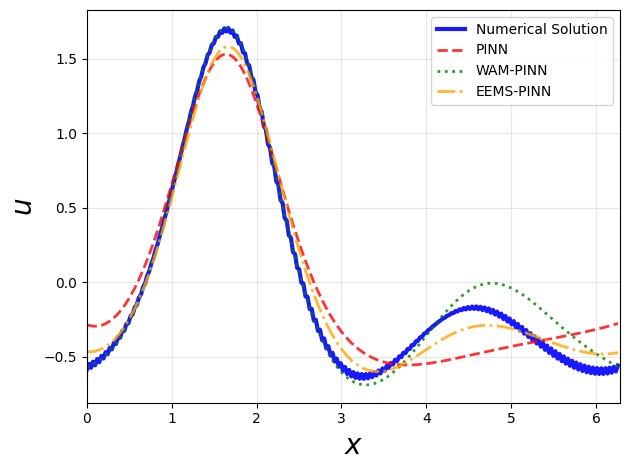}    
		\caption*{$t=0.4$}   
	\end{subfigure}  
        \begin{subfigure}{.24\textwidth}  
		\centering  
		\includegraphics[width=1.00\linewidth]{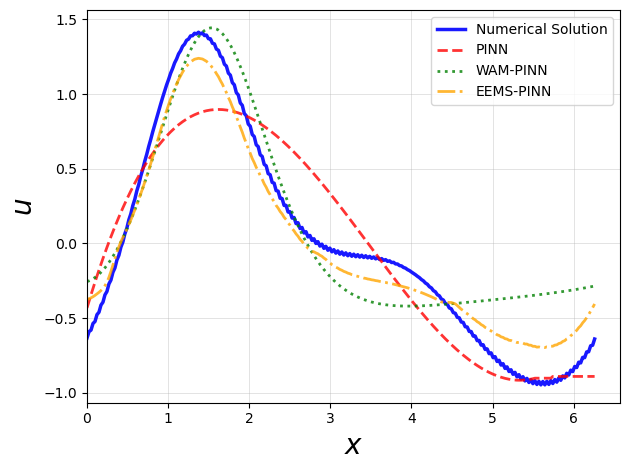}   
		\caption*{$t=0.6$}  
	\end{subfigure}  
	% 为整个图像阵列添加标题和标签  
	\caption{The exact solution and numerical solutions of PINN, WAM-PINN and EEMS-PINN of KdV equation (\ref{equ:kdv}) at time $t=0$, $t=0.2$, $t=0.4$ and $t=0.6$.} 
    % 整个图像阵列的标题  
	\label{fig:KdV_sol_time}  
\end{figure}

\begin{figure} [!htp]  
	\centering  
	
	% 第一行
	\begin{subfigure}{.32\textwidth}  
		\centering  
		\includegraphics[width=1.00\linewidth]{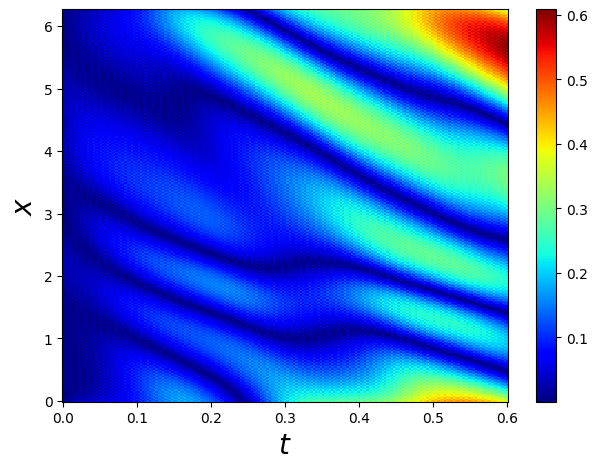}  
		%\caption*{PINN}  
	\end{subfigure}%  
	\begin{subfigure}{.32\textwidth}  
		\centering  
		\includegraphics[width=1.00\linewidth]{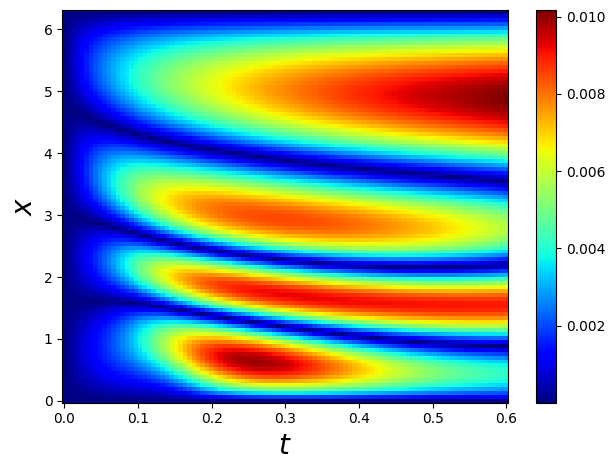}  
		%\caption*{WAM-PINN}    
	\end{subfigure}%  
		\begin{subfigure}{.32\textwidth}  
		\centering  
		\includegraphics[width=1.00\linewidth]{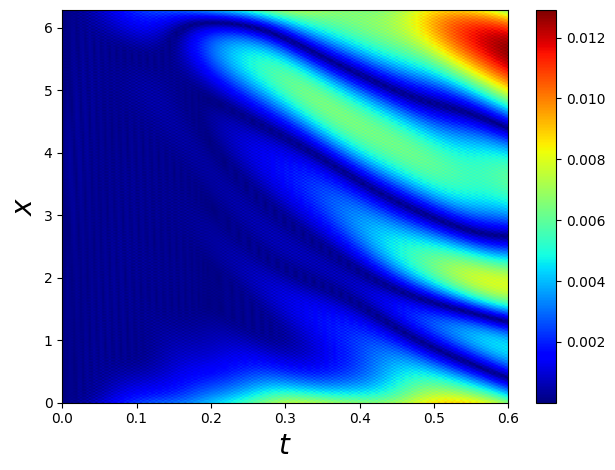}  
		%\caption*{EEMS-PINN}   
	\end{subfigure}%
    
	% 第二行
	\begin{subfigure}{.32\textwidth}  
		\centering  
		\includegraphics[width=1.00\linewidth]{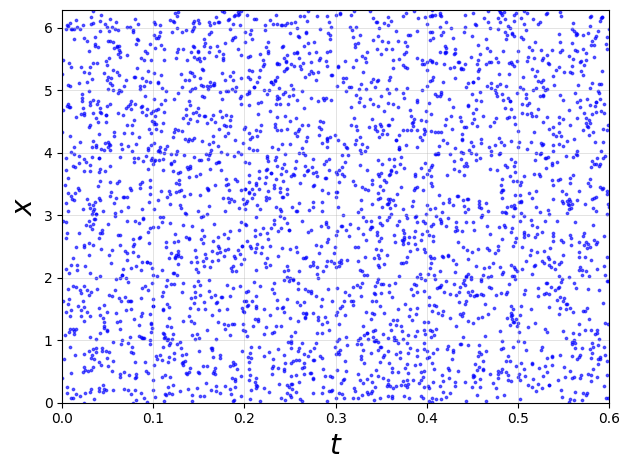}  
		\caption*{PINN}  
	\end{subfigure}%  
        \begin{subfigure}{.32\textwidth}  
		\centering  
		\includegraphics[width=1.00\linewidth]{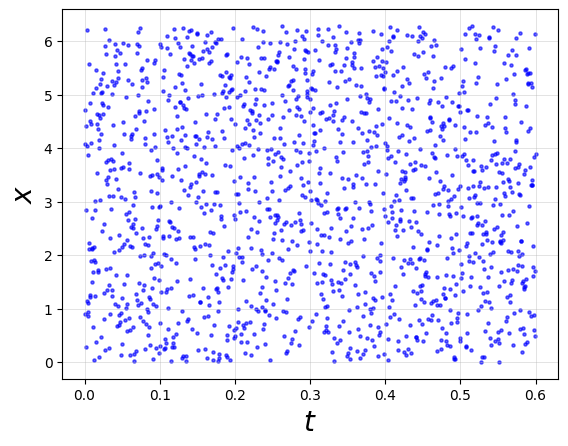}  
		\caption*{WAM-PINN}  
	\end{subfigure}%
	\begin{subfigure}{.32\textwidth}  
		\centering  
		\includegraphics[width=1.00\linewidth]{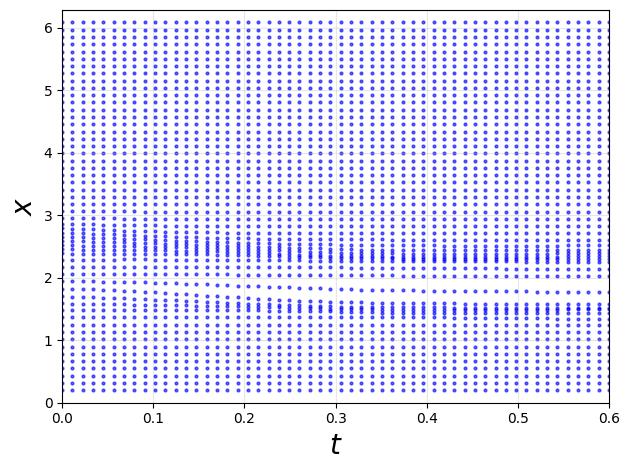}  
		\caption*{EEMS-PINN}  
	\end{subfigure}%
		% 为整个图像阵列添加标题和标签  
	\caption{(first row) The absolute errors of the estimated solutions; (second row) The distribution of sampling points for KdV equation (\ref{equ:kdv}) after one round mesh moving.}
	\label{fig:KdV_points} % 整个图像阵列的标签
\end{figure}

\begin{figure} [!htp]  
\centering  
\begin{subfigure}{.45\textwidth}  
		\centering  
		\includegraphics[width=1.00\linewidth]{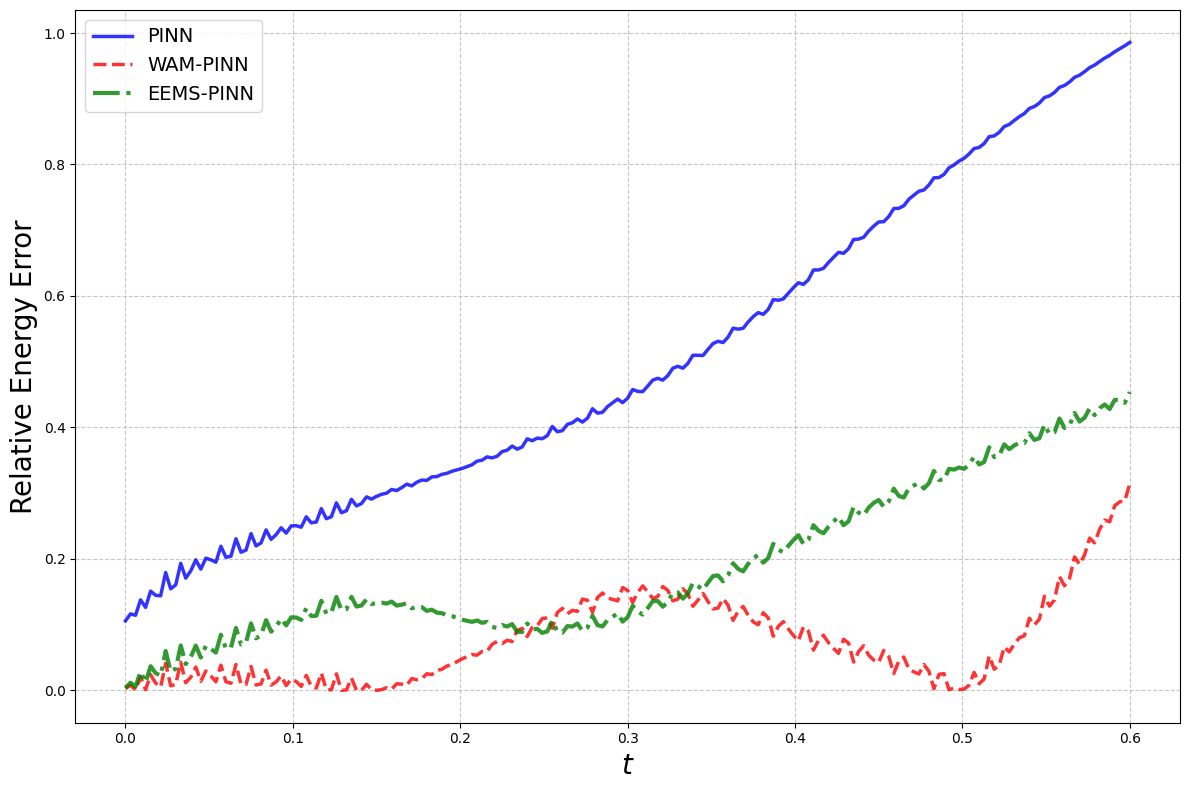}  
		\caption*{(a) Relative energy errors}  
\end{subfigure} 
 \begin{subfigure}{.45\textwidth}  
		\centering  
		\includegraphics[width=1.00\linewidth]{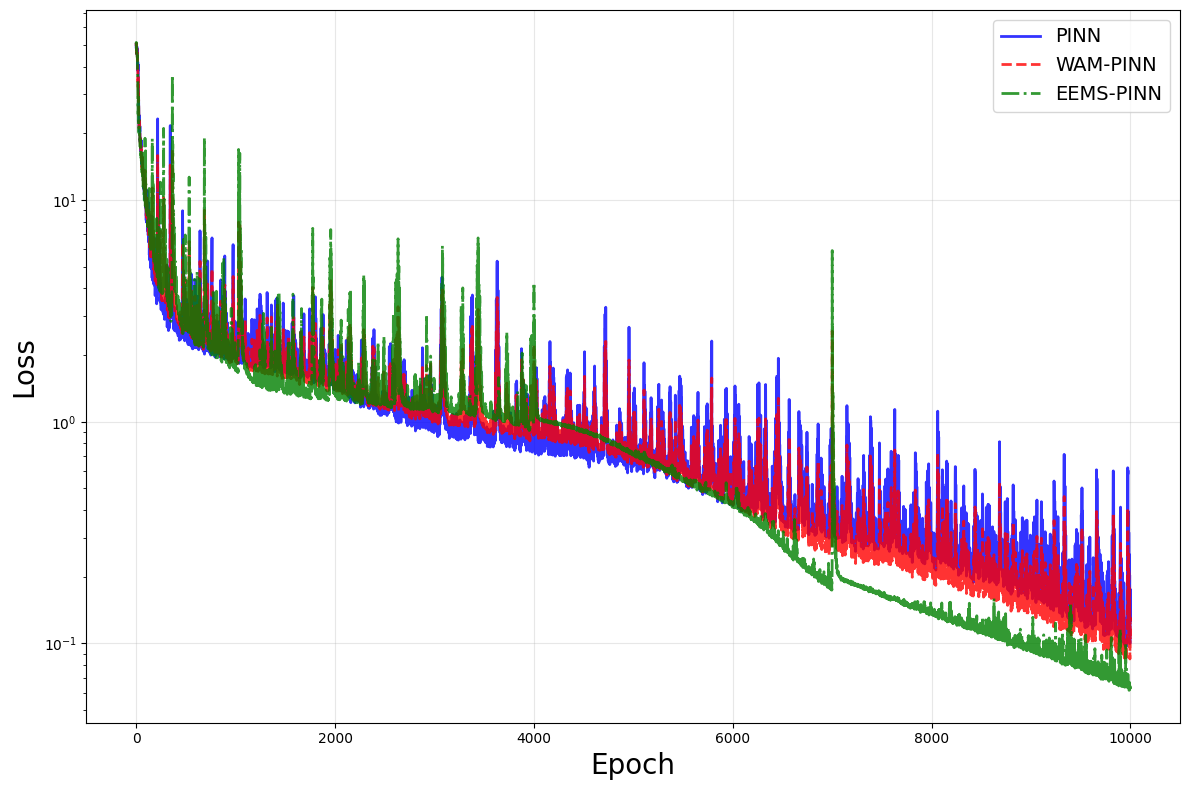}  
		\caption*{(b)  the PDE loss convergence}  
\end{subfigure}
	\caption{(a) The relative energy errors (in $\log$ scale); (b) the PDE loss convergence for KdV equation (\ref{equ:kdv}).}
    % 整个图像阵列的标题  
	\label{fig:kdv_convergence} % 整个图像阵列的标签
\end{figure}

Figures~\ref{fig:kdv_sol}--\ref{fig:kdv_convergence} present a comprehensive evaluation of solution accuracy and adaptive sampling for the KdV equation (\ref{equ:kdv}) under standardized training conditions ($N_p=3000$, $N_i=N_b=500$, final time $T=0.6$). The results demonstrate that EEMS-PINN maintains superior boundary approximation fidelity (Figure~\ref{fig:KdV_sol_time}). Table~\ref{Tab:KdV} shows EEMS-PINN achieves marginally better $L_2$ errors than other methods. While both EEMS-PINN and WAM-PINN maintain smaller energy errors than PINN, Figure~\ref{fig:kdv_convergence} reveals WWMS-PINN exhibits faster PDE loss convergence.

%In Figures~\ref{fig:kdv_sol}--\ref{fig:kdv_convergence}, we show a detailed comparison of solution accuracy and adaptive point distributions for KdV equation (\ref{equ:kdv}), evaluated under identical training conditions ($N_p=3000$, $N_i=N_b=500$) and to  the finial time $T=0.6$. From the Figure~\ref{fig:KdV_sol_time}, we can see that EEMS-PINN always shows better approximation near the boundaries. At the same time, EEMS-PINN achieves closer agreement with the exact solution compared to WAM-PINN, particularly in long-time simulations.  Figure~\ref{fig:KdV_points} demonstrates EEMS-PINN's enhanced ability to concentrate collocation points near the high energy density band regions. Table~\ref{Tab:KdV} indicates EEMS-PINN show slightly better $L_2$ errors among all the methods,  and Figure~\ref{fig:kdv_convergence} confirms that EEMS-PINN and WAM-PINN show smaller energy errors, while WWMS-PINN presents faster PDE loss convergence. 

\begin{table}[!htbp]
    \centering
\begin{tabular}{@{}lccc@{}}
\toprule
$N$ & 3000 & 10000 & 20000 \\ \midrule
PINN & $2.61\times 10^{-1}$& $8.39\times 10^{-2}$ & $6.30\times 10^{-2}$ \\
WAM-PINN & $1.25\times 10^{-1}$ & $6.03\times 10^{-2}$ & $5.41\times 10^{-2}$ \\
EEMS-PINN & $1.07\times 10^{-1}$ & $5.64\times 10^{-2}$ & $3.23\times 10^{-2} $\\
\bottomrule
\end{tabular}
 \caption{The relative $L_2$ errors of KdV equation (\ref{equ:kdv}).}
\label{Tab:KdV}
\end{table}

\subsection{Two-dimensional equations}

\noindent {\bf Example 5 (Linear wave equation)}  

We test the two-dimensional linear wave equation
\begin{equation}\label{equ:LW}
   u_{tt}=c^2(u_{x_1x_1}+u_{x_2x_2}),\ \  (x,t) \in \Omega \times [0, \infty)
\end{equation}
with $\Omega=[0,1]\times [0,1]$.
Then we have the analytical solution 
$$
u(x_1, x_2, t)= \sin(\pi x)\sin(\pi y)\cos(\sqrt{2}\pi ct) 
$$
with initial conditions
$$
u(x_1, x_2, 0)= \sin(\pi x)\sin(\pi y), \ u_t(x_1, x_2,0) = 0, 
$$
and the Dirichlet boundary conditions 
$$u(0,x_2,t) = u(1, x_2,t) = u(x_1, 0,t) = u(x_1, 1,t) = 0.$$
The energy of wave equation (\ref{equ:LW}) is given by
$$
\begin{aligned}
H(t)
& =\frac{1}{2} \iint_{} \left(u_t^2+ u_{x_1}^2 + u_{x_2}^2\right)\  {\rm d} {x_1}{\rm d} {x_2}.
\end{aligned}
$$

\begin{figure}[htbp]  
	\centering
	
	% t=0 时刻的精确解和数值解
	\begin{subfigure}{.24\textwidth}  
		\centering  
		\includegraphics[width=1.00\linewidth]{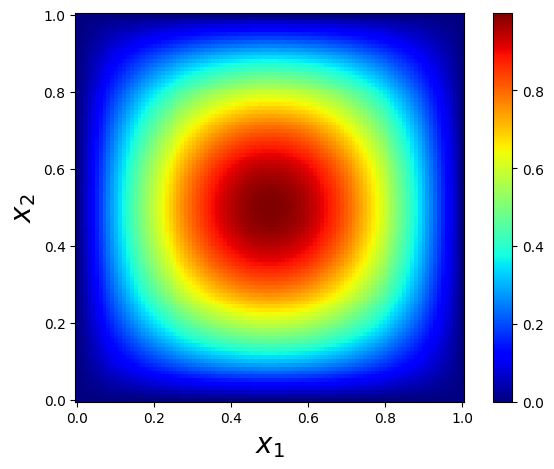}  
	\end{subfigure}%  
	\begin{subfigure}{.24\textwidth}  
		\centering  
		\includegraphics[width=1.00\linewidth]{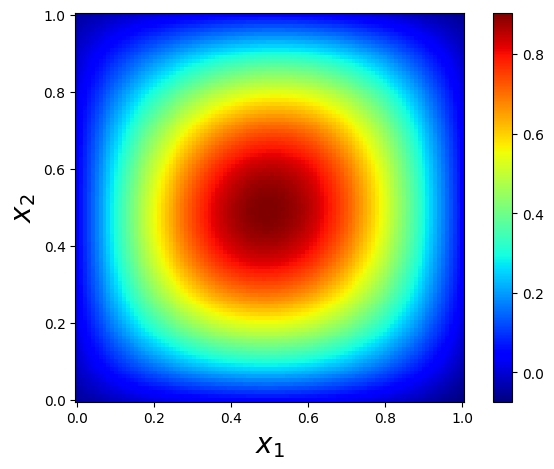}
	\end{subfigure}%
	\begin{subfigure}{.24\textwidth}  
		\centering  
		\includegraphics[width=1.00\linewidth]{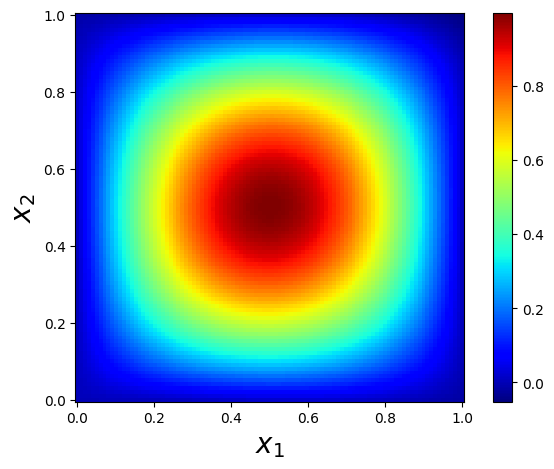} 
	\end{subfigure}%
	\begin{subfigure}{.24\textwidth}  
		\centering  
		\includegraphics[width=1.00\linewidth]{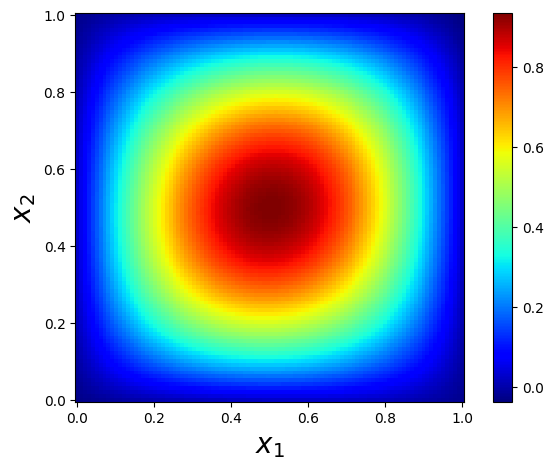}  
	\end{subfigure}
	\\
    % t=0.5 时刻的精确解和数值解
	\begin{subfigure}{.24\textwidth}  
		\centering  
		\includegraphics[width=1.00\linewidth]{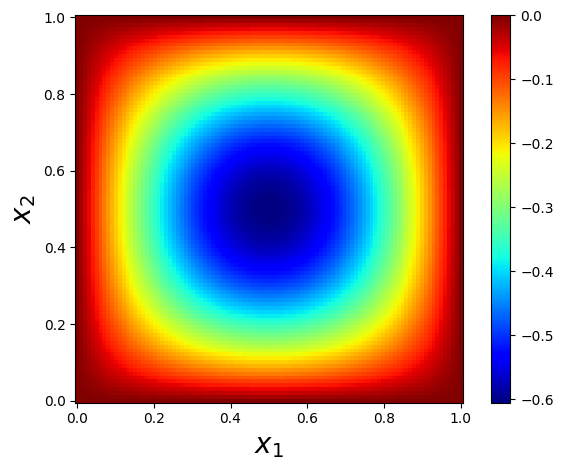}  
	\end{subfigure}%  
	\begin{subfigure}{.24\textwidth}  
		\centering  
		\includegraphics[width=1.00\linewidth]{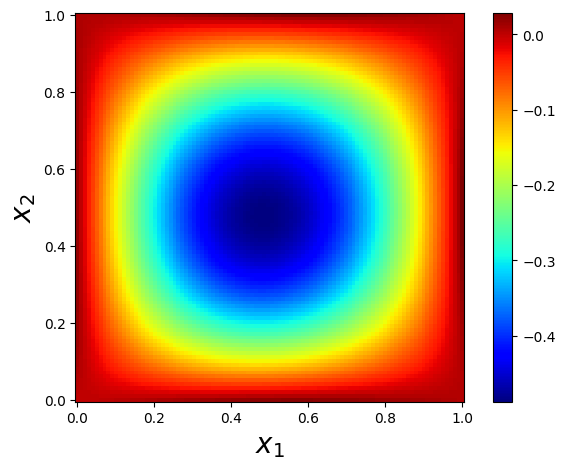}   
	\end{subfigure}%
	\begin{subfigure}{.24\textwidth}  
		\centering  
		\includegraphics[width=1.00\linewidth]{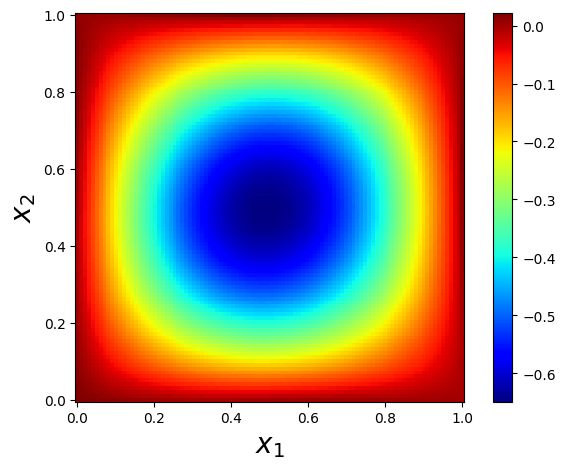} 
	\end{subfigure}%
	\begin{subfigure}{.24\textwidth}  
		\centering  
		\includegraphics[width=1.00\linewidth]{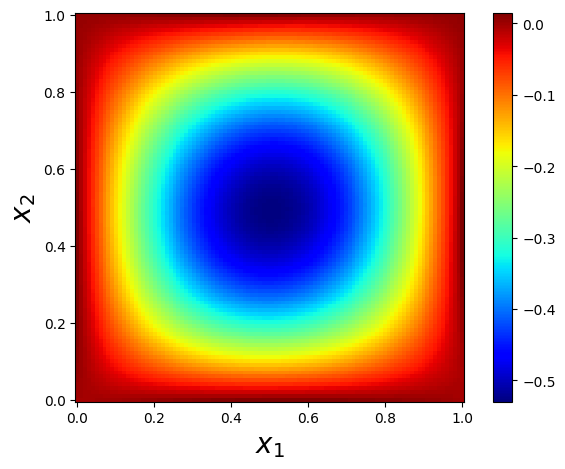}  
	\end{subfigure}
\\
	% t=1 时刻的精确解和数值解
	\begin{subfigure}{.24\textwidth}  
		\centering  
		\includegraphics[width=1.00\linewidth]{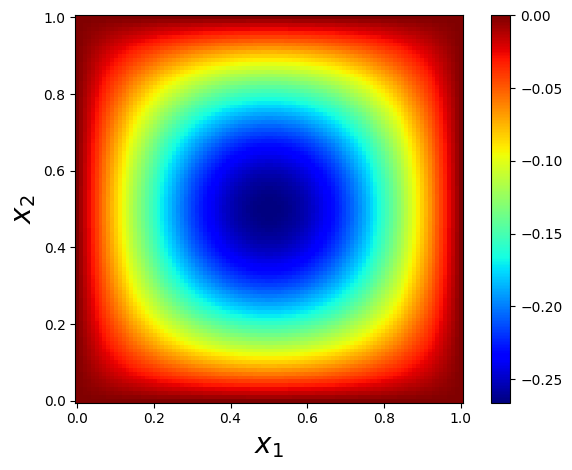}  
		\caption*{Exact}   
	\end{subfigure}%  
	\begin{subfigure}{.24\textwidth}  
		\centering  
		\includegraphics[width=1.00\linewidth]{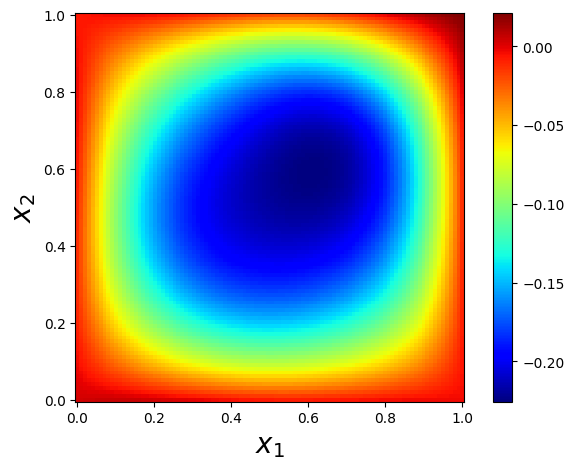} 
		\caption*{PINN}  
	\end{subfigure}%
	\begin{subfigure}{.24\textwidth}  
		\centering  
		\includegraphics[width=1.00\linewidth]{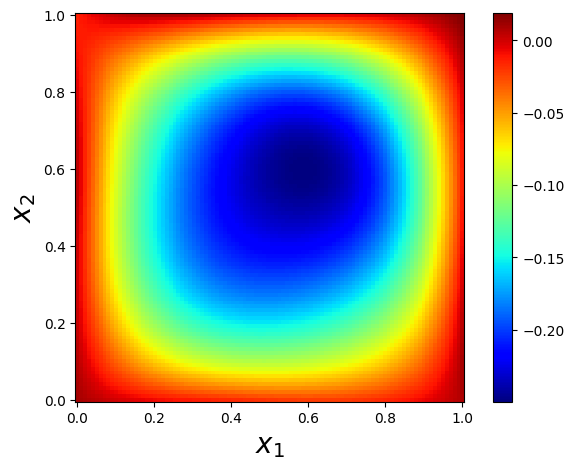} 
		\caption*{WAM-PINN}   
	\end{subfigure}%
	\begin{subfigure}{.24\textwidth}  
		\centering  
		\includegraphics[width=1.00\linewidth]{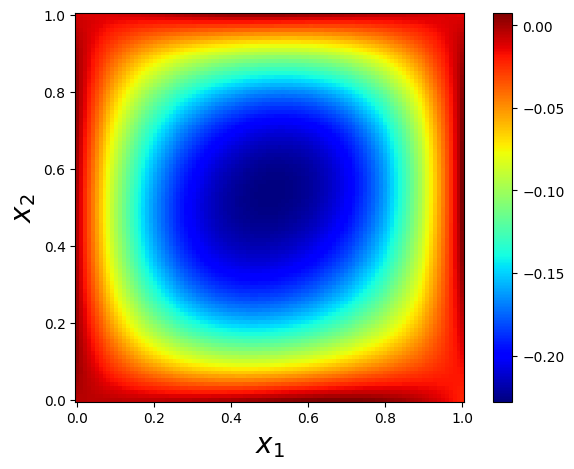}  
		\caption*{EEMS-PINN}  
	\end{subfigure}
	\caption{The exact solution and numerical solutions at $t=0$ (first row), $t=0.5$ (second row) and $t=1$ (last row) of PINN, WAM-PINN and EEMS-PINN for the linear wave equation (\ref{equ:LW}), respectively.} 
	\label{fig:Linearwave_sol}
\end{figure}

\begin{figure}[!htp]  
	\centering  
	
	% 第一行: t=0 时刻的误差
	\begin{subfigure}{.32\textwidth}  
		\centering  
		\includegraphics[width=1.00\linewidth]{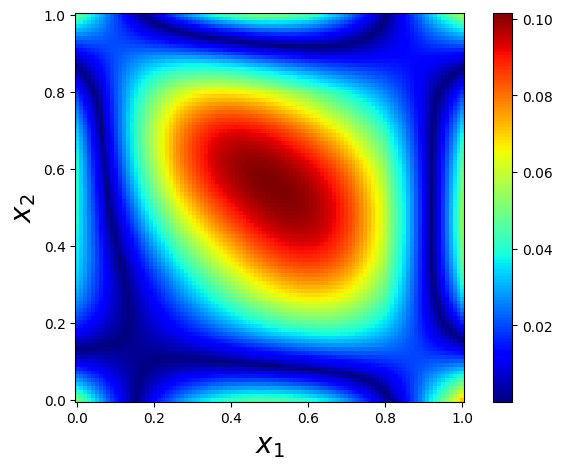}  
	\end{subfigure}%  
	\begin{subfigure}{.32\textwidth}  
		\centering  
		\includegraphics[width=1.00\linewidth]{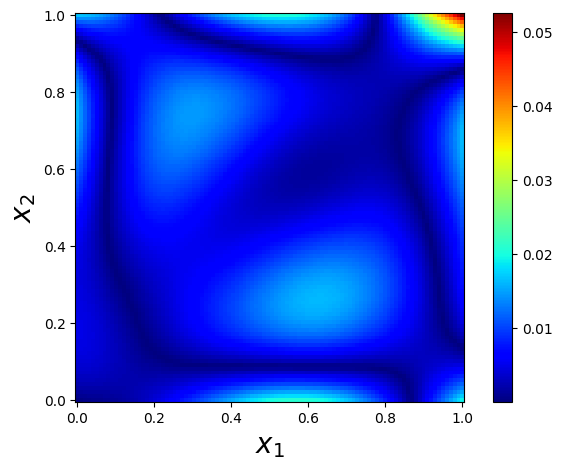}   
	\end{subfigure}%  
	\begin{subfigure}{.32\textwidth}  
		\centering  
		\includegraphics[width=1.00\linewidth]{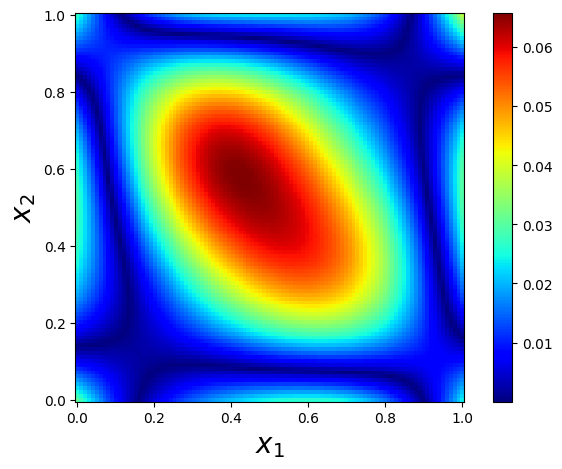}    
	\end{subfigure}%
    
	% 第二行: t=0.5 时刻的误差
	\begin{subfigure}{.32\textwidth}  
		\centering  
		\includegraphics[width=1.00\linewidth]{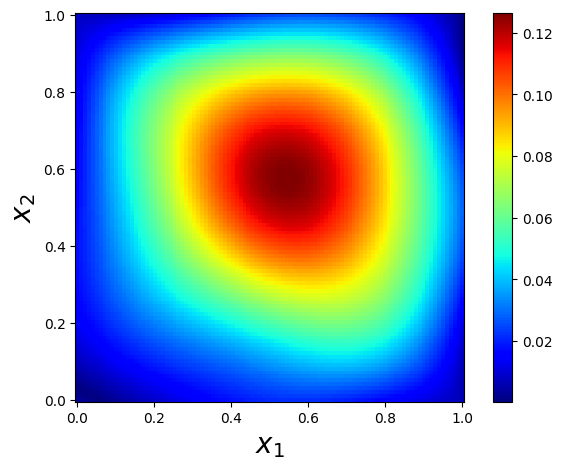}  
	\end{subfigure}%  
	\begin{subfigure}{.32\textwidth}  
		\centering  
		\includegraphics[width=1.00\linewidth]{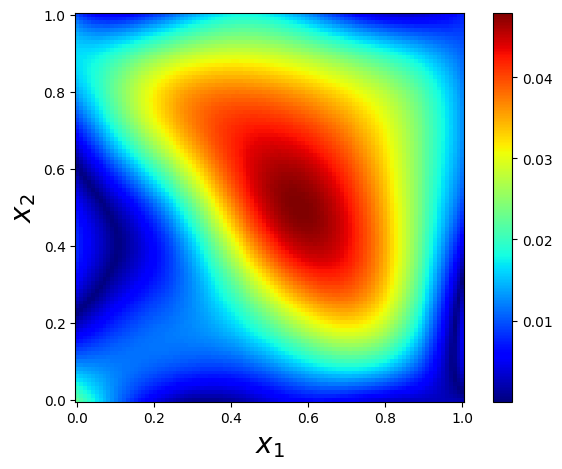}   
	\end{subfigure}%
	\begin{subfigure}{.32\textwidth}  
		\centering  
		\includegraphics[width=1.00\linewidth]{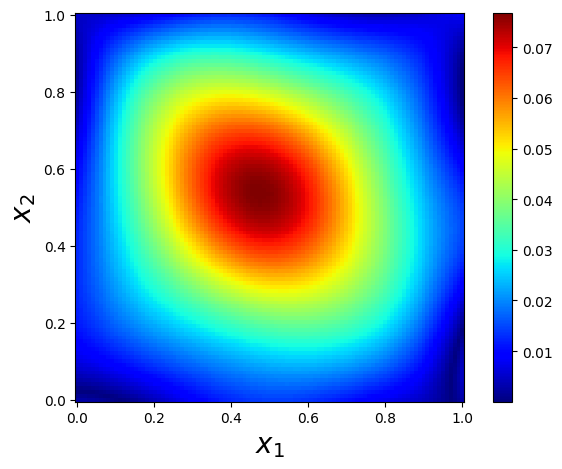}  
	\end{subfigure}%
	
	% 第三行: t=1 时刻的误差  
	\begin{subfigure}{.32\textwidth}  
		\centering  
		\includegraphics[width=1.00\linewidth]{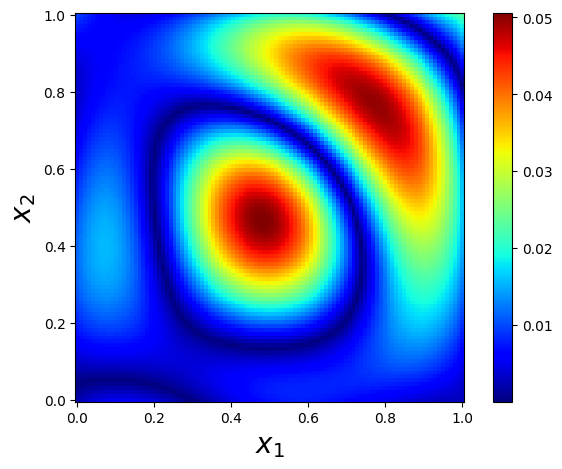}  
        \caption*{PINN}  
	\end{subfigure}%  
	\begin{subfigure}{.32\textwidth}  
		\centering  
		\includegraphics[width=1.00\linewidth]{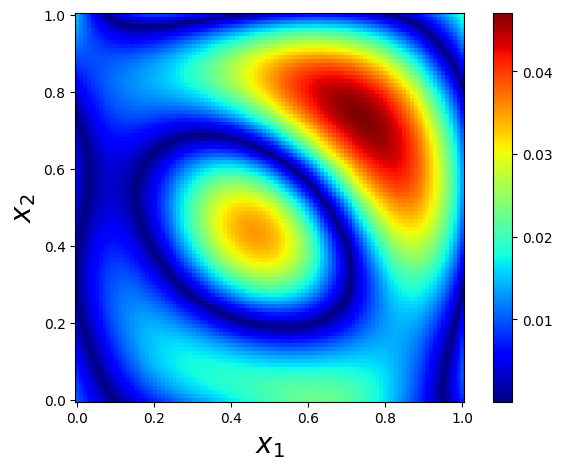}
        \caption*{WAM-PINN}  
	\end{subfigure}%
	\begin{subfigure}{.32\textwidth}  
		\centering  
		\includegraphics[width=1.00\linewidth]{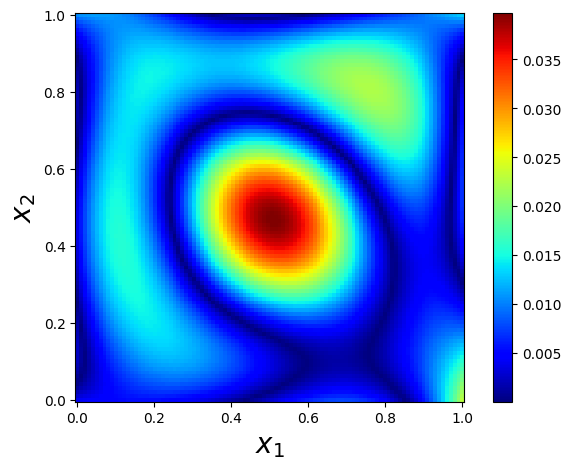}  
        \caption*{EEMS-PINN}
	\end{subfigure}
	\caption{Absolute errors of estimated solutions at  at $t=0$ (first row), $t=0.5$ (second row) and $t=1$ (last row) for the liner wave equation (\ref{equ:LW}) computed by PINN, WAM-PINN and EEMS-PINN, respectively.}
	\label{fig:linearwave_errors}
\end{figure}

In this experiment, we set $c=1$. 
For the linear wave equation (\ref{equ:LW}), the energy density function coincides with the gradient-based monitor function employed in WAM-PINN, leading to comparable qualitative behavior between WAM-PINN and EEMS-PINN as shown in Figures~\ref{fig:Linearwave_sol}--\ref{fig:linearwave_errors} when using $N_p=5000$ collocation points with $N_i=500$, $N_b=1000$. However, quantitative analysis in Table~\ref{Tab:Lwave} reveals EEMS-PINN's progressive advantage with increasing sample density: while all methods show error reduction from $N=5000$ to $N=20000$, EEMS-PINN achieves superior relative $L_2$ errors compared to WAM-PINN and PINN. %Notably, EEMS-PINN's error at $N=20000$ represents a 19\% improvement over WAM-PINN and 62\% over conventional PINN, demonstrating its enhanced capacity to utilize additional collocation points effectively.

%It can be noticed that the energy density function of the linear wave equation is the same as the gradient-based monitor function in WAM-PINN. So from Figures~\ref{fig:Linearwave_sol}--\ref{fig:linearwave_errors}, we can see that WAM-PINN and  EEMS-PINN present the same order of absolute solution errors for the linear wave equation  using $N_p=5000$ collocation points, $N_i=$ initial points and $N_b=$ boundary points. From Table \ref{Tab:Lwave}, the relative $L_2$ errors of EEMS-PINN is better than PINN when we increase the number of the samples. 

\begin{figure}[!htp]  
	\centering  
	
	% 第一行: t=0 时刻的采样点
	\begin{subfigure}{.32\textwidth}  
		\centering  
		\includegraphics[width=1.00\linewidth]{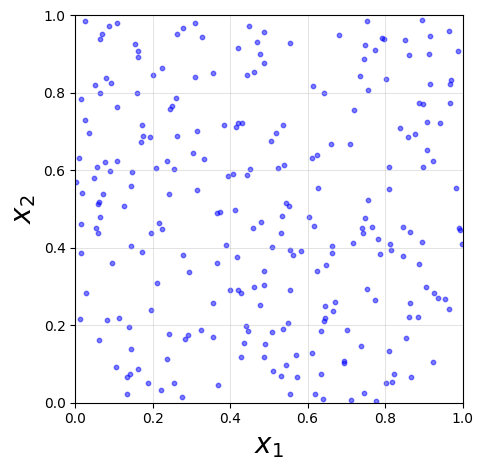}  
	\end{subfigure}%  
	\begin{subfigure}{.32\textwidth}  
		\centering  
		\includegraphics[width=1.00\linewidth]{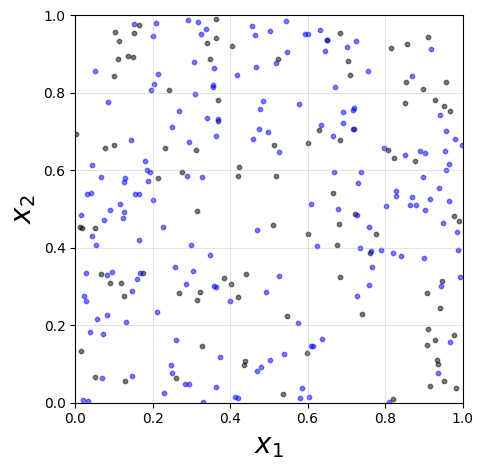}  
	\end{subfigure}%  
	\begin{subfigure}{.32\textwidth}  
		\centering  
		\includegraphics[width=1.00\linewidth]{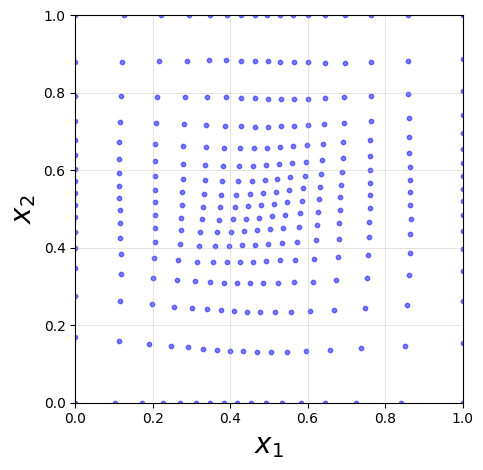}  
		%\caption*{EEMS-PINNs($t=0$)}   
	\end{subfigure}%
    
	% 第二行: t=0.5 时刻的采样点
	\begin{subfigure}{.32\textwidth}  
		\centering  
		\includegraphics[width=1.00\linewidth]{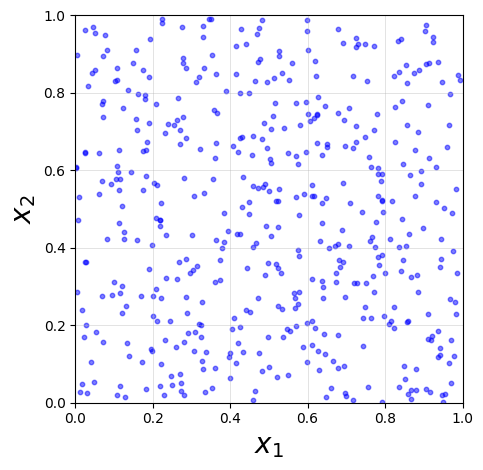}  
	%	\caption*{PINNs($t=0.5$)}  
	\end{subfigure}%  
	\begin{subfigure}{.32\textwidth}  
		\centering  
		\includegraphics[width=1.00\linewidth]{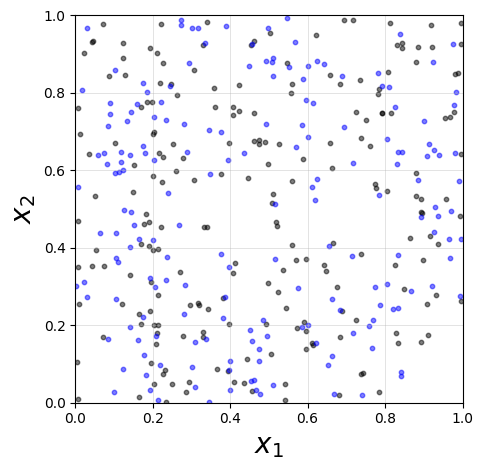}  
	%	\caption*{WAM-PINNs($t=0.5$)}  
	\end{subfigure}%
	\begin{subfigure}{.32\textwidth}  
		\centering  
		\includegraphics[width=1.00\linewidth]{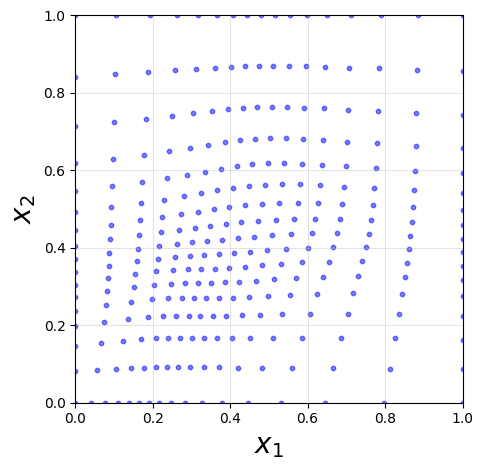}  
		%\caption*{EEMS-PINNs($t=0.5$)}  
	\end{subfigure}%
	
	% 第三行: t=1 时刻的采样点  
	\begin{subfigure}{.32\textwidth}  
		\centering  
		\includegraphics[width=1.00\linewidth]{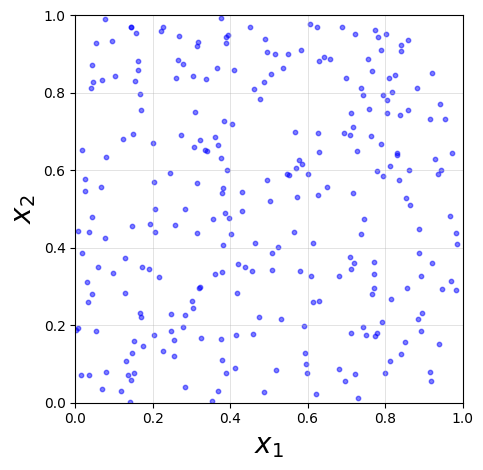}  
		\caption*{PINN}  
	\end{subfigure}%  
	\begin{subfigure}{.32\textwidth}  
		\centering  
		\includegraphics[width=1.00\linewidth]{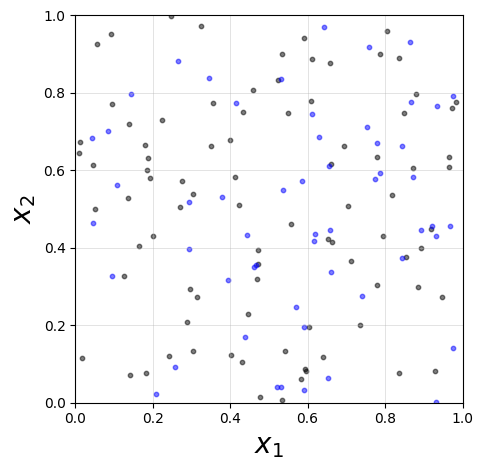}  
		\caption*{WAM-PINN}  
	\end{subfigure}%
	\begin{subfigure}{.32\textwidth}  
		\centering  
		\includegraphics[width=1.00\linewidth]{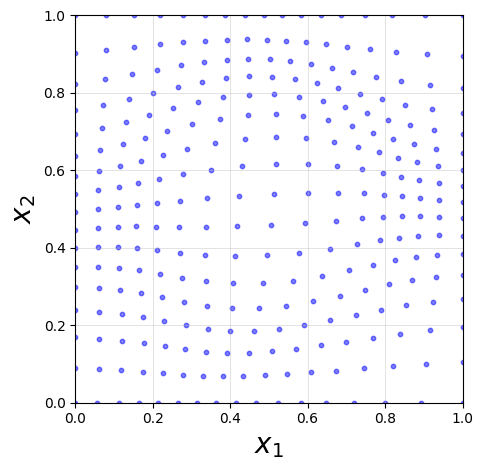}  
		\caption*{EEMS-PINN}  
	\end{subfigure}%
	
	\caption{Distribution of sampling points at  at $t=0$ (first row), $t=0.5$ (second row) and $t=1$ (last row) for the linear wave equation (\ref{equ:LW}) computed by PINN, WAM-PINN (after 10 round resampling) and EEMS-PINN (after 2 round mesh moving), respectively.}
	\label{fig:linearwave_points}
\end{figure}

Figure~\ref{fig:linearwave_points} demonstrates distinct adaptive sampling behaviors between methods for the linear wave equation. While the system's energy density exhibits relatively mild spatial variation compared to nonlinear cases, EEMS-PINN's point distributions effectively track these evolving energy patterns throughout the temporal domain ($t\in[0,1]$). In contrast, WAM-PINN's combined spatiotemporal sampling strategy produces inconsistent point distributions across time snapshots, ultimately yielding near-uniform spatial coverage. The EMMPDE framework achieves optimal mesh adaptation within just two iterations, demonstrating significantly faster convergence compared to WAM-PINN's stochastic resampling approach, which requires approximately 10 iterations to reach comparable sample optimization. 

\begin{figure} [!htp]  
\centering  
\begin{subfigure}{.45\textwidth}  
		\centering  
		\includegraphics[width=1.00\textwidth,height=0.8\textwidth]{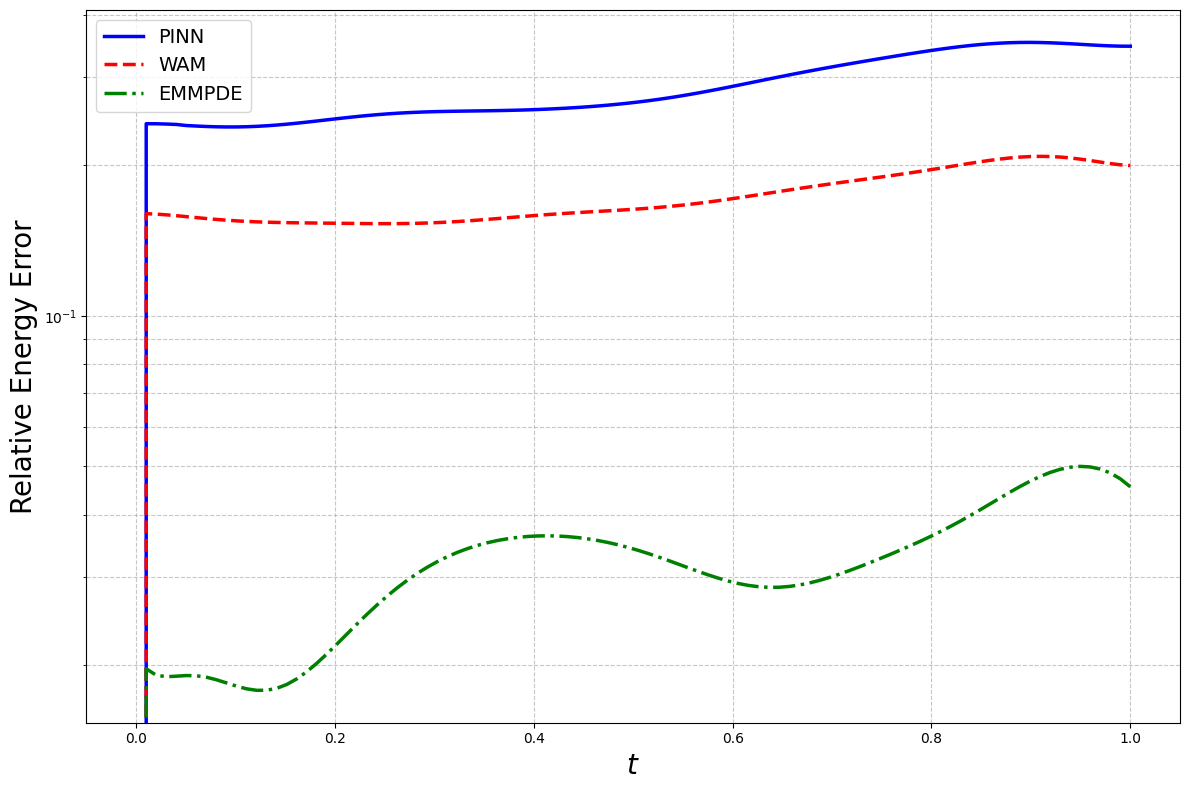}  
		\caption*{(a) Relative energy errors}  
\end{subfigure} 
 \begin{subfigure}{.45\textwidth}  
		\centering  
		\includegraphics[width=1.00\linewidth]{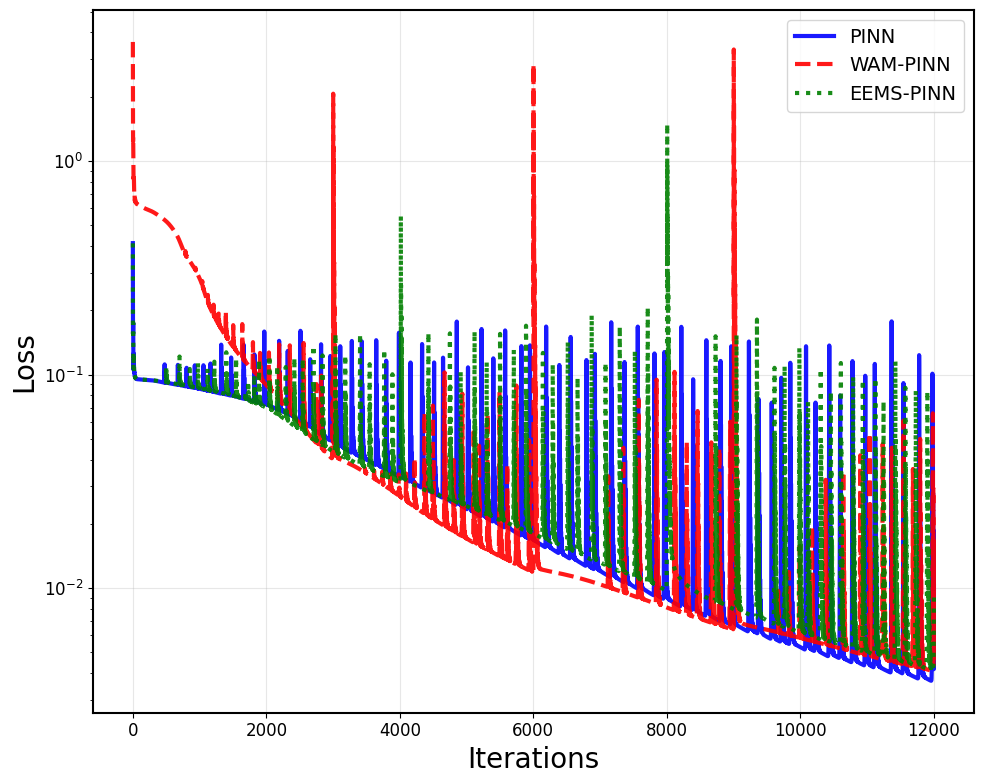}  
		\caption*{(b)  the PDE loss convergence}  
\end{subfigure}
	\caption{(a) The relative energy errors (in $\log$ scale); (b) the PDE loss convergence for linear wave equation (\ref{equ:LW}).} % 整个图像阵列的标题  
	\label{fig:LW_convergence} % 整个图像阵列的标签
\end{figure}

Figure~\ref{fig:LW_convergence} provides a broad comparison of numerical performance for the linear wave equation (\ref{equ:LW}). The results demonstrate that EEMS-PINN maintains superior energy conservation throughout the simulation---a direct consequence of the Energy-Equidistribution principle. While all methods ultimately converge to similar final PDE loss values, EEMS-PINN shows significantly faster convergence during early optimization stages. 

%Figure~\ref{fig:linearwave_points} provides a detailed comparison of the evolved collocation point distributions, revealing fundamental differences in the adaptive sampling behavior of each method. The energy density of the linear wave equation is not that steep compared to the nonlinear wave equations, but the points distribution of  EEMS-PINN can still capture the energy changes. We can see clearly how the points move during the time evolution. But in WAM-PINN, the movable points are uneven for different time snapshots  because the spatial and temporal variables are sampled together. And the movable collocation points $\mathcal{X}_{p}$ of WAM-PINN maintain a nearly uniform spatial distribution at different times. 
%EEMS-PINN demonstrates significantly more effective adaptation, with the points forming a concentrated band between $x=0$ and $x=5$.

%The relative energy errors and PDE loss convergence behavior of the results of the linear wave equation (\ref{equ:LW}) are compared across methods in Figure~\ref{fig:LW_convergence}.  Panel (a) demonstrates that EEMS-PINN maintains the best energy conservation due to the our Energy-Equidistion principle, outperforming PINN and WAM-PINN. Panel (b) reveals the training dynamics, where EEMS-PINN achieves faster convergence rate at the inital training stage but all the methods converges to basically same value in the end.

% t=0 时刻的解对比

% 误差对比图（所有时刻）

% 采样点分布对比图（所有时刻）

\begin{table}[!htbp]
    \centering
\begin{tabular}{@{}lccc@{}}
\toprule
$N$ & 5000 & 10000 & 20000 \\ \midrule
PINN & $4.72\times 10^{-2}$& $2.33\times 10^{-2}$ & $1.79\times 10^{-2}$ \\
WAM-PINN & $5.21\times 10^{-2}$ & $2.87\times 10^{-2}$ & $8.33\times 10^{-3}$ \\
EEMS-PINN & $3.27\times 10^{-2}$ & $9.60\times 10^{-3}$ & $6.87\times 10^{-3} $\\
\bottomrule
\end{tabular}
 \caption{The relative $L_2$ errors of linear wave
equation (\ref{equ:LW}).}
\label{Tab:Lwave}
\end{table}

\hspace*{\fill}

\noindent {\bf Example 6 (2D Sine-Gordon equation)}  

We now consider the two-dimensional Sine-Gordon equation 
\begin{align}\label{eq:NSineGordon_2d}
\begin{split}
    & u_{tt} - u_{x_1x_1} - u_{x_2x_2}  + \sin(u) = 0, \ (\boldsymbol{x},t) \in [-7,7]^2 \times [0,+\infty)
\end{split}   
\end{align}
with initial conditions
\begin{align*}
\begin{array}{ll}
&u(x_1, x_2, 0)=4 \arctan (\exp (x_1+x_2)), \  -7 \leqslant x_1, x_2 \leqslant 7, \\
&u_t(x_1, x_2, 0)=-\frac{4 \exp (x_1+x_2)}{1+\exp (2 x_1+2 x_2)}, \  -7 \leqslant x_1, x_2\leqslant 7,
\end{array}
\end{align*}
and the boundary conditions
\begin{align*}
&u_{x_1}=\frac{4 \exp (x_1+x_2+t)}{\exp (2 t)+\exp (2 x_1+2 x_2)},\ x_1=-7 \ \text {and} \ x_1=7, \ -7 \leqslant x_2 \leq 7, \ t>0, \\
& u_{x_2}=\frac{4 \exp (x_1+x_2+t)}{\exp (2 t)+\exp (2 x_1+2 x_2)}, \  x_2=-7 \ \text { and } \ x_2=7,\ -7 \leqslant x_1 \leqslant 7, \ t>0.
\end{align*}
The analytical solution is
$$
u(x_1, x_2, t)=4 \arctan (\exp (x_1+x_2-t)).
$$
The energy is given by
$$
H(t)=\frac{1}{2} \iint_{\Omega}\left[u_t^2+u_{x_1}^2+u_{x_2}^2+2(1-\cos u)\right] \mathrm{d} x_1\mathrm{d} x_2.
$$

\begin{figure}[!htbp]  
	\centering
	
	% t=0 时刻的精确解和数值解
	\begin{subfigure}{.24\textwidth}  
		\centering  
		\includegraphics[width=1.00\linewidth]{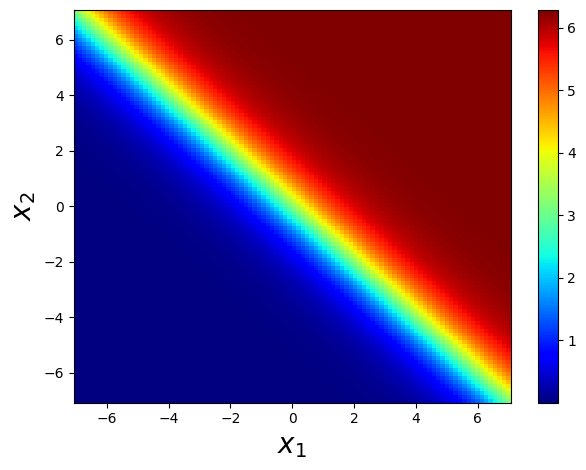}  
	\end{subfigure}%  
	\begin{subfigure}{.24\textwidth}  
		\centering  
		\includegraphics[width=1.00\linewidth]{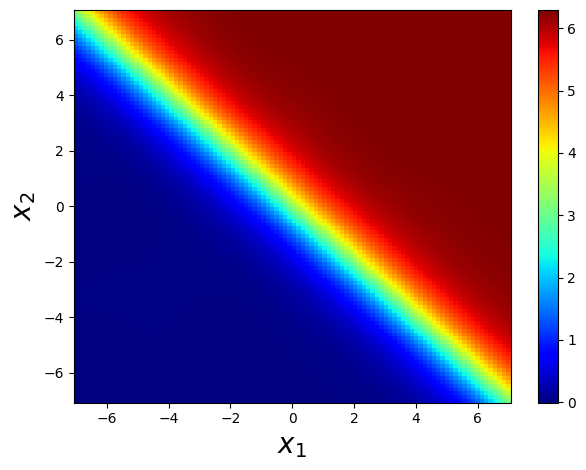}  
	\end{subfigure}%
	\begin{subfigure}{.24\textwidth}  
		\centering  
		\includegraphics[width=1.00\linewidth]{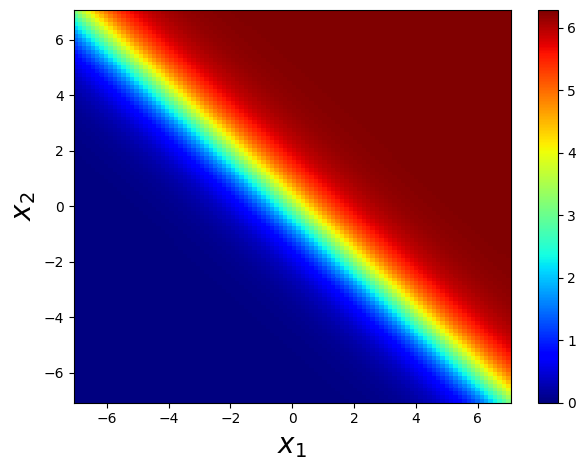}  
	\end{subfigure}%
	\begin{subfigure}{.24\textwidth}  
		\centering  
		\includegraphics[width=1.00\linewidth]{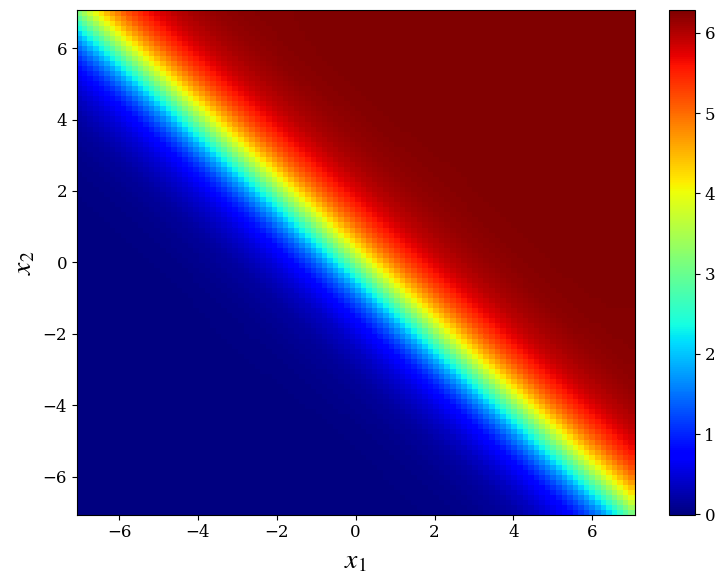}  
	\end{subfigure}
	\\
    
	% t=5 时刻的精确解和数值解
	\begin{subfigure}{.24\textwidth}  
		\centering  
		\includegraphics[width=1.00\linewidth]{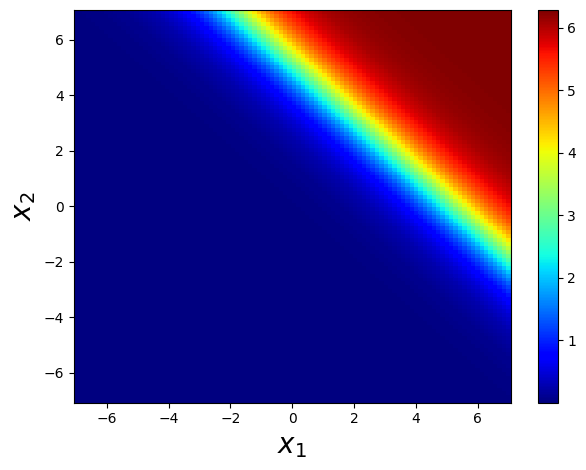}    
	\end{subfigure}%  
    \begin{subfigure}{.24\textwidth}  
		\centering  
		\includegraphics[width=1.00\linewidth]{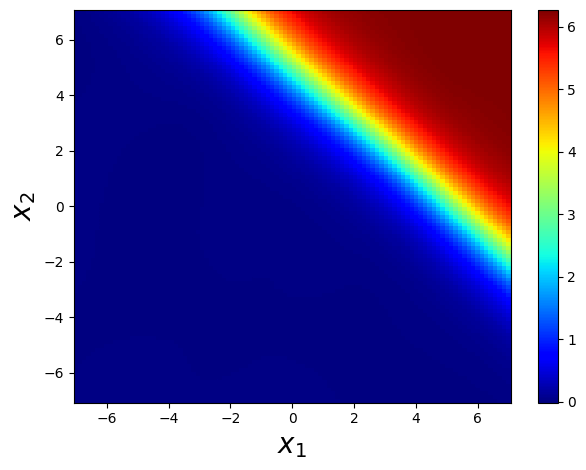}  
	\end{subfigure}%
	\begin{subfigure}{.24\textwidth}  
		\centering  
		\includegraphics[width=1.00\linewidth]{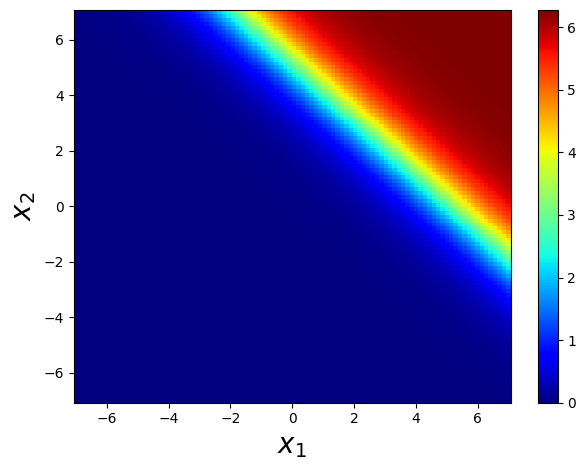}   
	\end{subfigure}%
	\begin{subfigure}{.24\textwidth}  
		\centering  
		\includegraphics[width=1.00\linewidth]{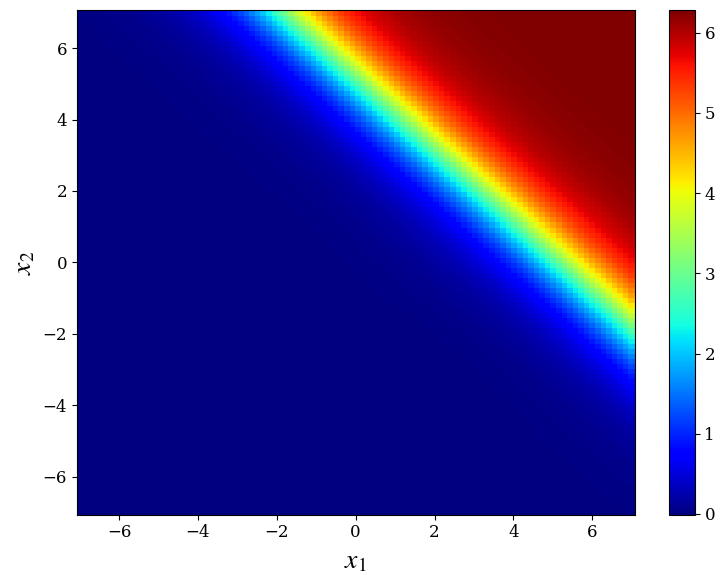}  
	\end{subfigure}
    \\
    
	% t=10 时刻的精确解和数值解
	\begin{subfigure}{.24\textwidth}  
		\centering  
		\includegraphics[width=1.00\linewidth]{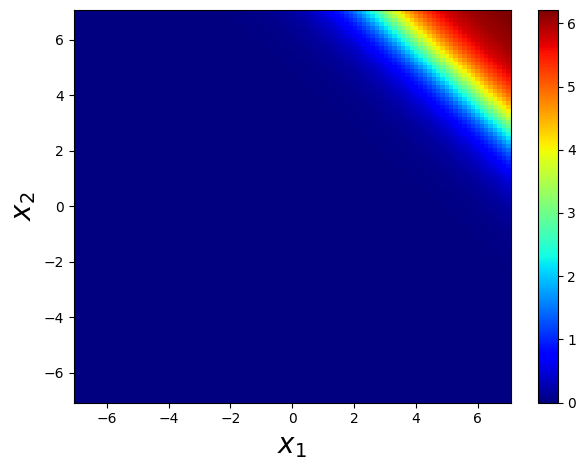}  
		\caption*{Exact}   
	\end{subfigure}%  
	\begin{subfigure}{.24\textwidth}  
		\centering  
		\includegraphics[width=1.00\linewidth]{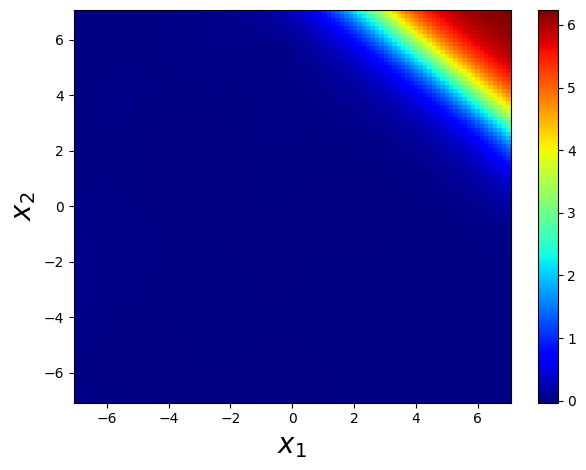}  
		\caption*{PINN}  
	\end{subfigure}%
	\begin{subfigure}{.24\textwidth}  
		\centering  
		\includegraphics[width=1.00\linewidth]{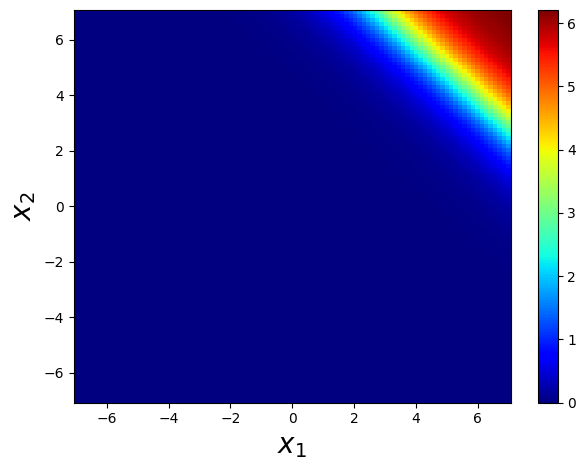} 
		\caption*{WAM-PINN}   
	\end{subfigure}%
	\begin{subfigure}{.24\textwidth}  
		\centering  
		\includegraphics[width=1.00\linewidth]{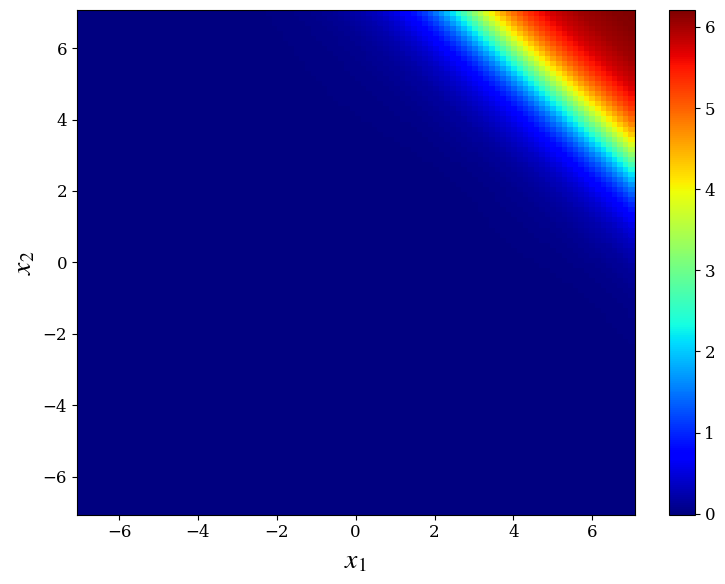}  
		\caption*{EEMS-PINN}  
	\end{subfigure}
	\caption{The exact solution and numerical solutions at at  at $t=0$ (first row), $t=5$ (second row) and $t=10$ (last row)  of PINN, WAM-PINN and EEMS-PINN for 2D Sine-Gordon equation (\ref{eq:NSineGordon_2d}).} 
	\label{fig:Sine2d_sol}
\end{figure}

\begin{figure}[!hbtp]  
	\centering  
	
	% 第一行: t=0 时刻的误差
	\begin{subfigure}{.32\textwidth}  
		\centering  
		\includegraphics[width=1.00\linewidth]{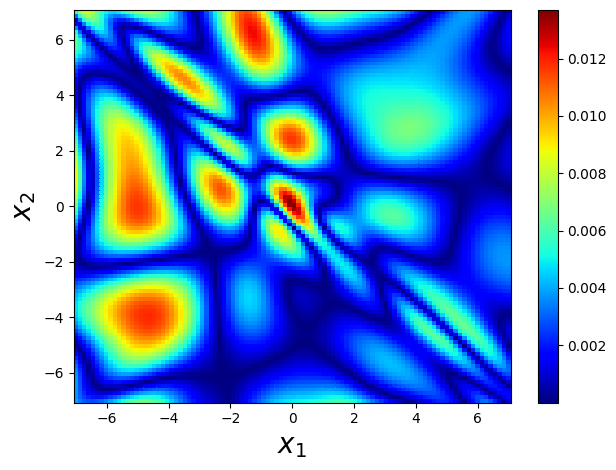}  
	\end{subfigure}%  
	\begin{subfigure}{.32\textwidth}  
		\centering  
		\includegraphics[width=1.00\linewidth]{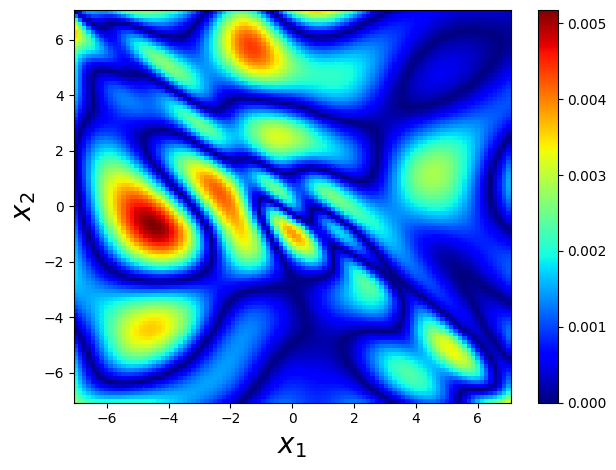}      
	\end{subfigure}%  
	\begin{subfigure}{.32\textwidth}  
		\centering  
		\includegraphics[width=1.00\linewidth]{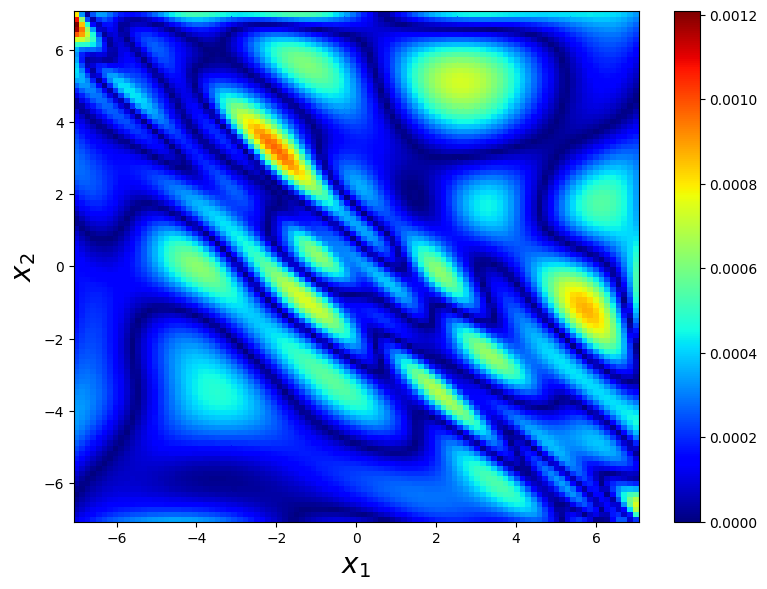}  
	\end{subfigure}%
    
	% 第二行: t=5 时刻的误差
	\begin{subfigure}{.32\textwidth}  
		\centering  
		\includegraphics[width=1.00\linewidth]{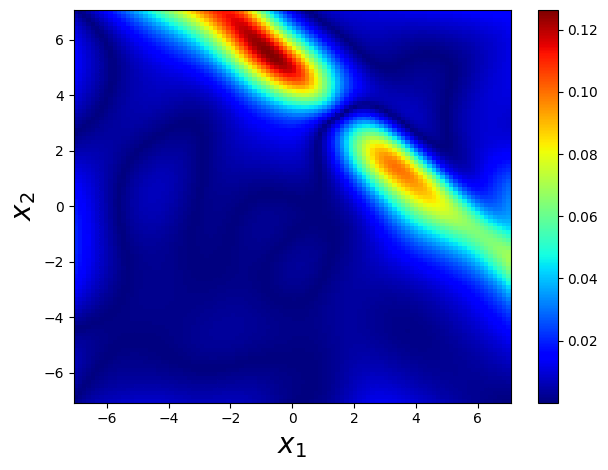}  
	\end{subfigure}%  
	\begin{subfigure}{.32\textwidth}  
		\centering  
		\includegraphics[width=1.00\linewidth]{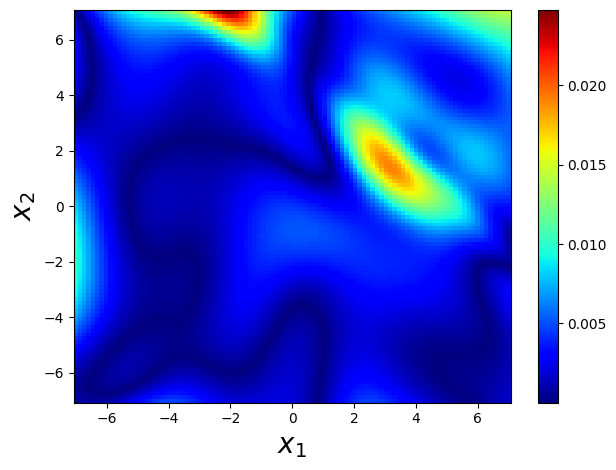}   
	\end{subfigure}%
	\begin{subfigure}{.32\textwidth}  
		\centering  
		\includegraphics[width=1.00\linewidth]{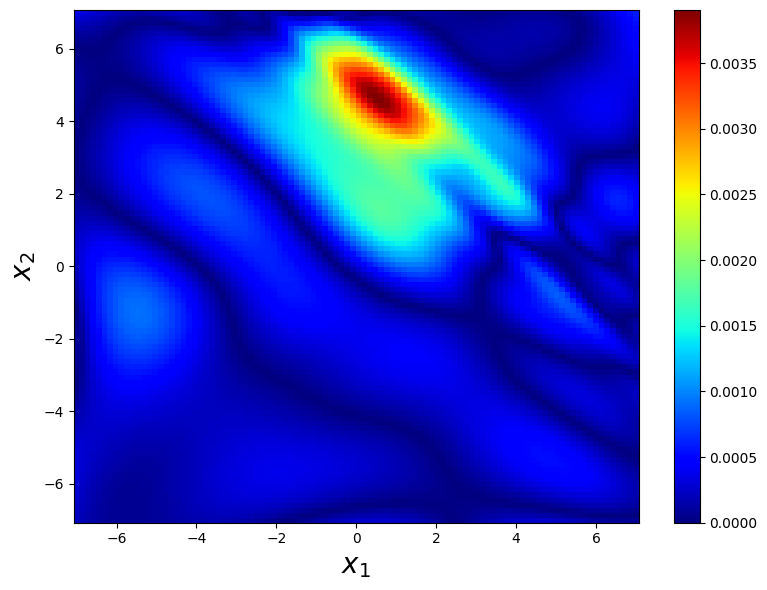}   
	\end{subfigure}%
	
	% 第三行: t=10 时刻的误差  
	\begin{subfigure}{.32\textwidth}  
		\centering  
		\includegraphics[width=1.00\linewidth]{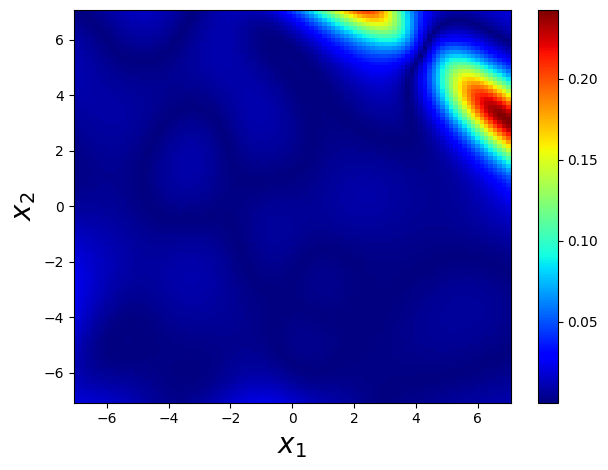}  
		\caption*{PINN}  
	\end{subfigure}%  
	\begin{subfigure}{.32\textwidth}  
		\centering  
		\includegraphics[width=1.00\linewidth]{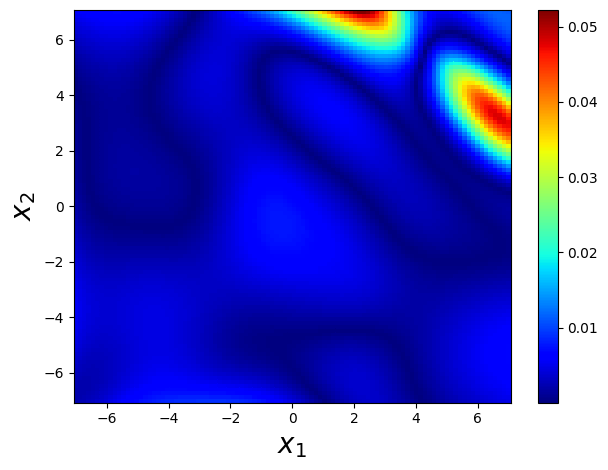}  
		\caption*{WAM-PINN}  
	\end{subfigure}%
	\begin{subfigure}{.32\textwidth}  
		\centering  
		\includegraphics[width=1.00\linewidth]{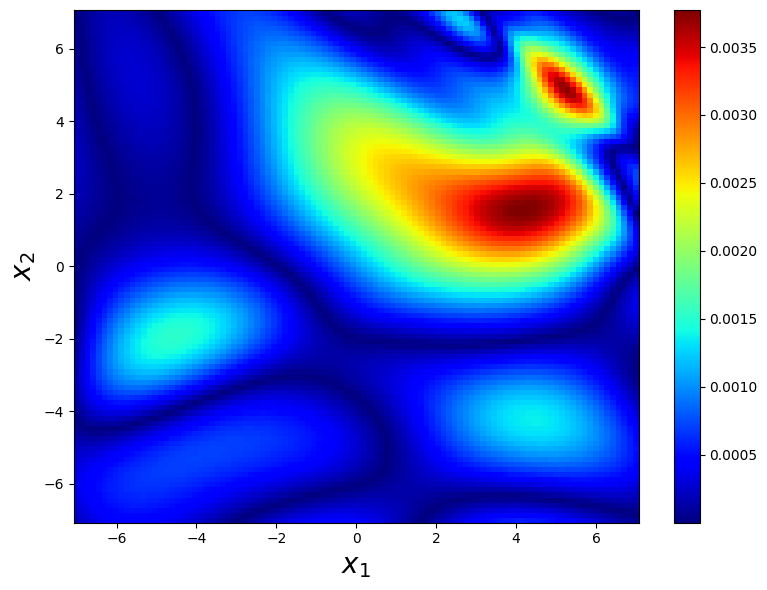}  
		\caption*{EEMS-PINN}  
	\end{subfigure}%
	
	\caption{Absolute errors of estimated solutions at at  at $t=0$ (first row), $t=5$ (second row) and $t=10$ (last row) for the 2D Sine-Gordon equation (\ref{eq:NSineGordon_2d}) computed by PINN, WAM-PINN and EEMS-PINN, respectively.}
	\label{fig:Sine2d_errors}
\end{figure}

Figures~\ref{fig:Sine2d_sol} and~\ref{fig:Sine2d_errors} present comparative solution and absolute errors of the 2D sine-Gordon equation at characteristic time snapshots ($t=0,5,10$) with $N_p=8000$ collocation points,  $N_i=400$ initial points and $N_b=400$  boundary points, demonstrating EEMS-PINN's consistent performance advantage. EEMS-PINN maintains maximum absolute errors of $\mathcal{O}(10^{-3})$ throughout the simulation, representing a tenfold improvement over WAM-PINN and a hundredfold enhancement compared to PINN. 

\begin{figure}[!htp]  
	\centering  
	
	% 第一行: t=0 时刻的采样点
	\begin{subfigure}{.32\textwidth}  
		\centering  
		\includegraphics[width=1.00\linewidth]{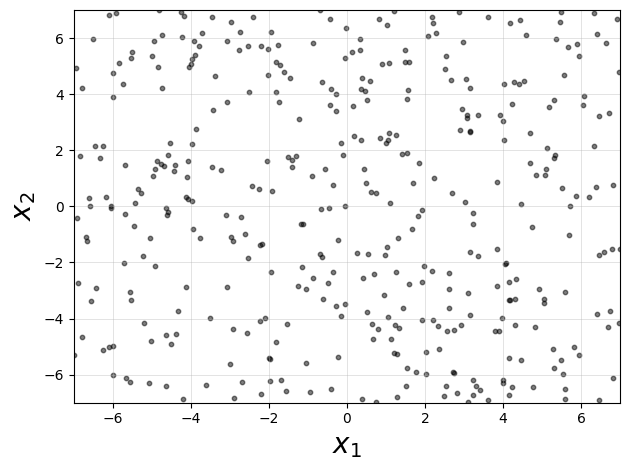}  
	\end{subfigure}%  
	\begin{subfigure}{.32\textwidth}  
		\centering  
		\includegraphics[width=1.00\linewidth]{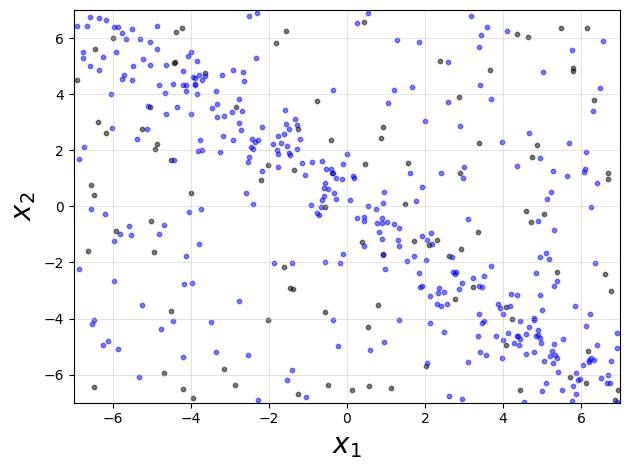}   
	\end{subfigure}%  
	\begin{subfigure}{.32\textwidth}  
		\centering  
		\includegraphics[width=1.00\linewidth]{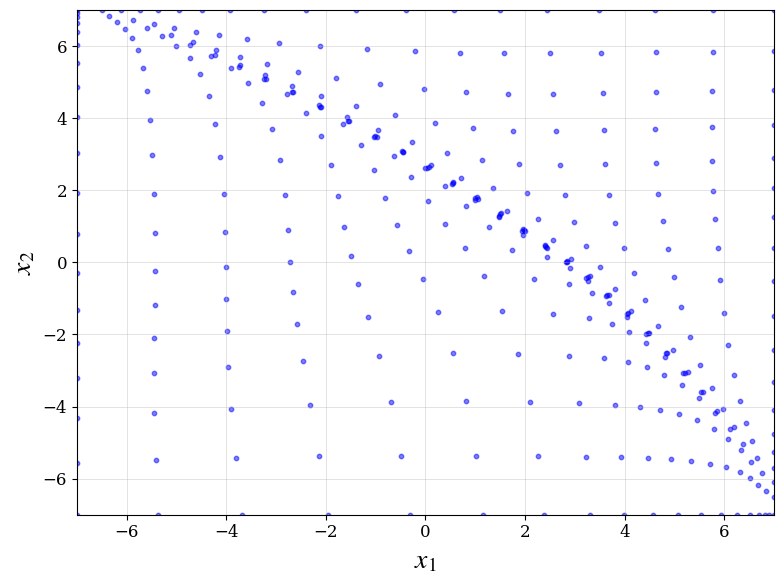}    
	\end{subfigure}%
    
	% 第二行: t=5 时刻的采样点
	\begin{subfigure}{.32\textwidth}  
		\centering  
		\includegraphics[width=1.00\linewidth]{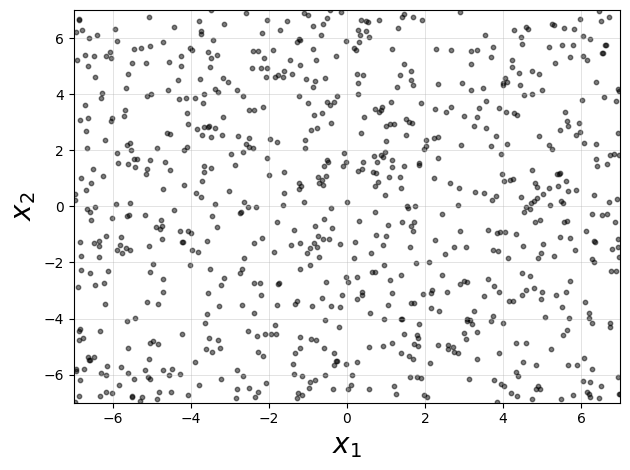}   
	\end{subfigure}%  
	\begin{subfigure}{.32\textwidth}  
		\centering  
		\includegraphics[width=1.00\linewidth]{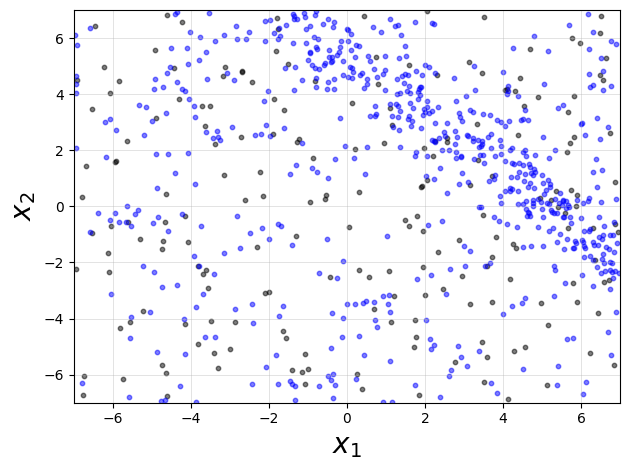}  
	\end{subfigure}%
	\begin{subfigure}{.32\textwidth}  
		\centering  
		\includegraphics[width=1.00\linewidth]{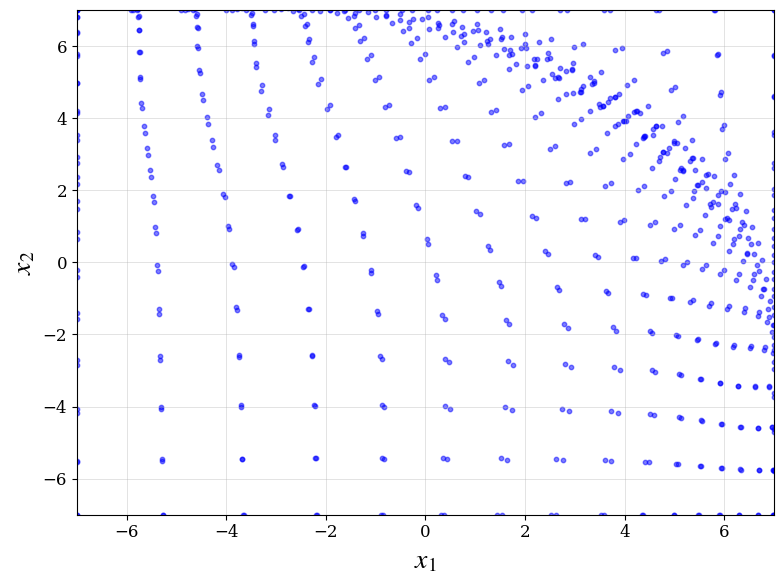}   
	\end{subfigure}%
	
	% 第三行: t=10 时刻的采样点  
	\begin{subfigure}{.32\textwidth}  
		\centering  
		\includegraphics[width=1.00\linewidth]{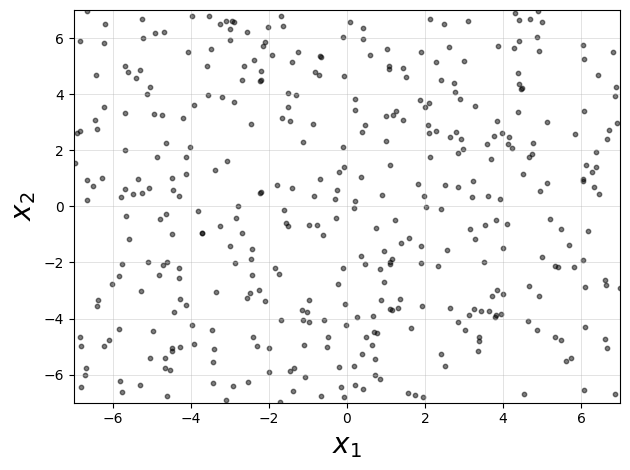}  
		\caption*{PINN}  
	\end{subfigure}%  
	\begin{subfigure}{.32\textwidth}  
		\centering  
		\includegraphics[width=1.00\linewidth]{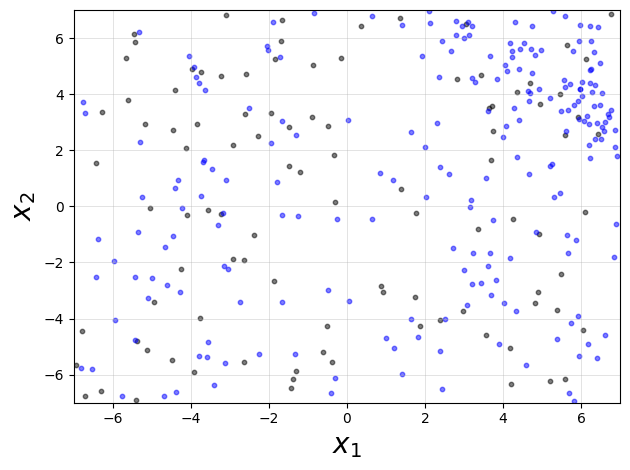}  
		\caption*{WAM-PINN}  
	\end{subfigure}%
	\begin{subfigure}{.32\textwidth}  
		\centering  
		\includegraphics[width=1.00\linewidth]{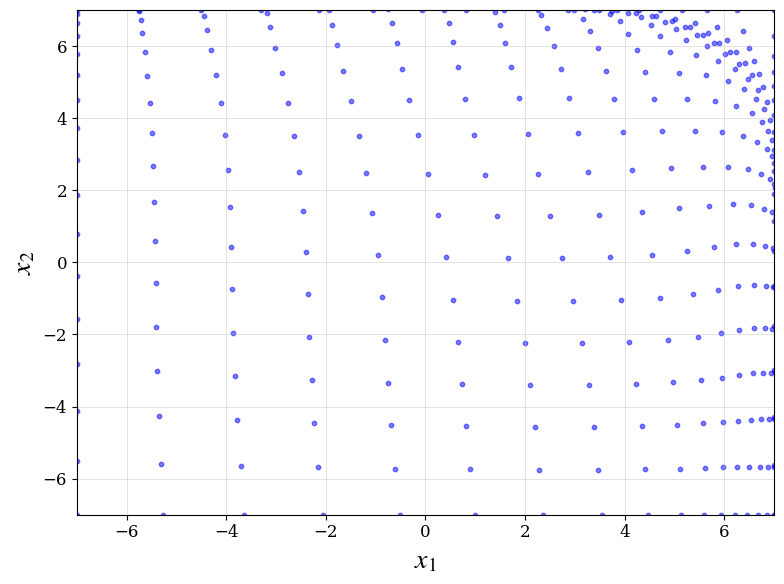}  
		\caption*{EEMS-PINN}  
	\end{subfigure}%
	
	\caption{Distribution of mobile collocation points at at  at $t=0$ (first row), $t=5$ (second row) and $t=10$ (last row) for the 2D Sine-Gordon equation (\ref{eq:NSineGordon_2d}) computed by PINN, WAM-PINN and EEMS-PINN after one round mesh moving, respectively.}
	\label{fig:Sine2d_points}
\end{figure}

Figure~\ref{fig:Sine2d_points} compares the evolved collocation point distributions, highlighting distinct adaptation behaviors across methods. While both WAM-PINN and EEMS-PINN effectively concentrate points near the propagating wavefront – corresponding to zones of extreme solution gradients and maximal energy accumulation. 

\begin{figure} [!htp]  
\centering  
\begin{subfigure}{.45\textwidth}  
		\centering  
		\includegraphics[width=1.00\textwidth,height=0.68\textwidth]{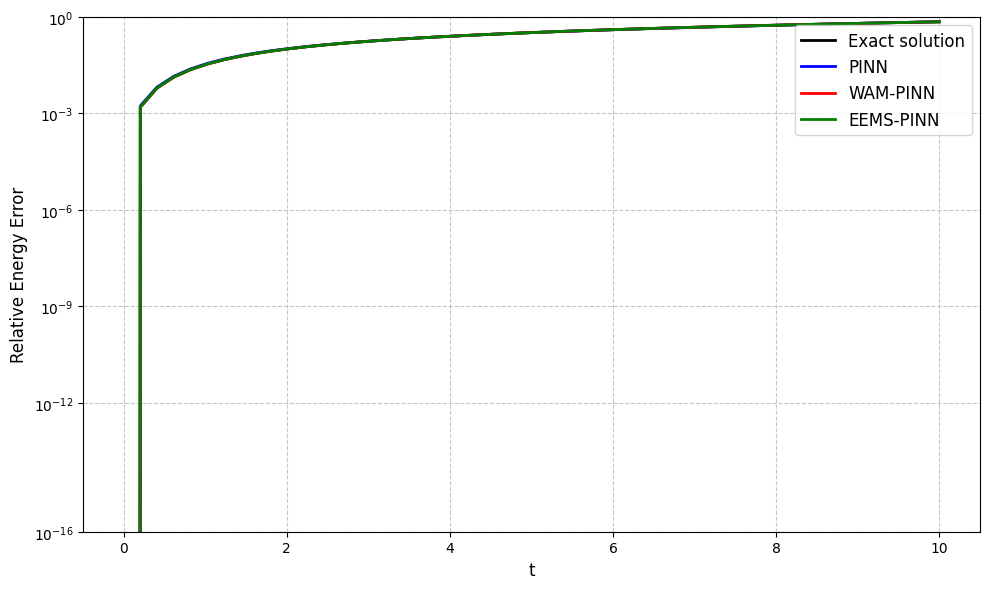}  
		\caption*{(a) Relative energy errors}  
\end{subfigure} 
 \begin{subfigure}{.45\textwidth}  
		\centering  
		\includegraphics[width=1.00\textwidth]{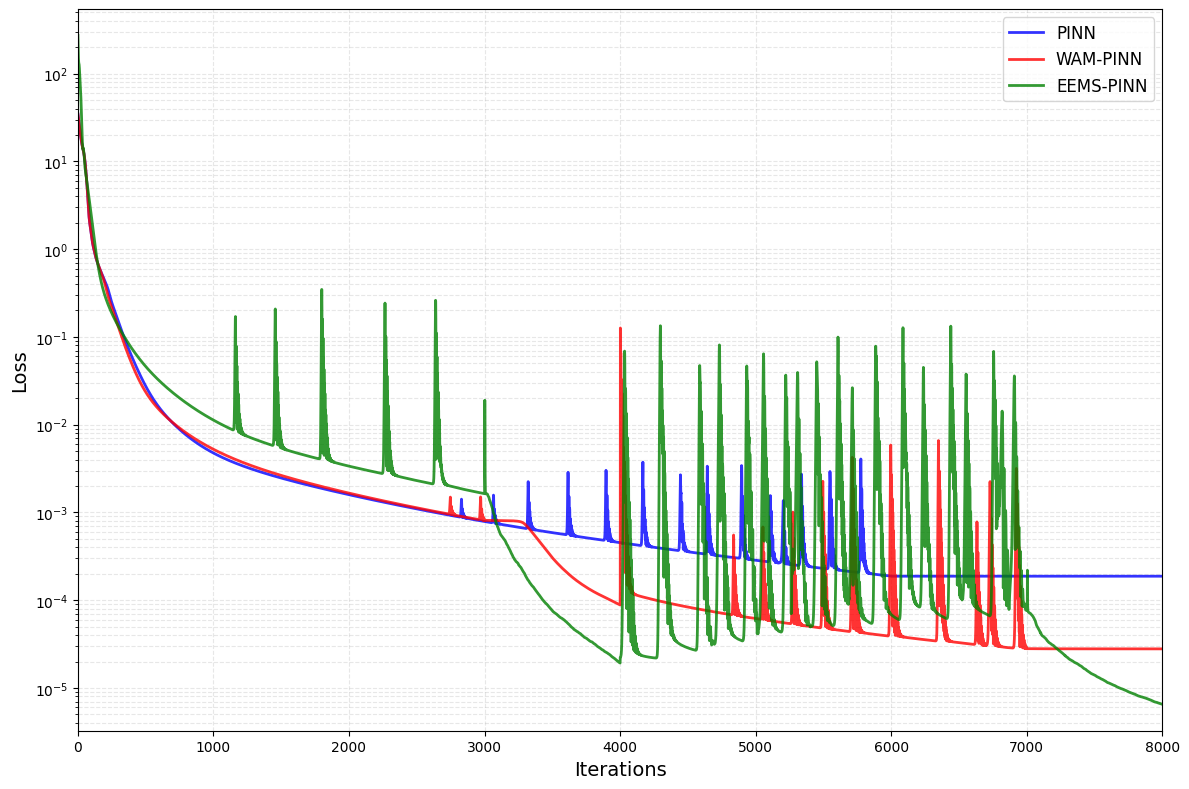}  
		\caption*{(b) the PDE loss convergence}  
\end{subfigure}
	\caption{(a) The relative energy errors (in $\log$ scale); (b) the PDE loss convergence for 2D Sine-Gordon equation (\ref{eq:NSineGordon_2d}).} % 整个图像阵列的标题  
	\label{fig:sine2d_convergence} % 整个图像阵列的标签
\end{figure}

The comparative analysis in Figure~\ref{fig:sine2d_convergence} demonstrates that while both EEMS-PINN and WAM-PINN achieve similar relative energy errors  throughout the simulation, EEMS-PINN exhibits significantly faster PDE loss convergence  during final optimization stages. Numerical results in Table~\ref{Tab:sine2d} further quantify EEMS-PINN's highly accurate approximations for the 2D Sine-Gordon system (\ref{eq:NSineGordon_2d}).

% 误差对比图（所有时刻）

% 采样点分布对比图（所有时刻）

\begin{table}[!htbp]
    \centering
    \begin{tabular}{@{}lccc@{}}
        \toprule
        $N$ & 1000 & 8000 & 27000 \\ 
        \midrule
        PINN& $2.13\times 10^{-2}$ & $9.72\times 10^{-3} $ & $5.59\times10^{-3} $\\        
        WAM-PINN& $3.83\times 10^{-3}$ & $2.97\times 10^{-3}$ & $8.58\times 10^{-4}$ \\
        EEMS-PINN& $8.72\times 10^{-4}$ & $2.90\times 10^{-4}$ & $8.13\times 10^{-5} $\\
        \bottomrule
    \end{tabular}
    \caption{The relative $L_2$ errors of 2D Sine-Gordon equation (\ref{eq:NSineGordon_2d}).}
    \label{Tab:sine2d}
\end{table}

\section{Conclusions}

%This paper presents a novel Energy-Equidistributed moving mesh strategy for solving conservative PDEs. By emphasizing the critical role of energy conservation in numerical simulations, we proposed the Energy-Equidistribution Principle (EEP) as a foundational framework to guide adaptive mesh refinement.  The energy-based moving mesh PDEs (EMMPDEs) were derived, each tailored to  dynamically move the collocation points in regions of significant energy variation. The integration of these strategies with PINNs resulted in the EEMS-PINNs framework.

%Numerical experiments across energy-conservative and -non-conservative PDEs demonstrated the superiority of the proposed method over approaches resampling from density function relying on traditional gradient-based monitor function. The energy-based monitor function ensures strict adherence to the Hamiltonian structure, reducing relative energy errors by orders of magnitude in long-term simulations. The EEMS-PINNs exhibited superior stability when handling rapidly propagating wavefronts (e.g., kink-antikink solutions with near-light-speed propagation).   

In this paper, we have introduced a novel Energy-Equidistributed moving mesh strategy for conservative PDEs, grounded in the EEP to rigorously preserve energy during adaptive mesh refinement. We derived the EMMPDEs that dynamically concentrate collocation points in regions of high energy variation, and integrate this framework with PINNs to develop EEMS-PINNs. Numerical experiments demonstrate clear advantages of the proposed method over gradient-based adaptive sampling approaches. The energy-driven framework rigorously enforces conservation laws while significantly improving solution accuracy. Performance remains robust across both conservative and non-conservative systems, maintaining stability throughout long-time simulations. 
%The numerical experiments confirm significant advantages over gradient-based adaptive sampling methods by rigorously enforcing conservation laws while significantly improving solution fidelity and demonstrating robust performance across both energy-conservative and non-conservative systems in long-term simulations. 

This research opens three key directions for further investigation. First, we will conduct rigorous theoretical analysis in EMMPDE formulation. Second, we plan to extend the Energy-Equidistributed framework to high-dimensional conservative systems, with targeted applications in image processing and turbulent fluid simulations. Finally, we will address the geometric challenges of implementing EEMS-PINNs on curved manifolds \cite{KOLASINSKI2020109097} by developing extrinsic sampling strategies that respect both the Hamiltonian dynamics and surface topology. 

%In the future research, firstly we will be doing further theoretical analysis to address the issues raised in the EMMPDE. And then we hope to apply the Energy-Equidistributed moving mesh strategy to conservative high-dimensional PDEs in a wider range of application domains, such as imaging processing and fluid dynamical simulations. At last, we will explore extending EEMS-PINNs to conservative systems  defined on curved surfaces, where developing effective sample-moving strategies for complex geometries remains a significant computational challenge \cite{KOLASINSKI2020109097}. 

%Future work will focus on extending EEMS-PINNs to stochastic Hamiltonian PDEs and conservative PDEs defined on surfaces. How to effectively moving the samples on complex geometry surfaces are still a challenging topic.  This methodology opens new avenues for simulating complex physical systems where energy conservation is paramount.

\bibliographystyle{plain}
\bibliography{main}
\end{document}